\documentclass[11pt]{amsart}
\usepackage[foot]{amsaddr}
\usepackage[letterpaper, total={6.5in,9in}, includefoot, centering]{geometry}
\usepackage{amsfonts}
\usepackage{amsmath,amsthm,amssymb,amscd}

\usepackage[mathscr]{eucal}
\usepackage{latexsym}
\usepackage{graphics}
\usepackage{verbatim}
\usepackage{upgreek}
\usepackage{overpic}
\usepackage[all,cmtip]{xy} % comm diagram
\usepackage{enumitem}
\usepackage{cancel}
\usepackage{tikz}
\usepackage{graphicx}
\graphicspath{{./images/}}
\usepackage[pagebackref]{hyperref}
\renewcommand*\backref[1]{\ifx#1\relax \else (Page #1) \fi}
\usepackage[font=scriptsize]{caption}

\usepackage[authoryear,sort&compress]{natbib}
\bibpunct{[}{]}{;}{n}{,}{,}

\makeatletter \@addtoreset{equation}{section}

\DeclareMathOperator{\ext}{ext}

\makeatother
\newtheorem{theorem}{Theorem}[section]

\newtheorem{remark}{Remark}[section]
\newtheorem{example}{Example}[section]
\newtheorem{prop}{Proposition}[section]

\DeclareMathOperator*{\sumint}{%
\mathchoice%
  {\ooalign{$\displaystyle\sum$\cr\hidewidth$\displaystyle\int$\hidewidth\cr}}
  {\ooalign{\raisebox{.14\height}{\scalebox{.7}{$\textstyle\sum$}}\cr\hidewidth$\textstyle\int$\hidewidth\cr}}
  {\ooalign{\raisebox{.2\height}{\scalebox{.6}{$\scriptstyle\sum$}}\cr$\scriptstyle\int$\cr}}
  {\ooalign{\raisebox{.2\height}{\scalebox{.6}{$\scriptstyle\sum$}}\cr$\scriptstyle\int$\cr}}
}

\begin{document}
\title{Multisymplectic Hamiltonian Variational Integrators}
\author{Brian Tran and Melvin Leok}
\address{Department of Mathematics, University of California, San Diego, 9500 Gilman Drive, La Jolla, CA 92093-0112, USA.}
\email{b3tran@ucsd.edu, mleok@ucsd.edu}
\allowdisplaybreaks

\begin{abstract}
Variational integrators have traditionally been constructed from the perspective of Lagrangian mechanics, but there have been recent efforts to adopt discrete variational approaches to the symplectic discretization of Hamiltonian mechanics using Hamiltonian variational integrators. In this paper, we will extend these results to the setting of Hamiltonian multisymplectic field theories. We demonstrate that one can use the notion of Type II generating functionals for Hamiltonian partial differential equations as the basis for systematically constructing Galerkin Hamiltonian variational integrators that automatically satisfy a discrete multisymplectic conservation law, and establish a discrete Noether's theorem for discretizations that are invariant under a Lie group action on the discrete dual jet bundle. In addition, we demonstrate that for spacetime tensor product discretizations, one can recover the multisymplectic integrators of Bridges and Reich, and show that a variational multisymplectic discretization of a Hamiltonian multisymplectic field theory using spacetime tensor product Runge--Kutta discretizations is well-defined if and only if the partitioned Runge--Kutta methods are symplectic in space and time. 
\end{abstract}

\maketitle

\tableofcontents

\section{Introduction}

Variational integrators have become an important class of geometric numerical integrators for the simulation of mechanical systems, and provides a systematic method of constructing symplectic integrators. The variational approach has numerous benefits, the first of which is that the resulting numerical integrators are automatically symplectic, and if they are group-invariant, then they satisfy a discrete Noether's theorem and preserve a discrete momentum map. In addition, it can be shown that the order of accuracy is related to the best approximation properties of the finite-dimensional function spaces and the order of the quadrature rule used to construct the variational integrator~\cite{HaLe2015}.

However, the variational integrator approach has traditionally been applied to Lagrangian formulations of mechanical systems, as summarized in \citet{MaWe2001}, and the development of Hamiltonian variational integrators has been less extensive. The notion of Hamiltonian variational integrators was first introduced in \citet{LaWe2006} as the dual formulation of a discrete constrained variational principle, but it did not provide an explicit characterization of the discrete Hamiltonian in terms of the continuous Hamiltonian and the corresponding discrete Noether's theorem, which was introduced in \citet{LeZh2009}. This involves constructing the exact Type II/Type III generating functions for the Hamiltonian flow of a mechanical system, which can be viewed as the analogue of Jacobi's solution of the Hamilton--Jacobi equation. The variational error analysis result for Hamiltonian variational integrators was established in \citet{ScLe2017}, and methods based on Taylor expansions were developed in \citet{ScShLe2017}.

Hamiltonian variational integrators also find application in discrete optimal control and discrete Hamilton--Jacobi theory, and it was shown in \citet{OhBlLe2010} that the Bellman equations of discrete optimal control are the lowest order approximation of a continuous optimal control problem arising from a particular choice of Hamiltonian variational integrator. The Poincar\'e transformed Hamiltonian was used independently by \citet{Ha1997} and \citet{Re1999} as a means of constructing time-adaptive symplectic integrators, and an adaptive approach based on Hamiltonian variational integrators was developed in \citet{DuScLe2020}. The Hamiltonian approach is necessary in this case as many monitor functions result in Poincar\'e transformed Hamiltonians that are degenerate, for which no Lagrangian analogue exists.

In the setting of Lagrangian and Hamiltonian partial differential equations, multisymplectic integrators that can be viewed as generalizations of symplectic integrators for mechanical systems to field theories were introduced from a Lagrangian perspective in \citet{MaPaSh1998}, and from the Hamiltonian, but non-variational perspective, in \citet{BrRe2001}. Our approach to constructing a variational description of multisymplectic integrators for Hamiltonian partial differential equations is based on the notion of generating functionals for multisymplectic relations that was introduced in \citet{VaLiLe2011}.

The advantage of the discrete variational principle approach is that it automatically yields multisymplectic integrators, and exhibit a discrete analogue of Noether's theorem. Furthermore, they naturally lend themselves to Galerkin discretizations that allow for the systematic construction of multisymplectic integrators by choosing a finite-dimensional approximation space for sections of the configuration bundle, and a numerical quadrature rule. In addition, group-invariant discretizations that exhibit a discrete Noether's theorem can be constructed from finite-dimensional approximation spaces that are equivariant with respect to the Lie symmetry group that generates the relevant momentum map.

\subsection{Lagrangian and Hamiltonian Variational Integrators}

Geometric numerical integration aims to preserve geometric conservation laws under discretization, and this field is surveyed in the monograph by \citet{HaLuWa2006}. Discrete variational mechanics \cite{MaWe2001, LeOh2008} provides a systematic method of constructing symplectic integrators. It is typically approached from a Lagrangian perspective by introducing the \textit{discrete Lagrangian}, $L_d:Q\times Q\rightarrow \mathbb{R}$, which is a Type I generating function of a symplectic map and approximates the \textit{exact discrete Lagrangian}, which is constructed from the Lagrangian $L: TQ \rightarrow \mathbb{R}$ as
\vspace*{-1.75ex}
\begin{equation}
L_d^E(q_0,q_1;h)=\ext_{\substack{q\in C^2([0,h],Q) \\ q(0)=q_0, q(h)=q_1}} \int_0^h L(q(t), \dot q(t)) dt,\label{exact_Ld_variational}
\end{equation}
which is equivalent to Jacobi's solution of the Hamilton--Jacobi equation.
The exact discrete Lagrangian generates the exact discrete-time flow map of a Lagrangian system, but, in general, it cannot be computed explicitly. Instead, this can be approximated by replacing the integral with a quadrature formula, and replacing the space of $C^2$ curves with a finite-dimensional function space.

Given a finite-dimensional function space $\mathbb{M}^n([0,h])\subset C^2([0,h],Q)$ and a quadrature formula $\mathcal{G}:C^2([0,h],Q)\rightarrow\mathbb{R}$, $\mathcal{G}(f)=h\sum_{j=1}^m b_j f(c_j h)\approx \int_0^h f(t) dt$, the \textit{Galerkin discrete Lagrangian} is
\[L_d(q_0,q_1)=\ext_{\substack{q\in \mathbb{M}^n([0,h]) \\ q(0)=q_0, q(h)=q_1}} \mathcal{G}(L(q, \dot q))=\ext_{\substack{q\in \mathbb{M}^n([0,h]) \\ q(0)=q_0, q(h)=q_1}} h \sum\nolimits_{j=1}^m b_j L(q(c_j h), \dot q(c_j h)).\]

Given a discrete Lagrangian $L_d$, the \textit{discrete Hamilton--Pontryagin principle} imposes the discrete second-order condition $q_k^1=q_{k+1}^0$ using Lagrange multipliers $p_{k+1}$, which yields a variational principle on $(Q\times Q)\times_Q T^*Q$,
\begin{equation*}
 \delta \left[\sum\nolimits_{k=0}^{n-1} L_d(q_k^0,q_k^1)+\sum\nolimits_{k=0}^{n-2}p_{k+1}(q_{k+1}^0-q_k^1)\right]=0.
 \end{equation*}
This in turn yields the \textit{implicit discrete Euler--Lagrange equations},
\begin{equation}\label{IDEL} q_k^1=q_{k+1}^0,\qquad p_{k+1}=D_2 L_d(q_k^0, q_k^1), \qquad p_{k}=-D_1 L_d(q_k^0, q_k^1),\end{equation}
where $D_i$ denotes the partial derivative with respect to the $i$-th argument. Making the identification $q_k=q_k^0=q_{k-1}^1$, we obtain the \textit{discrete Lagrangian map} and \textit{discrete Hamiltonian map} which are $F_{L_d}:(q_{k-1},q_k)\mapsto(q_k,q_{k+1})$ and $\tilde{F}_{L_d}:(q_k,p_k)\mapsto(q_{k+1},p_{k+1})$, respectively. The last two equations of \eqref{IDEL} define the \textit{discrete fiber derivatives}, $\mathbb{F}L_d^\pm:Q\times Q\rightarrow T^*Q$,
%\begin{align}
%\newsubeqblock
%\mysubeq \mathbb{F}L_d^+(q_k, q_{k+1})&=(q_{k+1},D_2 L_d(q_k, q_{k+1})), \label{FLd+}\\
%\mysubeq \mathbb{F}L_d^-(q_k, q_{k+1})&=(q_k,-D_1 L_d(q_k, q_{k+1})).\label{FLd-}
%\end{align}
\begin{align*}
\mathbb{F}L_d^+(q_k, q_{k+1})&=(q_{k+1},D_2 L_d(q_k, q_{k+1})),\\
\mathbb{F}L_d^-(q_k, q_{k+1})&=(q_k,-D_1 L_d(q_k, q_{k+1})).
\end{align*}
These two discrete fiber derivatives induce a single unique \textit{discrete symplectic form} $\Omega_{L_d}=(\mathbb{F}L_d^\pm)^*\Omega$, where $\Omega$ is the canonical symplectic form on $T^*Q$, and the discrete Lagrangian and Hamiltonian maps preserve $\Omega_{L_d}$ and $\Omega$, respectively.
%%%%%%%%%%%%%%
The discrete Lagrangian and Hamiltonian maps can be expressed as $F_{L_d}=(\mathbb{F}L_d^-)^{-1}\circ \mathbb{F}L_d^+$ and $\tilde{F}_{L_d}=\mathbb{F}L_d^+\circ (\mathbb{F}L_d^-)^{-1}$, respectively. This characterization allows one to relate the approximation error of the discrete flow maps to the approximation error of the discrete Lagrangian.

The variational integrator approach simplifies the numerical analysis of symplectic integrators. The task of establishing the geometric conservation properties and order of accuracy of the discrete Lagrangian map $F_{L_d}$ and discrete Hamiltonian map $\tilde{F}_{L_d}$ reduces to the simpler task of verifying certain properties of the discrete Lagrangian $L_d$ instead.

\begin{theorem}[Discrete Noether's theorem (Theorem 1.3.3 of \cite{MaWe2001})] If a discrete Lagrangian $L_d$ is invariant under the diagonal action of $G$ on $Q\times Q$, then the single unique \textit{discrete momentum map}, $\mathbf{J}_{L_d}=(\mathbb{F}L_d^\pm)^*\mathbf{J}$, is invariant under the discrete Lagrangian map $F_{L_d}$, i.e., $F_{L_d}^* \mathbf{J}_{L_d}=\mathbf{J}_{L_d}$.
\end{theorem}

\begin{theorem}[Variational error analysis (Theorem 2.3.1 of \cite{MaWe2001})\label{thm_variational_error_analysis}]
If a discrete Lagrangian $L_d$ approximates the exact discrete Lagrangian $L_d^E$ to order $p$, i.e., $ L_d(q_0, q_1;h)=L_d^E(q_0,q_1;h)+\mathcal{O}(h^{p+1}),$
then the discrete Hamiltonian map $\tilde{F}_{L_d}$ is an order $p$ accurate one-step method.
\end{theorem}

The bounded energy error of variational integrators can be understood by performing backward error analysis, which then shows that the discrete flow map is approximated with exponential accuracy by the exact flow map of the Hamiltonian vector field of a modified Hamiltonian~\cite{BeGi1994,Ta1994}.

Given a degenerate Hamiltonian, where the Legendre transform $\mathbb{F}H:T^*Q\rightarrow TQ$, $(q,p)\mapsto (q, \frac{\partial H}{\partial p})$, is noninvertible, there is no equivalent Lagrangian formulation. Thus, a characterization of variational integrators directly in terms of the continuous Hamiltonian is desirable. This is achieved by considering the Type II analogue of Jacobi's solution, given by
\[H_d^{+,E}(q_k, p_{k+1})=  \ext_{\substack{(q, p) \in
C^2([t_k,t_{k+1}],T^*Q)\\q(t_k)=q_k, p(t_{k+1})=p_{k+1}}} \Big[ p(t_{k+1}) q (t_{k+1}) - \int_{t_k}^{t_{k+1}} \left[ p \dot{q}-H(q, p) \right]
dt \Big].
\]
A computable Galerkin discrete Hamiltonian $H_d^+$ is obtained by choosing a finite-dimensional function space and a quadrature formula,
\[H_d^+(q_0,p_1)=\ext_{\substack{q\in \mathbb{M}^n([0,h]) \\ q(0)=q_0\\(q(c_j h), p(c_j h))\in T^*Q}} \left[p_1 q(t_1) - h \sum\nolimits_{j=1}^m b_j [ p(c_j h)\dot q(c_j h)-H(q(c_j h),p(c_j h)) ]\right].\]
Interestingly, the Galerkin discrete Hamiltonian does not require a choice of a finite-dimensional function space for the curves in the momentum, as the quadrature approximation of the action integral only depend on the momentum values $p(c_j h)$ at the quadrature points, which are determined by the extremization principle. In essence, this is because the action integral does not depend on the time derivative of the momentum $\dot{p}$. As such, both the Galerkin discrete Lagrangian and the Galerkin discrete Hamiltonian depend only on the choice of a finite-dimensional function space for curves in the position, and a quadrature rule. It was shown in Proposition 4.1 of \citep{LeZh2009} that when the Hamiltonian is hyperregular, and for the same choice of function space and quadrature rule, they induce equivalent numerical methods.

The \textit{Type II discrete Hamilton's phase space variational principle} states that 
\[
\delta \left\{p_N q_N -\sum_{k=0}^{N-1}\left[p_{k+1} q_{k+1} -H_d^+(q_k,
 p_{k+1})\right]\right\}=0,
\]
for discrete curves in $T^*Q$ with fixed $(q_0, p_N)$ boundary conditions. This yields the \textit{discrete Hamilton's equations}, which are given by
\begin{equation}\label{dHamiltonEq}
q_{k+1} = D_2 H_d^+(q_k, p_{k+1}), \qquad p_k=D_1 H_d^+(q_k, p_{k+1}).
\end{equation}

Given a discrete Hamiltonian $H_d^+$, we introduce the \textit{discrete fiber derivatives} (or discrete Legendre transforms), $\mathbb{F}^+ H_d^+$,
\begin{align*}
\mathbb{F}^+H_d^+&: (q_0,p_1)\mapsto (D_2 H_d^+(q_0,p_1),p_1),\\
\mathbb{F}^-H_d^+&: (q_0,p_1)\mapsto (q_0 , D_1 H_d^+ (q_0,p_1)).
\end{align*}
The discrete Hamiltonian map can be expressed in terms of the discrete fiber derivatives,
$$ \tilde{F}_{H_d^+}(q_0,p_0) = \mathbb{F}^+H_d^{+} \circ (\mathbb{F}^-H_d^+)^{-1} (q_0,p_0) = (q_1,p_1) ,$$

Similar to the Lagrangian case, we have a discrete Noether's theorem and variational error analysis result for Hamiltonian variational integrators.

\begin{theorem}[Discrete Noether's theorem (Theorem 5.3 of \cite{LeZh2009})] Let $\Phi^{T^*Q}$ be  the cotangent lift action of  the action $\Phi$ on the configuration manifold $Q$.  If the generalized discrete Lagrangian $R_d(q_0, q_1, p_1)=p_1 q_1 -H_d^+(q_0, p_1)$ is invariant under the cotangent lifted action $\Phi^{T^*Q}$, then the discrete Hamiltonian map $\tilde{F}_{H_d^+}$ preserves the momentum map, i.e., $\tilde{F}_{H_d^+}^* \mathbf{J}=\mathbf{J}$.
\end{theorem}

\begin{theorem}[Variational error analysis (Theorem 2.2 of \cite{ScLe2017})]
If a discrete Hamiltonian $H_d^+$ approximates the exact discrete Hamiltonian $H_d^{+,E}$ to order $p$, i.e., $ H_d^+(q_0, p_1;h)=H_d^{+,E}(q_0,p_1;h)+\mathcal{O}(h^{p+1})$, then the discrete Hamiltonian map $\tilde{F}_{H_d^+}$ is an order $p$ accurate one-step method.
\end{theorem}

It should be noted that there is an analogous theory of discrete Hamiltonian variational integrators based on Type III generating functions $H_d^-(p_0,q_1)$. 

\begin{remark}\label{IntrinsicHamiltonianVariationalIntegrators} It should be noted that the current construction of Hamiltonian variational integrators is only valid on vector spaces and local coordinate charts as it involves Type II/Type III generating functions $H_d^+(q_k. p_{k+1})$, $H_d^-(p_k, q_{k+1})$, which depend on the position at one boundary point, and the momentum at the other boundary point. However, this does not make intrinsic sense on a manifold, since one needs the base point in order to specify the corresponding cotangent space. One possible approach to constructing an intrinsic formulation of Hamiltonian variational integrators is to start with discrete Dirac mechanics~\cite{LeOh2008}, and consider a generating function $E_d^+(q_k, q_{k+1}, p_{k+1})$, $E_d^-(q_k, p_k, q_{k+1})$, that depends on the position at both boundary points and the momentum at one of the boundary points. This approach can be viewed as a discretization of the generalized energy $E(q,v,p)=\langle p,v\rangle - L(q,v)$, in contrast to the Hamiltonian $H(q,p)=\ext_{v}\langle p,v\rangle - L(q,v)=\left.\langle p,v\rangle - L(q,v)\right|_{p=\frac{\partial L}{\partial v}}$.
\end{remark}

\subsection{Multisymplectic Hamiltonian Field Theory}
While classical field theories can be viewed as an infinite-dimensional Hamiltonian system with time as the independent variable (see, for example, \citet{AbMa1978}), we will adopt the multisymplectic formulation with spacetime as the independent variables, which has been extensively studied in, for example, \citet{GoIsMaMo1998, GoIsMaMo2004}, \citet{MaSh1999}, \citet{MaPeShWe2001}. The description of multisymplectic classical field theories in the literature is traditionally formulated in the Lagrangian setting or in the Hamiltonian setting via the covariant Legendre transform to pass between the two settings. However, as we are interested in constructing variational integrators purely within the Hamiltonian setting, we will outline the necessary ingredients of multisymplectic Hamiltonian field theory in this section, without the use of the Lagrangian framework or the covariant Legendre transform.

Consider a trivial vector bundle $E = X \times Q \rightarrow X$ over an oriented spacetime $X$ (although we will refer to $X$ as spacetime with evolutionary Hamiltonian PDEs in mind, $X$ could be either Riemannian or Lorentzian), with volume form denoted $d^{n+1}x$. Let $\Theta$ be the Cartan form on the dual jet bundle $J^1E^*$, which has coordinates $(x^\mu, \phi^A, p, p^A_\mu)$, where $x^\mu$ are the coordinates on spacetime, $\phi^A$ are the coordinates on $Q$, and $p$ and $p^A_\mu$ are the coordinates of the affine map on the jet bundle, $v^A_\mu \mapsto (p + p_A^\mu v^A_\mu)d^{n+1}x$. Define the restricted dual jet bundle $\widetilde{J^1E^*}$ as the quotient of $J^1E^*$ by horizontal one-forms; this space is coordinatized by $(x^\mu,\phi^A,p^A_\mu)$ and is the relevant configuration bundle for a Hamiltonian field theory; we interpret $\phi^A$ as the value of the field and $p_A^\mu$ as the associated momenta in the direction $x^\mu$. The dual jet bundle can be viewed as a bundle over the restricted bundle, $\mu : J^1E^* \rightarrow \widetilde{J^1E^*}$ (see \citet{LePrRoVi2017}). Let $H \in C^\infty(\widetilde{J^1E^*})$ be the Hamiltonian of our theory. This defines a section of $\mu$, in coordinates $\tilde{H}(x^\mu,\phi^A,p_A^\mu) = (x^\mu, \phi^A, -H, p_A^\mu)$ or using the projections $\pi^{j,k}$ from the bundle of $(j+k)$-forms on $E$ to the subbundle of $j$-horizontal, $k$-vertical forms, this can be defined as the set of $z \in J^1E^*$ such that $\pi^{n+1,0}(z) = - H(\pi^{n,1}(z))d^{n+1}x$. Using this section, one can pullback the Cartan form to a form on the restricted bundle,
$$ \Theta_H = \tilde{H}^*\Theta = p_A^\mu d\phi^A \wedge d^nx_\mu - H d^{n+1}x. $$
We then define the action $S^U$ (relative to an arbitrary region $U \subset X$) as a functional on the sections of $\widetilde{J^1E^*}$ (viewed as a bundle over spacetime),
\begin{equation}\label{Hamiltonian Action}
S^U[\phi,p] = \int_U (\phi,p)^* \Theta_H.
\end{equation}
Hamilton's principle states that this action is stationary for compactly supported vertical variations, i.e.,
$$ 0 = dS^U[\phi,p]\cdot V = \int_U(\phi,p)^* i_V d\Theta_H + \underbrace{\int_{\partial U} (\phi,p)^* i_V\Theta_H}_{=0,\ V \Subset U}. $$
Since $U$ is arbitrary, for a sufficiently smooth solution, this gives the strong form of Hamilton's equations, $(\phi,p)^*i_V \Omega_H = 0$, where we defined the multisymplectic form $\Omega_H = -d\Theta_H$. In coordinates, for $V = \delta\phi^A \partial/\partial \phi^A + \delta p_A^\mu \partial/\partial p_A^\mu$, these equations read
$$ \delta \phi^A (\partial_\mu p_A^\mu + \frac{\partial H}{\partial \phi^A})d^{n+1}x + \delta p_A^\mu (-\partial_\mu\phi^A + \frac{\partial H}{\partial p_A^\mu})d^{n+1}x = 0. $$
Since this must hold for $\delta \phi^A, \delta p_A^\mu$ independent, this gives the De Donder--Weyl equations
\begin{subequations}
\begin{align} \label{DDW1} 
\partial_\mu p_A^\mu &= - \frac{\partial H}{\partial \phi^A},\\
\partial_\mu\phi^A &= \frac{\partial H}{\partial p_A^\mu}. \label{DDW2}
\end{align}
\end{subequations}
To write these equations as a multi-Hamiltonian system, define $z^A = (\phi^A, p^0_A, \dots, p^n_A)^T$; it is clear that the De Donder--Weyl equations can be written as
$$ \underbrace{\begin{pmatrix} 0 & -1 & 0 & \dots & 0 \\ 1 & 0 & 0 & \dots & 0 \\ 0 & 0 & 0 & \dots & 0 \\ \vdots & \vdots & \vdots & \ddots & \vdots \\ 0 & 0 & 0 & 0 & 0 \end{pmatrix}}_{\equiv K^0} \partial_0 z^A + \dots + \underbrace{\begin{pmatrix} 0 & 0 & \dots & 0 & -1 \\ 0 & \ddots & \dots & 0 & 0 \\ \vdots & \vdots & \ddots & \vdots & \vdots \\ 0 & 0 & 0 & \ddots & 0 \\ 1 & 0 & 0 & 0 & 0 \end{pmatrix}}_{\equiv K^n} \partial_n z^A = \nabla_{z^A} H, $$
or $K^0 \partial_0z^A + \dots + K^n \partial_n z^A = \nabla_{z^A}H$, where the matrices $K^\mu$ are $(n+2)\times (n+2)$ skew-symmetric matrices which have value $-1$ in the $(0,\mu+1)$ entry and $1$ in the $(\mu+1,0)$ entry (we are indexing the matrices from $0$ to $n+1$), and $0$ everywhere else. This form of the equations was studied in \citet{Br1997}. We can associate to each of these matrices a degenerate two-form on the restricted dual jet bundle,
$$ \omega^\mu \equiv \sum_A d(z^A)^T \otimes K^\mu dz^A = (-dp^\mu_A \otimes d\phi^A + d\phi^A \otimes dp^\mu_A) = d\phi^A \wedge dp^\mu_A. $$
For simplicity of notation, we will implicitly suppress the duality pairing between $(\phi^A)_A$ (valued in $Q$) and $(p^\mu_A)_A$ (valued in $Q^*$) and write this as $\omega^\mu = d\phi \wedge dp^\mu$ (throughout, we will suppress this duality pairing, e.g. $p^\mu\phi \equiv p^\mu_A \phi^A$). Hamilton's equations $(\phi,p)^*i_V\Omega_H = 0$ can then be written as $\omega^\mu(\partial_\mu z, V) = 0$ (sum over $\mu$), which relates the multisymplectic structure $\Omega_H$ to $(n+1)$-pre-symplectic structures $\{\omega^\mu\}$. 

\begin{remark} The multisymplectic structure is more fundamental, since the $\omega^\mu$ were constructed via a particular coordinate representation. In fact, as discussed in \citet{MaSh1999}, the $\omega^\mu$ are a particular coordinate decomposition of the multisymplectic form; in general, the $\omega^\mu$ are not intrinsic unless the dual jet bundle is trivial, although their combination as the multisymplectic form is intrinsic. Since we will utilize Cartesian coordinates on a rectangular mesh for discretization and we will assume trivial bundles for the discrete theory, these coordinate representatives will be simpler to deal with and correspond to the current literature on multisymplectic Hamiltonian integrators. It would be interesting to investigate variational discretizations of field theories where the dual jet bundle is not trivial; in this setting, utilizing the multisymplectic structure is more fundamental.
\end{remark}

\textbf{Multisymplecticity and the Boundary Hamiltonian.} The above Hamiltonian system admits a notion of conserving multisymplecticity, which generalizes the usual notion of symplecticity. In particular, let $V,W$ be two first variations, i.e., vector fields whose flows map solutions of Hamilton's equations again to solutions; then, for any region $U \subset X$, one has the multisymplectic form formula:
\begin{equation}\label{MFF}
\int_{\partial U} (\phi,p)^*(i_V i_W \Omega_H) = 0,
\end{equation}
which follows from $d^2S^U[\phi,p]\cdot(V,W) = 0$ for a solution $(\phi,p)$ of Hamilton's equations. In coordinates where $V = \delta\phi^A \partial/\partial \phi^A + \delta p_A^\mu \partial/\partial p_A^\mu$ and $V = \delta y^A \partial/\partial \phi^A + \delta \pi_A^\mu \partial/\partial p_A^\mu,$ this reads
\begin{align*}
0 &= \int_{\partial U} (\phi,p)^*(i_V i_W \Omega_H) = \int_{\partial U} (\delta \phi^A \delta \pi_A^\mu - \delta y^A \delta p_A^\mu)|_{(\phi,p)} d^nx_\mu = \int_{\partial U} \omega^\mu|_{(\phi,p)}(V,W) d^nx_\mu.
\end{align*}
Applying Stokes' theorem and noting that $U$ is arbitrary, the strong form of the multisymplectic form formula can be expressed $\partial_\mu \omega^\mu = 0$, which holds when evaluated on two first variations at a solution of Hamilton's equations $(\phi,p)$. In terms of our coordinate representation of Hamilton's equations, by taking the exterior derivative of Hamilton's equations, a first variation is a vector field $V$ which satisfies
$$ K^0 dz_0(V) + \dots + K^n dz_n(V) = (D_{zz}H) dz(V),  $$
where $z_\mu \equiv \partial_\mu z.$ One of the aims of this paper is to construct variational integrators for multi-Hamiltonian PDEs which admit a discrete analog of the multisymplectic conservation law for a suitably defined discrete notion of first variations.

Analogous to how the Type II generating functions are utilized in the construction of Galerkin Hamiltonian variational integrators (see \citet{LeZh2009}), we will utilize the boundary Hamiltonian introduced in \citet{VaLiLe2011}, which will act as a generalized Type II generating functional. Consider a domain $U \subset X$ and partition the boundary $\partial U = A \cup B$; we supply fixed field boundary values $\varphi_A$ on $A$ and fixed normal momenta $\pi_B$ on $B$. The boundary Hamiltonian is defined as a functional on these boundary values
\begin{align}\label{Boundary Hamiltonian}
H_{\partial U}(\varphi_A,\pi_B) &= \ext\Big[ \int_B p^\mu \phi d^nx_\mu - \int_U (\phi,p)^* \Theta_H  \Big] \\ \nonumber
&= \ext\Big[ \int_B p^\mu \phi d^nx_\mu - \int_U (p^\mu \partial_\mu\phi - H(\phi,p) ) d^{n+1}x \Big],
\end{align}
where one extremizes over all fields $(\phi,p)$ satisfying the fixed boundary conditions along $A$ and $B$. 

An extremizer of the above expression restricted to the aforementioned boundary conditions satisfies the De Donder--Weyl equations, which follows from
\begin{align*}
\delta &\Big[ \int_B p^\mu \phi d^nx_\mu - \int_U (p^\mu \partial_\mu\phi - H(\phi,p) ) d^{n+1}x \Big] 
\\ &= \int_B \cancel{\delta p^\mu} \phi d^nx_\mu + \int_B p^\mu \delta \phi d^nx_\mu - \int_U (\delta p^\mu \partial_\mu\phi + p^\mu \partial_\mu \delta\phi - \frac{\partial H(\phi,p)}{\partial p^\mu}\delta p^\mu - \frac{\partial H(\phi,p)}{\partial \phi}\delta \phi  ) d^{n+1}x
\\ &= \int_B p^\mu \delta \phi d^nx_\mu - \int_{\partial U = A \cup B} p^\mu \delta\phi d^nx_\mu  - \int_U (\delta p^\mu \partial_\mu\phi - \partial_\mu p^\mu \delta\phi - \frac{\partial H(\phi,p)}{\partial p^\mu}\delta p^\mu - \frac{\partial H(\phi,p)}{\partial \phi}\delta \phi  ) d^{n+1}x \\
&= - \int_{\partial U = A} p^\mu \cancel{\delta\phi} d^nx_\mu  - \int_U \left[(\partial_\mu\phi - \frac{\partial H(\phi,p)}{\partial p^\mu})\delta p^\mu - (\partial_\mu p^\mu  + \frac{\partial H(\phi,p)}{\partial \phi})\delta\phi \right] d^{n+1}x,\\
\end{align*}
where we used $\delta p^\mu|_{B} = 0 = \delta\phi|_A$.

This is a Type II generating functional in the sense that it generates the boundary values for the field along $B$ (denoted $\phi|_B$) and the normal momenta along $A$ (denoted $p^n|_A$),
\begin{equation}\label{Generating Functional}
 \frac{\delta H_{\partial U}}{\delta \varphi_A} = -p^n|_A, \quad \frac{\delta H_{\partial U}}{\delta \pi_B} = \phi|_B.
\end{equation}
To obtain (\ref{Generating Functional}), perform an analogous computation as the one above (take the variation, integrate by parts, and use that the internal field satisfies the De Donder--Weyl equations), which gives
$$ d H_{\partial U}(\varphi_A,\pi_B)\cdot (\delta\varphi_A,\delta\pi_B) = \int_B \delta \pi_B \cdot \phi|_B - \int_A \pi|_A \cdot \delta\varphi_A; $$
i.e., (\ref{Generating Functional}). Note that the generating relation (\ref{Generating Functional}) only determines the normal component of the momentum along $A$; this is consistent with the De Donder--Weyl equation (\ref{DDW1}), since it only specifies $\partial_\mu p^\mu$.

Since an extremizer of $H_{\partial U}(\varphi_A,\pi_B)$ satisfies the De Donder--Weyl equations, it satisfies the multisymplectic form formula. Since the multisymplectic form formula is expressed as an integral over $\partial U$ and the generating functional gives us the field values on $\partial U$, $(\varphi,\pi) = (\varphi_A,\varphi_B,\pi_A,\pi_B)$, the above generating map (\ref{Generating Functional}) is multisymplectic in the sense 
$$ \int_{\partial U}\omega^\mu|_{(\varphi,\pi)}(V,W) d^nx_\mu = 0, $$
for first variations $V$ and $W$.

We will utilize a discrete approximation of the boundary Hamiltonian and its property as a generating functional to construct variational integrators which are naturally multisymplectic.

\textbf{Noether's Theorem.} Another important conservative property of Hamiltonian systems arises from symmetries. Suppose there is a smooth group action of $G$ on the restricted dual jet bundle which leaves the action $S^U$ invariant. Let $\tilde{\xi}$ denote the infinitesimal generator vector field for $\xi \in \mathfrak{g}$ associated to this action. For a solution $(\phi,p)$ of Hamilton's equations, one has
$$ 0 = \pounds_{\tilde{\xi}}S^U[\phi,p] = dS^U[\phi,p]\cdot \tilde{\xi} = \int_U (\phi,p)^* i_{\tilde{\xi}} d\Theta_H + \int_{\partial U}(\phi,p)^* i_{\tilde{\xi}}\Theta_H. $$
Note that the term involving the integral over $U$ vanishes, even though $\tilde{\xi}$ is not necessarily compactly supported in $U$, since Hamilton's equations hold pointwise ($U$ is arbitrary). Hence, Noether's theorem in this setting is the statement
\begin{equation}\label{Noether's Theorem}
\int_{\partial U}(\phi,p)^* i_{\tilde{\xi}}\Theta_H = 0.
\end{equation}
In the discrete setting, we will be particularly concerned with vertical variations (where the group action on the base space $X$ is the identity). In this case, we can write the above in coordinates as
\begin{equation}\label{Noether's Theorem coordinates}
\int_{\partial U} p^\mu (i_{\tilde{\xi}}d\phi) d^nx_\mu = 0.
\end{equation}
We will see that if there is a group action on the discrete analog of the restricted dual jet bundle which leaves the discrete action (the generalized discrete Lagrangian) invariant, then there is a discrete analog of Noether's theorem, equation (\ref{Noether's Theorem coordinates}).

\subsection{Multisymplectic Integrators for Hamiltonian PDEs}\label{Multisymplectic Integrators Subsection}
Consider the class of Hamiltonian PDEs,
\begin{equation}\label{HamiltonianPDEs}
K^0 z_0 + \dots + K^n z_n = \nabla_z H(z),
\end{equation}
with independent variable $x = (x^0,\dots,x^n) \in \mathbb{R}^{n+1}$, dependent variable $z: \mathbb{R}^n \rightarrow \mathbb{R}^m$, each $K^\mu$ is an $m \times m$ skew-symmetric matrix, and the Hamiltonian $H: \mathbb{R}^m \rightarrow \mathbb{R}$ is sufficiently smooth. 

Defining a two-form for each $K^\mu$, $\omega^\mu(U,V) = \langle K^\mu U,V\rangle$ (with respect to an inner product $\langle\cdot,\cdot\rangle$ on $\mathbb{R}^m$), the equation (\ref{HamiltonianPDEs}) admits the multisymplectic conservation law
\begin{equation}\label{Multisymplectic Conservation Law}
\partial_\mu \omega^\mu(U,V) = 0,
\end{equation}
for any pair of first variations $U,V$ satisfying the variational equation
$$ K^0 dz_0 + \dots + K^n dz_n = D_{zz}H(z). $$

As we saw, the De Donder--Weyl equations, which arose from the variational principle applied to the Hamiltonian action (\ref{Hamiltonian Action}), are an example of a Hamiltonian PDE in the form (\ref{HamiltonianPDEs}). From our variational perspective, the action and variational principle are more fundamental, as opposed to the field equations (\ref{HamiltonianPDEs}). However, as shown by \citet{Ch2005}, the Hamiltonian system (\ref{HamiltonianPDEs}) arises from the variational principle, so there is no loss of generality working with the formulation based on the Hamiltonian action (\ref{Hamiltonian Action}).

For the Hamiltonian system (\ref{HamiltonianPDEs}), a multisymplectic integrator is defined in \citet{BrRe2001} to be a method
$$ K^0 \partial_0^{i_0\dots i_n}z_{i_0\dots i_n} + \dots + K^n \partial_n^{i_0\dots i_n} z_{i_0\dots i_n} = (\nabla_zS(z_{i_0\dots i_n}))_{i_0\dots i_n},$$
where $\partial_\mu^{i_0\dots i_n}$ is a discretization of $\partial_\mu$, such that a discrete analog of equation (\ref{Multisymplectic Conservation Law}) holds,
$$ \partial_\mu^{i_0\dots i_n} \omega^\mu(U_{i_0\dots i_n},V_{i_0\dots i_n}) = 0,$$
when evaluated on discrete first variations $U_{i_0\dots i_n},V_{i_0\dots i_n}$ satisfying the discrete variational equations
$$ K^0 \partial_0^{i_0\dots i_n}dz_{i_0\dots i_n} + \dots + K^n \partial_n^{i_0\dots i_n} dz_{i_0\dots i_n} = d\Big((\nabla_zS(z_{i_0\dots i_n}))_{i_0\dots i_n}\Big). $$
We will see that the variational integrators that we construct will automatically satisfy a discrete multisymplectic conservation law, as a consequence of the Type II variational principle. Furthermore, we will show in Section \ref{Multisymplecticity Revisited} that this discrete multisymplectic conservation law reproduces the Bridges and Reich notion of multisymplecticity. 

\begin{example}
An example of a multisymplectic integrator in $1+1$ spacetime dimensions is the centered Preissman scheme,
$$ K^0 \frac{z^1_{1/2} -z^0_{1/2}}{\Delta t} + K^1 \frac{z^{1/2}_1 - z^{1/2}_0}{\Delta x} = \nabla_z H\Bigl(z^{1/2}_{1/2}\Bigr), $$
where $z^0_{1/2} = \frac{1}{2}(z^0_0 + z^0_1),$ etc. and $z^{1/2}_{1/2} = \frac{1}{4}(z^1_1 + z^0_1 + z^1_0 + z^0_0)$. As noted in \citet{Re2000b}, this can be obtained from a cell-vertex finite volume discretization on a rectangular grid, or alternatively, as observed in \citet{Re2000}, it is an example of a multisymplectic Gauss--Legendre collocation method, in the case of one collocation point. Furthermore, the multisymplectic Gauss--Legendre collocation methods are members of a larger class of multisymplectic integrators,  the multisymplectic partitioned Runge--Kutta methods (see, for example, \citet{HoLiSu2006}, \citet{RyMcFr2007}). In Section \ref{MPRK subsection}, we will derive the class of multisymplectic partitioned Runge--Kutta methods within our variational framework.
\end{example}

\subsection{Main Contributions}
In this paper, we introduce a variational construction of multisymplectic Hamiltonian integrators utilizing a discrete approximation of the boundary Hamiltonian and the corresponding Type II variational principle. Although variational integrators have been extensively studied in the setting of Lagrangian PDEs, where they have been used to construct robust and flexible numerical methods for nonlinear elasticity~\cite{LeMaOrWe2003}, collision and impact dynamics for continuum mechanics~\cite{DeGBRa2016}, and geometrically exact beam dynamics~\cite{LeSaLe2021}, the variational perspective has not been studied in the setting of integrators for Hamiltonian PDEs. 

This paper serves as a stepping stone in constructing variational integrators in the Hamiltonian PDE setting. Our hope is that, by introducing a variational perspective in the setting of integrators for Hamiltonian PDEs, the well-developed techniques and machinery of variational integrators for Lagrangian PDEs can be analogously developed on the Hamiltonian side. It should be noted that the theory in this paper relies on a trivial configuration bundle, since the notion of a boundary Hamiltonian is only intrinsic in the case that the bundle is trivial. Analogous to an intrinsic approach to variational integrators for Hamiltonian mechanics, outlined in Remark \ref{IntrinsicHamiltonianVariationalIntegrators}, one possible approach for constructing an intrinsic formulation of multisymplectic integrators is to start with a discrete notion of a multi-Dirac structure (for details on multi-Dirac structures in classical field theories, see  \citet{VaYoLeMa2012}) and discretize the variational principle utilizing the generalized energy as a generating functional; we will investigate this in future work. 

In Section \ref{Discrete Hamiltonian Field Theory}, we begin by developing a discrete notion of Hamiltonian field theory, the discrete boundary Hamiltonian, and the corresponding Type II variational principle. Subsequently, we specialize to the case of a spacetime tensor product rectangular mesh which allows us to give an explicit characterization of the equations resulting from the Type II variational principle. We prove discrete analogues of multisymplecticity and Noether's theorem for these equations. In Section \ref{GHVI}, we utilize a Galerkin approximation of the action to complete the discretization of the boundary Hamiltonian. Subsequently, in Section \ref{MPRK subsection}, we utilize a particular choice of Galerkin approximation to derive the class of multisymplectic partitioned Runge--Kutta methods. In Section \ref{Multisymplecticity Revisited}, we reinterpret the discrete multisymplectic conservation law as one that is naturally associated to the difference equations which approximate the De Donder--Weyl equations. Finally, in Section \ref{Numerical Section}, we provide a numerical example which allows us to visualize multisymplecticity as symplecticity in the spatial and temporal directions for the class of sine--Gordon soliton solutions.

\section{Multisymplectic Hamiltonian Variational Integrators}

\subsection{Discrete Hamiltonian Field Theory}\label{Discrete Hamiltonian Field Theory}
We will discuss our construction of a discrete boundary Hamiltonian for the general case of an arbitrary mesh and subsequently study the particular case of a rectangular mesh where the variational equations can be written explicitly. Let $X \subset \mathbb{R}^{n+1}$ be a polygonal domain and $\mathcal{T}(X)$ an associated mesh. In general, a discrete configuration bundle consists of a choice of finite element space taking values in the fiber $Q$ that is subordinate to the mesh $\mathcal{T}(X)$. To be more concrete, for every mesh element $\bigtriangleup \in \mathcal{T}(X)$, we introduce nodes $x_i \in \bigtriangleup, i \in I$, and parametrize the finite element space by the fiber value at each node. A multisymplectic variational integrator based on finite elements was developed from the Lagrangian perspective in \citet{Ch2008}.

The discrete analog of the configuration bundle, on an element by element level, is the base space $\{x_i\}_{i \in I}$ with fiber $Q$ over each node; the total space is $\{x_i\}_{i \in I} \times Q$ and a section is a map from each node to $Q$, denoted $\phi_i \in Q$. Analogously, the discrete analog of the restricted dual jet bundle is $\{x_i\}_{i \in I} \times Q \times (Q^*)^{n+1}$, where a section is specified by $\phi_i \in Q, p^\mu_{i} \in Q^*$. Let $S_d^{\bigtriangleup}[\phi_i,p^\mu_i]$ be some discrete approximation of the action $S^\bigtriangleup[\phi,p]$. As in the discussion of the boundary Hamiltonian (\ref{Boundary Hamiltonian}), partition the boundary of the element $\partial\bigtriangleup = A \cup B$ and let $\sumint_B \pi_B \varphi_B$ be some discrete approximation to the boundary integral $\int_B p^\mu \phi d^nx_\mu$, depending only on the field and normal momenta boundary values on the nodes $x_i \in B$, which we denoted $\varphi_B$ and $\pi_B$ respectively. Define the discrete boundary Hamiltonian
$$ H_d^{\partial \bigtriangleup} (\varphi_A,\pi_B) = \ext_{ \stackrel{\phi_i \in Q, p^\mu_{i} \in Q^*}{\phi|_A = \varphi_A, p^n|_B = \pi_B} } \Big[ \sumint \pi_B\varphi_B - S_d^\bigtriangleup[\phi_i,p^\mu_i] \Big], $$
where $p^n|_B$ denotes the normal component of the momenta along $B$. Repeat the above construction for each $\bigtriangleup \in \mathcal{T}(X)$; partitioning the boundaries $\partial \bigtriangleup = A(\bigtriangleup) \cup B(\bigtriangleup)$ and the boundary of the full region $\partial X = A(X) \cup B(X)$ (where $A(X) = \cup_{\Delta \in \mathcal{T}(X)} (A(\Delta) \cap \partial X)$ and $B(X) = \cup_{\Delta \in \mathcal{T}(X)} (B(\Delta) \cap \partial X)$). Define the discrete action sum
$$ S_d[\{\varphi_{A(\bigtriangleup)},\pi_{B(\bigtriangleup)}\}_{\bigtriangleup \in \mathcal{T}(X)} ] = \sumint_{B(X)} \pi_{B(X)}\varphi_{B(X)} - \sum_{\bigtriangleup \in \mathcal{T}(X)} \Big[ \sumint_{B(\bigtriangleup)} \pi_{B(\bigtriangleup)}\varphi_{B(\bigtriangleup)} -  H_d^{\partial \bigtriangleup}(\varphi_A,\pi_B) \Big]. $$
The Type II variational principle $\delta S_d = 0$ (subject to variations of $\varphi$ vanishing along $A(X)$ and variations of $\pi$ vanishing along $B(X)$) gives a set of (generally coupled) maps $(\varphi_{A(\bigtriangleup)}, \pi_{B(\bigtriangleup)}) \mapsto (\varphi_{B(\bigtriangleup)}, \pi_{A(\bigtriangleup)})$ in analogy with the generating functional relation, equation (\ref{Generating Functional}). In the case of finite element spaces which are not parametrized by the nodal values, we evaluate the discrete boundary Hamiltonian on the discrete space of boundary data induced by the choice of mesh and discrete configuration bundle, and extremize the expressions above over the finite elements that satisfy the prescribed boundary conditions. This is the most general form of our multisymplectic Hamiltonian variational integrator.

\textbf{Spacetime Tensor Product Rectangular Mesh.} Now, consider the particular case of a rectangular domain $X$ and an associated rectangular mesh $\mathcal{T}(X)$. For simplicity and clarity in the notation, we will focus on the case of $1+1$ spacetime dimensions, although higher dimensions can be treated similarly (we treat the case of higher dimensions in Appendix \ref{HigherDimensions}).

Consider a rectangle $[t, t + \Delta t] \times [x, x + \Delta x] = \Box \in \mathcal{T}(X)$. Introduce nodes on the intervals $\{t_1 = t, t_2, \dots , t_{s-1}, t_s = t + \Delta t\}$ and $\{x_1 = x, x_2, \dots, x_{\sigma-1}, x_\sigma = x + \Delta x\}$ (as we will introduce in the next section for Galerkin Hamiltonian variational integrators, these nodes correspond to quadrature points along the time and space intervals). The discrete base space is $X_d=\{(t_i,x_j)\mid i=1,\ldots, s,\ j=1,\ldots, \sigma\}$, the discrete configuration bundle is $X_d \times Q$, where a section is map from each node $(t_i,x_j)$ to $(t_i,x_j, \phi_{ij})$, where $\phi_{ij} \in Q$. Analogously, the discrete restricted dual jet bundle is $X_d \times Q \times (Q^*)^{2}$, where a section is specified by $\phi_{ij} \in Q, p^\mu_{ij} \in Q^*$. Let $S^\Box_d [\phi_{ij},p^\mu_{ij}]$ be some discrete approximation to $S^\Box[\phi,p]$ (we will explicitly construct such a discrete approximation in the next section using Galerkin techniques and quadrature). Partitioning the boundary $\partial\Box = A(\Box) \cup B(\Box)$, the discrete boundary Hamiltonian is given by
\begin{equation}\label{Discrete Boundary Hamiltonian}
H^{\partial\Box}_d (\varphi_{A(\Box)},\pi_{B(\Box)}) = \ext_{\stackrel{\phi_{ij} \in Q, p^\mu_{ij} \in Q^*}{\phi|_{A(\Box)} = \varphi_{A(\Box)}, p^n|_{B(\Box)}= \pi_{B(\Box)} }} \Big[ \sumint_{B(\Box)} \pi_{B(\Box)}\varphi_{B(\Box)} - S^\Box_d[\phi_{ij},p^\mu_{ij}] \Big],
\end{equation}
where $\varphi_{A(\Box)}$ denotes the boundary values on $A(\Box)$, i.e., at nodes $(t_i,x_j) \in A$ (and similarly for $\pi$). The discrete action sum is
$$ S_d\big[\{\varphi_{A(\Box)}, \pi_{B(\Box)}\}_{\Box \in \mathcal{T}(X)}\big] = \sumint_{B(X)} \pi_{B(X)}\varphi_{B(X)} - \sum_{\Box \in \mathcal{T}(X)} \Big[ \sumint_{B(\Box)} \pi_{B(\Box)} \varphi_{B(\Box)} - H_d^{\partial \Box}(\varphi_{A(\Box)},\pi_{B(\Box)}) \Big]. $$

\begin{figure}[h]
\begin{center}
\includegraphics[width=90mm]{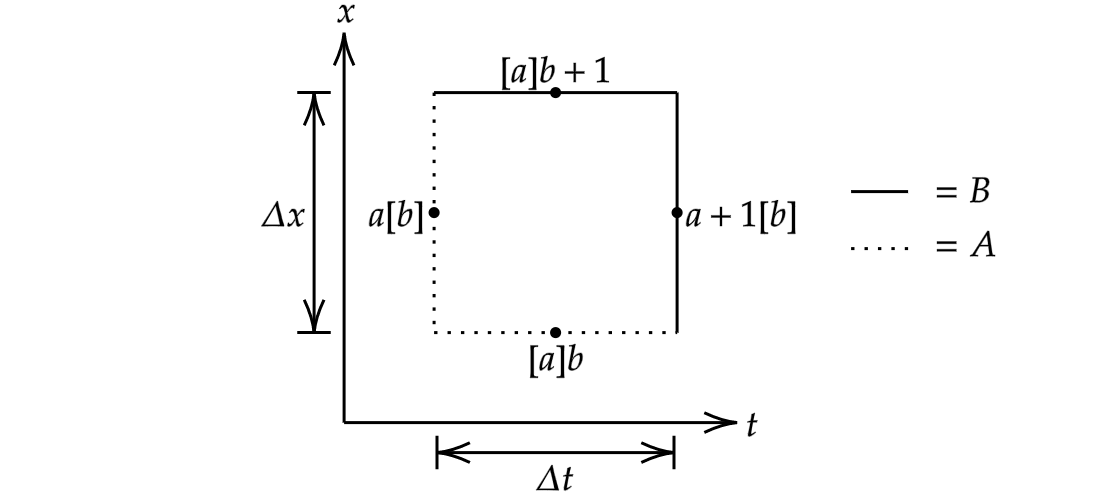}
\caption{Schematic for one quadrature point along each edge of $\Box_{ab} \in \mathcal{T}(X)$.}
\label{onequad}
\end{center}
\end{figure}

Recall the Type II variational principle $\delta S_d = 0$ gives a set of maps $(\varphi_{A(\Box)},\pi_{B(\Box)}) \mapsto (\varphi_{B(\Box)},\pi_{A(\Box)})$. To give a more explicit characterization of these maps, let us introduce a quadrature approximation of the boundary integral over $B$. First, consider the simple case of one quadrature point along each edge of $\Box_{ab} = [t_0+a\Delta t, t_0 + (a+1)\Delta t] \times [x_0+b\Delta x, x_0 + (b+1)\Delta x]$, where $\mathcal{T}(X) = \{ \Box_{ab} \}_{a,b}$. Let $\varphi_{[a]b}$ denote the field boundary value at the quadrature point along the bottom edge $(t_a,t_a+\Delta t)\times \{x_b\}$ (where we orient our axes such that time is horizontal and space is vertical) and $\varphi_{a[b]}$ denote its value at the quadrature point along the left edge $\{t_a\} \times (x_b,x_b+\Delta x)$ (and similarly $\varphi_{[a]b+1}$ for the top edge, $\varphi_{a+1[b]}$ for the right edge). We take $A$ to be the bottom and left edges, and $B$ to the top and right edges. The normal momenta through the top edge is the momenta associated to the $x$ direction (at the quadrature point), which we denote $\pi^1_{[a]b+1}$, and the normal momenta through the right edge is the momenta associated to the $t$ direction, which we denote $\pi^0_{a+1[b]}$. Since we only have one quadrature point along each edge, the quadrature weight for the temporal edge is $\Delta t$ and similarly for the spatial edge is $\Delta x$. See Figure \ref{onequad}.

Then, the boundary integral can be approximated
\begin{align*}
\int_B p^\mu \phi d^nx_\mu &= \int_{t_a}^{t_{a+1}} (p^1\phi)|_{x = x_{b+1}} dt + \int_{x_b}^{x_{b+1}} (p^0\phi)_{t = t_{a+1}} dx \\ 
&\approx \pi^1_{[a]b+1}\varphi_{[a]b+1} \Delta t + \pi^0_{a+1[b]}\varphi_{a+1[b]} \Delta x \equiv \sumint_B \pi_B\varphi_B.
\end{align*}
The associated discrete boundary Hamiltonian is
$$ H_d^+(\varphi_{[a]b},\varphi_{a[b]},\pi^1_{[a]b+1},\pi^0_{a+1[b]}) = \ext \Big( \pi^1_{[a]b+1}\varphi_{[a]b+1} \Delta t + \pi^0_{a+1[b]}\varphi_{a+1[b]} \Delta x - S_d^{\Box_{ab}}[\phi,p]\Big), $$
where the $+$ specifies that we chose $B$ to be in the forward direction (in the direction of increasing temporal and spatial values), analogous to the notion of discrete right Hamiltonian in discrete mechanics. Again, we extremize over $\phi,p$ satisfying the boundary conditions (note we have not given an explicit construction for such a $S_d^{\Box_{ab}}$ yet; see Section \ref{GHVI}).

\begin{prop}\label{Discrete Hamilton's Equations Prop}
The Type II variational principle $\delta S_d = 0$, subject to variations of $\varphi$ vanishing along $A(X)$ and variations of $\pi$ vanishing along $B(X)$, yields the following,
\begin{subequations}
\begin{align}\label{Discrete Hamilton's Equations, one quad, 1}
\pi^1_{[a]b} &= \frac{1}{\Delta t}D_1H_d^{+}(\varphi_{[a]b},\varphi_{a[b]},\pi^1_{[a]b+1},\pi^0_{a+1[b]}), \\
\pi^0_{a[b]} &= \frac{1}{\Delta x}D_2H_d^{+}(\varphi_{[a]b},\varphi_{a[b]},\pi^1_{[a]b+1},\pi^0_{a+1[b]}),\\
\varphi_{[a]b+1} &= \frac{1}{\Delta t}D_3H_d^{+}(\varphi_{[a]b},\varphi_{a[b]},\pi^1_{[a]b+1},\pi^0_{a+1[b]}),\\
\label{Discrete Hamilton's Equations, one quad, 2}
\varphi_{a+1[b]} &= \frac{1}{\Delta x} D_4H_d^{+}(\varphi_{[a]b},\varphi_{a[b]},\pi^1_{[a]b+1},\pi^0_{a+1[b]}),
\end{align}
\end{subequations}
where $D_i$ denotes differentiation with respect to the $i^{th}$ argument. We refer to these equations as the discrete forward Hamilton's equations (in the case of one quadrature point). Note that these equations define a map $(\varphi_{A},\pi_B) = (\varphi_{[a]b},\varphi_{a[b]},\pi^1_{[a]b+1},\pi^0_{a+1[b]}) \mapsto (\varphi_B,\pi_A) = (\varphi_{[a]b+1},\varphi_{a+1[b]},\pi^1_{[a]b},\pi^0_{a[b]}). $ 
\begin{proof}
Recall the full mesh $\mathcal{T}(X) = \{\Box_{ab}\}_{a,b}$; say $a=0,\ldots,N-1$, and $b=0,\ldots,M-1$ (so that $X = [t_0,t_0 + N \Delta t] \times [x_0, x_0 + M \Delta x]$). $B(X)$ consists of the forward edges of $X$, i.e., 
$$ B(X) = \Big( [t_0, t_0 + N \Delta t] \times \{x_0 + M \Delta x\} \Big) \cup \Big( \{t_0 + N \Delta t\} \times [x_0,x_0 + M \Delta x] \Big). $$
Consider the discrete action sum
\begin{align*}
\hspace{-2cm} S_d&[\{\varphi_{A(\Box)},\pi_{B(\Box)}\}]\\
&= \sumint_{B(X)} \pi_{B(X)}\varphi_{B(X)} - \sum_{\Box \in \mathcal{T}(X)} \Big[ \sumint_{B(\Box)} \pi_{B(\Box)}\varphi_{B(\Box)} - H_d^+(\varphi_{A(\Box)},\pi_{B(\Box)}) \Big]\\
&= \sum_{a=0}^{N-1} \pi^1_{[a]M}\varphi_{[a]M} \Delta t + \sum_{b=0}^{M-1} \pi^0_{N[b]}\varphi_{N[b]} \Delta x \\
&\qquad\ -\sum_{a,b=0}^{N-1,M-1} \Big[\pi^1_{[a]b+1}\varphi_{[a]b+1} \Delta t + \pi^0_{a+1[b]}\varphi_{a+1[b]} \Delta x -  H_d^+(\varphi_{[a]b},\varphi_{a[b]},\pi^1_{[a]b+1},\pi^0_{a+1[b]}) \Big] \\
&= \underbrace{-\sum_{a,b=0}^{N-1,M-2} \pi^1_{[a]b+1}\varphi_{[a]b+1} \Delta t}_{\equiv(a)}\ \underbrace{ - \sum_{a,b = 0}^{N-2,M-1} \pi^0_{a+1[b]}\varphi_{a+1[b]} \Delta x}_{\equiv (b)}\\
&\qquad+ \underbrace{\sum_{a,b = 0}^{N-1,M-1}  H_d^+(\varphi_{[a]b},\varphi_{a[b]},\pi^1_{[a]b+1},\pi^0_{a+1[b]})}_{\equiv(c)}.
\end{align*}
The Type II variational principle states $0 = \delta S_d = \delta (a) + \delta (b) + \delta (c)$, subject to variations of $\varphi$ vanishing along $A(X)$ (i.e., $\delta \varphi_{[a]0} = 0 = \delta\varphi_{0[b]}$) and variations of $\pi$ vanishing along $B(X)$ (i.e., $\delta \pi^0_{N[b]} = 0 = \delta \pi^1_{[a]M} $). Compute the variations of $(a),(b),(c)$ keeping only the independent variations $\delta \varphi_{[a]b}$, $\delta\varphi_{a[b]}$, $\delta \pi^0_{a[b]}$, $\delta \pi^1_{[a]b}$ not required to vanish by the boundary conditions (note such vanishing variations will only appear in $(c)$).
\begin{align*}
\delta (a) &= - \Delta t \sum_{a=0}^{N-1} \sum_{b=0}^{M-2} \Big(\varphi_{[a]b+1} \delta \pi^1_{[a]b+1} + \pi^1_{[a]b+1} \delta\varphi_{[a]b+1} \Big) \\
&=- \Delta t \sum_{a=0}^{N-1} \sum_{b=0}^{M-2} \varphi_{[a]b+1} \delta \pi^1_{[a]b+1} - \Delta t \sum_{a=0}^{N-1} \sum_{b=1}^{M-1} \pi^1_{[a]b} \delta\varphi_{[a]b}, \\
\delta (b) &= - \Delta x \sum_{a=0}^{N-2} \sum_{b=0}^{M-1} \Big( \varphi_{a+1[b]} \delta \pi^0_{a+1[b]} + \pi^0_{a+1[b]}\delta\varphi_{a+1[b]}\Big) \\
&=- \Delta x \sum_{a=0}^{N-2} \sum_{b=0}^{M-1}\varphi_{a+1[b]} \delta \pi^0_{a+1[b]} - \Delta x \sum_{a=1}^{N-1} \sum_{b=0}^{M-1} \pi^0_{a[b]}\delta\varphi_{a[b]}.
\end{align*}
For brevity, denote $H_d^+[a,b] \equiv H_d^+(\varphi_{[a]b},\varphi_{a[b]},\pi^1_{[a]b+1},\pi^0_{a+1[b]})$. Compute
\begin{align*}
\delta (c) &= \sum_{a,b = 0}^{N-1,M-1} \Big( D_1H_d^+[a,b] \delta\varphi_{[a]b} + D_2H_d^+[a,b] \delta\varphi_{a[b]} + D_3H_d^+[a,b] \delta \pi^1_{[a]b+1} + D_4H_d^+[a,b] \delta \pi^0_{a+1[b]} \Big) \\
&=\sum_{a = 0}^{N-1} \sum_{b=0}^{M-1}  D_1H_d^+[a,b] \delta\varphi_{[a]b} + \sum_{a = 0}^{N-1} \sum_{b=0}^{M-1} D_2H_d^+[a,b] \delta\varphi_{a[b]} \\
&\qquad+ \sum_{a = 0}^{N-1} \sum_{b=0}^{M-1} D_3H_d^+[a,b] \delta \pi^1_{[a]b+1} + \sum_{a = 0}^{N-1} \sum_{b=0}^{M-1} D_4H_d^+[a,b] \delta \pi^0_{a+1[b]}.
\end{align*}
Note in the first double sum above, $\delta\varphi_{[a]0} = 0$ so we remove the $b=0$ terms. In the second double sum, $\delta\varphi_{0[b]} = 0$ so we remove the $a= 0$ terms. In the third double sum above, $\delta \pi^1_{[a]M} = 0$ so we remove the $b = M-1$ terms. In the fourth double sum above, $\delta \pi^0_{N[b]} = 0$ so we remove the $a = N-1$ terms. This gives,
\begin{align*}
\delta (c) &= \sum_{a = 0}^{N-1} \sum_{b=1}^{M-1}  D_1H_d^+[a,b] \delta\varphi_{[a]b} + \sum_{a = 1}^{N-1} \sum_{b=0}^{M-1} D_2H_d^+[a,b] \delta\varphi_{a[b]} \\
&\qquad+ \sum_{a = 0}^{N-1} \sum_{b=0}^{M-2} D_3H_d^+[a,b] \delta \pi^1_{[a]b+1} + \sum_{a = 0}^{N-2} \sum_{b=0}^{M-1} D_4H_d^+[a,b] \delta \pi^0_{a+1[b]}.
\end{align*}
Putting everything together, we have 
\begin{align*}
0 &= \delta S_d = \delta(a) + \delta (b) + \delta (c) \\
&= \sum_{a = 0}^{N-1} \sum_{b=1}^{M-1} (-\Delta t\ \pi^1_{[a]b} +  D_1H_d^+[a,b]) \delta\varphi_{[a]b} + \sum_{a = 1}^{N-1} \sum_{b=0}^{M-1} (-\Delta x\  \pi^0_{a[b]} + D_2H_d^+[a,b]) \delta\varphi_{a[b]} \\
&\qquad+\sum_{a = 0}^{N-1} \sum_{b=0}^{M-2}(- \Delta t\ \varphi_{[a]b+1} + D_3H_d^+[a,b] )\delta \pi^1_{[a]b+1} + \sum_{a = 0}^{N-2} \sum_{b=0}^{M-1}(-\Delta x\ \varphi_{a+1[b]} + D_4H_d^+[a,b] ) \delta \pi^0_{a+1[b]}.
\end{align*}
The variations in the above expression are all independent, so this gives (\ref{Discrete Hamilton's Equations, one quad, 1})-(\ref{Discrete Hamilton's Equations, one quad, 2}).
\end{proof}
\end{prop}

\textbf{Discrete Multisymplecticity.} Analogous to the continuum case, we define a discrete first variation as a vector field such that the above equations (\ref{Discrete Hamilton's Equations, one quad, 1})-(\ref{Discrete Hamilton's Equations, one quad, 2}) still hold when evaluated at the level of the exterior derivative, e.g. for equation (\ref{Discrete Hamilton's Equations, one quad, 1}),
$$d\pi^1_{[a]b} = \frac{1}{\Delta t} d \Big( D_1H_d^{+}(\varphi_{[a]b},\varphi_{a[b]},\pi^1_{[a]b+1},\pi^0_{a+1[b]})\Big).$$
and similarly for the others. As we saw in the continuum theory, the map generated by the boundary Hamiltonian implies the multisymplectic form formula, since the multisymplectic form formula can be expressed over the boundary $\partial U$. Since we constructed a discrete approximation to the boundary Hamiltonian before enforcing the variational principle, we would naturally expect a discrete notion of multisymplecticity to arise as well. Furthermore, in the continuum theory, multisymplecticity follows from $d^2=0$ applied to the boundary Hamiltonian, evaluated on first variations. As we will see, our discrete multisymplectic form formula follows from computing $d^2 = 0$ applied to the discrete boundary Hamiltonian, in analogy with the continuum theory.
\begin{prop}
The discrete forward Hamilton's equations (\ref{Discrete Hamilton's Equations, one quad, 1})-(\ref{Discrete Hamilton's Equations, one quad, 2}) are multisymplectic, in the sense that for a solution of the discrete forward Hamilton's equations,
\begin{equation*}
\Delta t\ d\varphi_{[a]b+1} \wedge d\pi^1_{[a]b+1}   - \Delta t\ d\varphi_{[a]b} \wedge d\pi^1_{[a]b} + \Delta x\ d\varphi_{a+1[b]} \wedge d\pi^0_{a+1[b]}  - \Delta x\ d\varphi_{a[b]} \wedge d\pi^0_{a[b]}  = 0,
\end{equation*}
evaluated on discrete first variations.
\begin{proof}
In what follows, $H_d^+$ will be evaluated at $(\varphi_{[a]b},\varphi_{a[b]},\pi^1_{[a]b+1},\pi^0_{a+1[b]})$. Compute
\begin{align*}
0 &= d^2H_d^+ = d\Big( D_1H_d^+ d\varphi_{[a]b} + D_2H_d^+ d\varphi_{a[b]} + D_3H_d^+ d\pi^1_{[a]b+1} + D_4H_d^+ d\pi^0_{a+1[b]}  \Big) \\
&= d(D_1H_d^+) \wedge d\varphi_{[a]b} + d(D_2H_d^+)\wedge d\varphi_{a[b]} + d(D_3H_d^+)\wedge d\pi^1_{[a]b+1} + d(D_4H_d^+)\wedge d\pi^0_{a+1[b]}.
\end{align*}
Then, by our definition of discrete first variations, we have 
\begin{align*}
d(D_1H_d^+) &= \Delta t\ d\pi^1_{[a]b},\\
d(D_2H_d^+) &= \Delta x\ d\pi^0_{a[b]},\\
d(D_3H_d^+) &= \Delta t\ d\varphi_{[a]b+1}, \\
d(D_4H_d^+) &= \Delta x\ d\varphi_{a+1[b]}.
\end{align*}
Substituting these expressions into the equation for $d^2 H_d^+$ yields
\begin{align*}
0 &= d(D_1H_d^+) \wedge d\varphi_{[a]b} + d(D_2H_d^+)\wedge d\varphi_{a[b]} + d(D_3H_d^+)\wedge d\pi^1_{[a]b+1} + d(D_4H_d^+)\wedge d\pi^0_{a+1[b]} \\
&= \Delta t\ d\pi^1_{[a]b} \wedge d\varphi_{[a]b} + \Delta x\ d\pi^0_{a[b]}\wedge d\varphi_{a[b]} + \Delta t\ d\varphi_{[a]b+1} \wedge d\pi^1_{[a]b+1} + \Delta x\ d\varphi_{a+1[b]} \wedge d\pi^0_{a+1[b]} \\
&= - \Delta t\ d\varphi_{[a]b} \wedge d\pi^1_{[a]b} - \Delta x\ d\varphi_{a[b]} \wedge d\pi^0_{a[b]} + \Delta t\ d\varphi_{[a]b+1} \wedge d\pi^1_{[a]b+1} + \Delta x\ d\varphi_{a+1[b]} \wedge d\pi^0_{a+1[b]} \\
&= \Delta t\ d\varphi_{[a]b+1} \wedge d\pi^1_{[a]b+1}   - \Delta t\ d\varphi_{[a]b} \wedge d\pi^1_{[a]b} + \Delta x\ d\varphi_{a+1[b]} \wedge d\pi^0_{a+1[b]}  - \Delta x\ d\varphi_{a[b]} \wedge d\pi^0_{a[b]}.
\end{align*}
\end{proof}
\end{prop}

\begin{remark}
Recall that $\omega^\mu = d\varphi \wedge d\pi^\mu$. Observe that if we divide the above discrete multisymplectic form formula by $\Delta t \Delta x$, it is just a first-order finite difference approximation of $\partial_\mu \omega^\mu = 0$.

Furthermore, it is clear that the above equation is precisely quadrature applied to the multisymplectic form formula $\int_{\partial \Box } \omega^\mu|_{(\varphi,\pi)}(\cdot,\cdot) d^nx_\mu = 0$. 

Finally, we note that a discrete notion of multisymplecticity holds in the more general setting described at the beginning of Section \ref{Discrete Hamiltonian Field Theory}. In the more general setting, discrete multisymplecticity is interpreted as $d^2H_d^{\partial\Delta} = 0$ (when evaluated on first variations), which reduces to the ``usual" notion of multisymplecticity in the spacetime tensor product case. 
\end{remark}

\begin{figure}[h]
\begin{center}
\includegraphics[width=90mm]{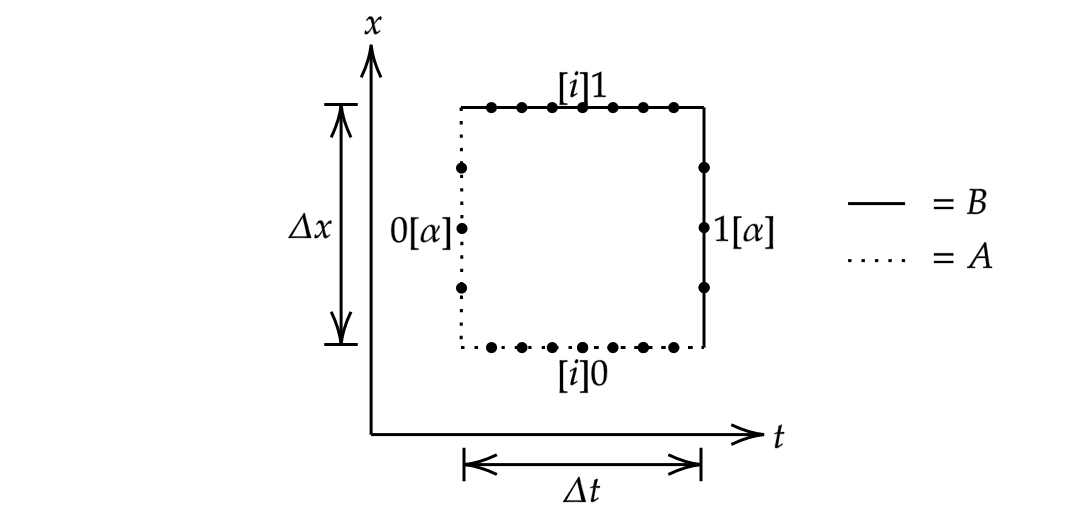}
\caption{Schematic for multiple quadrature points along each edge of $\Box \in \mathcal{T}(X)$.}
\label{multquad}
\end{center}
\end{figure}

\textbf{General Quadrature Approximation.} From here, the generalization to multiple quadrature points is straight-forward. For simplicity, we take the bottom-left vertex of $\Box \in \mathcal{T}(X)$ to be $(0,0)$. Then, $\Box = [0,\Delta t] \times [0,\Delta x]$. In the temporal direction, introduce quadrature points $c_i \in [0,1]$, $i= 1,\dots,s$, and associated quadrature weights $b_i$; we normalize these such that $\sum_ib_i = 1$ (for both $c_i$ and $b_i$, we'll have to explicitly include a factor of $\Delta t$ later) and without loss of generality, we assume each $b_i \neq 0$. Similarly, for the spatial direction, introduce quadrature points $\tilde{c}_\alpha$, $\alpha = 1, \dots, \sigma$ and the associated non-zero weights $\tilde{b}_\alpha$ (normalized as before). Let $\varphi_{[i]0} = \varphi(c_i \Delta t,0)$, $\varphi_{0[\alpha]} = (0,\tilde{c}_\alpha \Delta x)$, $\varphi_{[i]1} = \varphi(c_i \Delta t,\Delta x)$, $\varphi_{1[\alpha]} = (\Delta t, \tilde{c}_\alpha \Delta x)$. Similarly define $\pi^0_{0[\alpha]}$, $\pi^1_{[i]0}$, $\pi^0_{1[\alpha]}$, $\pi^1_{[i]1}$. As before, we take $B$ to be the part of the boundary in the forward direction. See Figure \ref{multquad}.

Then, use quadrature to approximate the boundary integral:
\begin{align*}
\int_B p^\mu \phi d^nx_\mu &= \int_{0}^{\Delta t} (p^1\phi)|_{x = \Delta x} dt + \int_{0}^{\Delta x} (p^0\phi)_{t = \Delta t} dx \\ 
& \approx \sum_{i=1}^s \Delta t\ b_i \pi^1_{[i]1}\varphi_{[i]1} + \sum_{\alpha = 1}^\sigma \Delta x \tilde{b}_\alpha p^0_{1[\alpha]} \varphi_{1[\alpha]} \equiv \sumint_B \pi_B\varphi_B.
\end{align*}
The associated discrete boundary Hamiltonian is
$$ H_d^+(\{\varphi_{[i]0},\varphi_{0[\alpha]},\pi^1_{[i]1},\pi^0_{1[\alpha]}\}_{i,\alpha}) = \ext \Big( \sum_{i=1}^s \Delta t\ b_i \pi^1_{[i]1}\varphi_{[i]1} + \sum_{\alpha = 1}^\sigma \Delta x \tilde{b}_\alpha \pi^0_{1[\alpha]} \varphi_{1[\alpha]} - S_d^{\Box_{ab}}[\phi,p]\Big). $$
\begin{prop}
The discrete forward Hamilton's equations arising from the Type II variational principle are
\begin{subequations}
\begin{alignat}{3}
\label{Discrete Hamilton's Equations, full quad, 1}
\pi^1_{[i]0} &= \frac{1}{b_i \Delta t}D_{1,i}H_d^{+}(\{\varphi_{[j]0},\varphi_{0[\beta]},\pi^1_{[j]1},\pi^0_{1[\beta]}\}_{j,\beta}  ),&\qquad i&=1,\dots,s, \\
\pi^0_{0[\alpha]} &= \frac{1}{\tilde{b}_\alpha \Delta x}D_{2,\alpha}H_d^{+}(\{\varphi_{[j]0},\varphi_{0[\beta]},\pi^1_{[j]1},\pi^0_{1[\beta]}\}_{j,\beta}),& \alpha &=1,\dots,\sigma, \\
\varphi_{[i]1} &= \frac{1}{b_i \Delta t}D_{3,i} H_d^{+}(\{\varphi_{[j]0},\varphi_{0[\beta]},\pi^1_{[j]1},\pi^0_{1[\beta]}\}_{j,\beta}),&  i&=1,\dots,s, \\
\label{Discrete Hamilton's Equations, full quad, 2}
\varphi_{1[\alpha]} &= \frac{1}{\tilde{b}_\alpha \Delta x} D_{4,\alpha} H_d^{+}(\{\varphi_{[j]0},\varphi_{0[\beta]},\pi^1_{[j]1},\pi^0_{1[\beta]}\}_{j,\beta}),& \alpha &= 1,\dots,\sigma,
\end{alignat}
\end{subequations}
where $D_{1,i} \equiv \partial/\partial \varphi_{[i]0}$, $D_{2,\alpha} \equiv \partial/\partial \varphi_{0[\alpha]}$, $D_{3,i} \equiv \partial/\partial \pi^1_{[i]1}$, $D_{4,\alpha} \equiv \partial/\partial \pi^0_{1[\alpha]}$. Furthermore, a solution of the discrete forward Hamilton's equations (\ref{Discrete Hamilton's Equations, full quad, 1})-(\ref{Discrete Hamilton's Equations, full quad, 2}) satisfies the discrete multisymplectic conservation law,
\begin{equation}\label{Discrete MFF}
\sum_{i=1}^s \Delta t\ b_i\Big( d\varphi_{[i]1} \wedge d\pi^1_{[i]1} - d\varphi_{[i]0}\wedge d\pi^1_{[i]0} \Big) + \sum_{\alpha=1}^\sigma \Delta x\ \tilde{b}_\alpha \Big( d\varphi_{1[\alpha]}\wedge d\pi^0_{1[\alpha]} - d\varphi_{0[\alpha]}\wedge d\pi^0_{0[\alpha]} \Big) = 0,
\end{equation}
evaluated on discrete first variations.
\begin{proof}
The proof follows similarly to the case of one quadrature point, Proposition \ref{Discrete Hamilton's Equations Prop}. Namely, the discrete forward Hamilton's equations follow from the Type II variational principle $\delta S_d = 0$ subject to variations of $\varphi$ vanishing along $A(X)$ and variations of $\pi$ vanishing along $B(X)$. The discrete multisymplectic conservation law follows from 
$$d^2H_d^+(\{\varphi_{[j]0},\varphi_{0[\beta]},\pi^1_{[j]1},\pi^0_{1[\beta]}\}_{j,\beta} )=0.$$
\end{proof}
\end{prop}
As in the case of one quadrature point, the discrete multisymplectic conservation law is the given quadrature rule applied to $\int_{\partial \Box }\omega^\mu|_{(\varphi,\pi)}(\cdot,\cdot) d^nx_\mu = 0.$

\begin{remark}
The above discrete forward Hamilton's equations were defined on $\Box = [0,\Delta t] \times [0,\Delta x]$. For $\Box_{ab} = [t_a, t_a + \Delta t] \times [x_b, x_b + \Delta x]$, shift the indices $0,1$ appropriately to $a,a+1$ and $b,b+1$, i.e., $\varphi_{[i]0} \rightarrow \varphi_{[i]b}$, $\varphi_{[i]1} \rightarrow \varphi_{[i]b+1}$, $\varphi_{0[\alpha]} \rightarrow \varphi_{a[\alpha]}$, $\varphi_{1[\alpha]} \rightarrow \varphi_{a+1[\alpha]}$ and similarly for the momenta.
\end{remark}

\textbf{Boundary Conditions and Solution Method.} Recall that the discrete forward Hamilton's equations produce a map $(\varphi_{A(\Box)},\pi_{B(\Box)}) \mapsto (\varphi_{B(\Box},\pi_{A(\Box)})$ for each $\Box \in \mathcal{T}(X)$. However, depending on the boundary conditions that we supply on $\partial X$, the actual realization of these maps may be different (in that the boundary conditions determine the variables in $(\varphi_{A(\Box)},\pi_{B(\Box)}) \mapsto (\varphi_{B(\Box},\pi_{A(\Box)})$ that we implicitly solve for). The key point is that we must specify the field value or the normal momenta along four edges (and the edges may repeat, such as supplying field values and normal momenta on the same edge; see the discussion of evolutionary systems below). This will depend on whether the Hamiltonian PDE we are considering is stationary or evolutionary. 

\begin{figure}[h]
\begin{center}
\includegraphics[width=70mm]{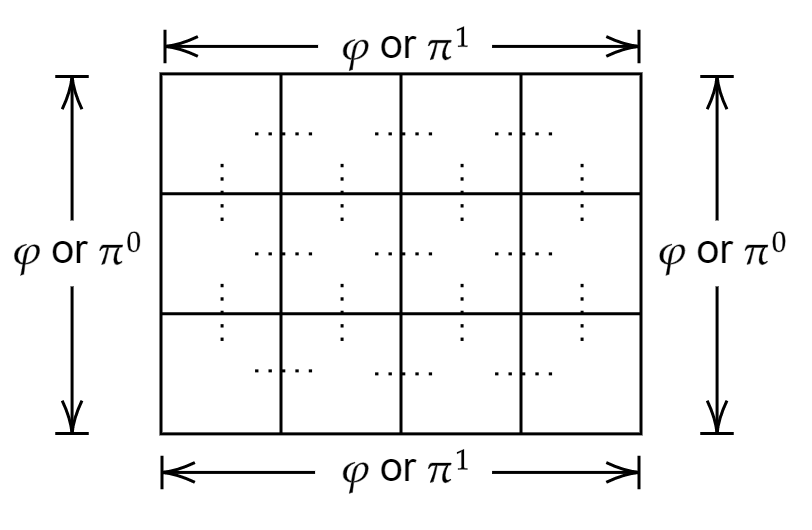}
\caption{Coupling of all of the discrete forward Hamilton's equations for stationary Hamiltonian PDEs; dashed lines along interior edges denote field and normal momenta continuity.}
\label{elliptic}
\end{center}
\end{figure}

Consider a stationary system (e.g., an elliptic system). Then, along $\partial X$, we can specify either Dirichlet boundary conditions, given by the field value $\varphi$, or Neumann boundary conditions, given by the normal momenta value $\pi$. If we supply such boundary conditions, then each $\Box \in \mathcal{T}(X)$ either has two edges with supplied boundary conditions (those on the corners of $X$), has one edge with supplied boundary conditions (those on the edges of $X$), or no supplied boundary conditions (those on the interior). However, the field values and normal momenta values have to be the same along interior edges, which makes up the other required degrees of freedom (recall, we need to specify the field value or normal momenta along four edges). This couples all of the implicit maps $(\varphi_{A(\Box)},\pi_{B(\Box)}) \mapsto (\varphi_{B(\Box},\pi_{A(\Box)})$ together, so that the solution must be solved simultaneously for every $\Box \in \mathcal{T}(X)$. See Figure \ref{elliptic}.

\begin{figure}[h]
\begin{center}
\includegraphics[width=70mm]{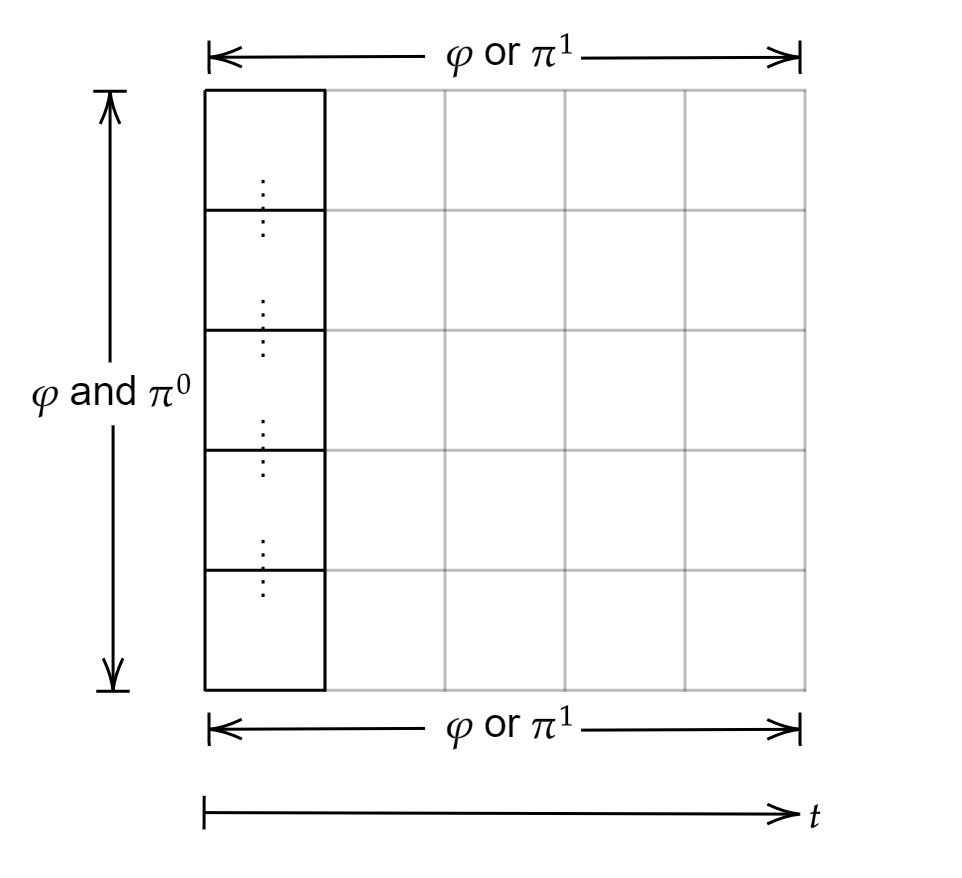}
\caption{Coupling of the discrete forward Hamilton's equations in the same time slice for evolutionary Hamiltonian PDEs; dashed lines along interior edges denote field and normal momenta continuity.}
\label{hyperbolic}
\end{center}
\end{figure}

For an evolutionary system (e.g., a hyperbolic system), we specify the initial conditions at $t=0$, which consist of both the field and normal momenta value ($\pi^0$). On the spatial boundaries, we can either supply Dirichlet or Neumann conditions as above. The continuity of field and normal momenta on the interior edges couples the maps $(\varphi_{A(\Box)},\pi_{B(\Box)}) \mapsto (\varphi_{B(\Box},\pi_{A(\Box)})$ together for each $\Box$ in the same time slice and produces the remaining required degrees of freedom. Hence, one solves these coupled equations on the first time slice which supplies new initial conditions for the subsequent timeslice; one then continues this process recursively for each time step, thereby allowing the discrete solution to be computed in a time marching fashion. See Figure \ref{hyperbolic}.

\begin{remark}\textbf{Solvability.}
It should be noted that the map $(\varphi_A,\pi_B) \mapsto (\varphi_B,\pi_A)$ defined by the discrete forward Hamilton's equations are always well-defined, as can be seen explicitly from the equations (\ref{Discrete Hamilton's Equations, full quad, 1})-(\ref{Discrete Hamilton's Equations, full quad, 2}). This is a property of the (discrete) generating functional and is agnostic to the specific Hamiltonian in question. However, as discussed above, with regard to constructing a numerical method, the implementation of the method in general involves implicitly inverting the relation $(\varphi_A,\pi_B) \mapsto (\varphi_B,\pi_A)$ for the desired variables. For example, if one specifies Neumann boundary conditions on all of $\partial X$ for a stationary system, then the numerical method is given by solving for the map $(\pi_A,\pi_B) \mapsto (\varphi_A,\varphi_B)$ implicitly from the map $(\varphi_A,\pi_B) \mapsto (\varphi_B,\pi_A)$. As another example, for an evolutionary problem, if one specifies Neumann spatial boundary conditions and specifies initial conditions (with both $\varphi$ and $\pi^0$), then the numerical method is given by solving for the map $(\varphi_A,\pi_A) \mapsto (\varphi_B,\pi_B)$ implicitly from the map $(\varphi_A,\pi_B) \mapsto (\varphi_B,\pi_A)$. As these two examples indicate, in general, the form of the map necessary to implement the method is highly dependent on the type of Hamiltonian, as well as the supplied boundary conditions. As such, a discussion of the well-definedness of the implemented map is beyond the scope of this paper, since such a discussion would be highly dependent on the type of problem and boundary conditions, and the functional analytic tools needed in each case would differ drastically. 

We will outline the general argument, although the specifics are left to future work. Note that equations (\ref{Discrete Hamilton's Equations, full quad, 1})-(\ref{Discrete Hamilton's Equations, full quad, 2}) can be written formally as
\begin{subequations}
\begin{align}
\pi_A = \hat{D}_{\varphi_A}H_d^+(\varphi_A,\pi_B), \label{formalDiscreteEqns a} \\
\varphi_B = \hat{D}_{\pi_B}H_d^+(\varphi_A,\pi_B), \label{formalDiscreteEqns b}
\end{align}
\end{subequations}
where $\hat{D}$ denotes the differentiation operators in (\ref{Discrete Hamilton's Equations, full quad, 1})-(\ref{Discrete Hamilton's Equations, full quad, 2}) (and appropriately scaled by the quadrature weights). Showing that one can invert the relations (\ref{formalDiscreteEqns a})-(\ref{formalDiscreteEqns b}) for the implemented map would then rest on an implicit function theorem type argument, for a sufficiently small $\Box \subset X$. The derivatives of the equations (\ref{formalDiscreteEqns a})-(\ref{formalDiscreteEqns b}) would then involve second derivatives of $H_d^+$, so hyperregularity would prove crucial in such a proof. For degenerate Hamiltonians, some form of constraint or gauge-fixing would be necessary to complete the proof. We aim to explore issues dealing with solvability in future work, as well as related issues such as error analysis, which is again highly dependent on the specific class of Hamiltonians and boundary conditions considered.
\end{remark}

\textbf{Discrete Noether's Theorem.} 
In the continuum theory, we saw that for a vertical group action on the restricted dual jet bundle which leaves the action invariant, there is an associated Noether conservation law (\ref{Noether's Theorem coordinates}) for solutions of Hamilton's equations.

In the discrete setting, suppose there is a differentiable and vertical $G$  action on the discrete restricted dual jet bundle $\{t_i,x_j\} \times Q \times (Q^*)^2$ (relative to $\Box \in \mathcal{T}(X)$) which leaves invariant the generalized discrete Lagrangian
\begin{align*}
R^{\Box}_d(\varphi_{A(\Box)},\varphi_{B(\Box)},\pi_{B(\Box)}) &= \sumint_{B(\Box)} \pi_{B(\Box)}\varphi_{B(\Box)} - H_d^+(\varphi_{A(\Box)},\pi_{B(\Box)}) \\ &= \sum_{i=1}^s \Delta t\ b_i \pi^1_{[i]1}\varphi_{[i]1} + \sum_{\alpha = 1}^\sigma \Delta x \tilde{b}_\alpha p^0_{1[\alpha]} \varphi_{1[\alpha]} - H_d^+(\{\varphi_{[i]0},\varphi_{0[\alpha]},\pi^1_{[i]1},\pi^0_{1[\alpha]}\}_{i,\alpha}).
\end{align*}
\begin{prop}
If the generalized discrete Lagrangian is invariant under a differentiable and vertical $G$ action on the discrete restricted dual jet bundle, then a solution of the discrete forward Hamilton's equations (\ref{Discrete Hamilton's Equations, full quad, 1})-(\ref{Discrete Hamilton's Equations, full quad, 2})  admits a discrete analog of Noether's theorem:
\begin{align}\label{Discrete Noether's Theorem}
\sum_i \Delta t &\ b_i \pi^1_{[i]1} i_{\tilde{\xi}}d\varphi_{[i]1} + \sum_\alpha \Delta x\ \tilde{b}_\alpha \pi^0_{1[\alpha]} i_{\tilde{\xi}}d\varphi_{1[\alpha]} \\
\nonumber &- \sum_i \Delta t\ b_i \pi^1_{[i]0} i_{\tilde{\xi}}d\varphi_{[i]0} - \sum_\alpha \Delta x\ \tilde{b}_\alpha \pi^0_{0[\alpha]} i_{\tilde{\xi}}d\varphi_{0[\alpha]} = 0,
\end{align}
where $\tilde{\xi}$ is the infinitesimal generator associated with $\xi \in \mathfrak{g}$.
\begin{proof}
For brevity, we will omit the arguments of $R^{\Box}_d$ and $H_d^+$ (refer to the definition of $R^{\Box}_d$ above). Since the generalized discrete Lagrangian is invariant under the $G$ action, that means that the directional derivative in the direction of the infinitesimal generator vanishes,
\begin{align*}
0 &= dR_d^{\Box} \cdot \tilde{\xi}\\
&= \sum_i \Delta t\ b_i i_{\tilde{\xi}}d(\pi^1_{[i]1}\varphi_{[i]1}) + \sum_\alpha \Delta x\ \tilde{b}_\alpha i_{\tilde{\xi}}d(\pi^0_{1[\alpha]} \varphi_{1[\alpha]}) \\
&\qquad - \sum_i \Big(D_{1,i} H_d^+ i_{\tilde{\xi}} d\varphi_{[i]0} + D_{3,i} H_d^+ i_{\tilde{\xi}} d\pi^1_{[i]1}\Big) - \sum_\alpha \Big(D_{2,\alpha} H_d^+ i_{\tilde{\xi}} d\varphi_{0[\alpha]} + D_{4,\alpha} H_d^+ i_{\tilde{\xi}} d\pi^0_{1[\alpha]} \Big) \\
&= \sum_i \Delta t\ b_i (\cancelto{(1)}{i_{\tilde{\xi}} d\pi^1_{[i]1}\varphi_{[i]1}} + \pi^1_{[i]1} i_{\tilde{\xi}}d\varphi_{[i]1}) + \sum_\alpha \Delta x\ \tilde{b}_\alpha (\cancelto{(2)}{i_{\tilde{\xi}}d\pi^0_{1[\alpha]} \varphi_{1[\alpha]}} + \pi^0_{1[\alpha]} i_{\tilde{\xi}}d\varphi_{1[\alpha]}) \\
&\qquad- \sum_i \Delta t\ b_i \Big(\pi^1_{[i]0} i_{\tilde{\xi}} d\varphi_{[i]0} + \cancelto{(1)}{\varphi_{[i]1} i_{\tilde{\xi}} d\pi^1_{[i]1}\Big)} - \sum_\alpha \Delta x\ \tilde{b}_\alpha \Big(\pi^0_{0[\alpha]} i_{\tilde{\xi}} d\varphi_{0[\alpha]} + \cancelto{(2)}{\varphi_{1[\alpha]} i_{\tilde{\xi}} d\pi^0_{1[\alpha]}} \Big).
\end{align*}
\end{proof}
\end{prop}
\begin{remark}
Note that the above looks like quadrature applied to the continuous Noether's theorem, 
$$\int_{\partial \Box} p^\mu (i_{\tilde{\xi}}d\phi)d^nx_\mu = 0$$
(with the caveat that, in the continuum case, $G$ acts on the restricted dual jet bundle, whereas in the discrete case, $G$ acts on the discrete restricted dual jet bundle). One can obtain such a $G$-invariant $R_d$ via $G$-equivariant interpolation (see \citet{LeZh2009} and \citet{Le2019}), in which case, the discrete Noether theorem is precisely quadrature applied to Noether's theorem.

Also, note that a discrete Noether's theorem holds in the more general setting described at the beginning of Section \ref{Discrete Hamiltonian Field Theory}. In the more general setting, the discrete Noether's theorem is interpreted as $d R_d^{\Delta}\cdot \tilde{\xi} = 0$ (for a $G$-invariant generalized discrete Lagrangian), which reduces to the ``usual" coordinate notion of the discrete Noether's theorem, equation (\ref{Discrete Noether's Theorem}), in the spacetime tensor product case. 
\end{remark}

\begin{remark}
Another way to interpret this discrete Noether's theorem is to view the map determined by the discrete forward Hamilton's equations, $(\varphi_{A(\Box)},\pi_{B(\Box)}) \mapsto (\varphi_{B(\Box)}, \pi_{A(\Box)})$, as implicitly defining a  forward map $F_{H^+_d}: (\varphi_{A(\Box)},\pi_{A(\Box)}) \mapsto (\varphi_{B(\Box)},\pi_{B(\Box)})$. For some subset $S$ of $\partial \Box$, define the discrete (Hamiltonian) Cartan form (at a solution of the discrete forward Hamilton's equations)
\begin{equation}\label{Discrete Cartan Form}
\Theta_d^{S} = \sum_{(t_k,x_l) \in S} \beta_{kl} \pi^n_{kl} d\varphi_{kl},
\end{equation}
where $\pi^n$ denotes the normal component of the momenta and $\beta^{kl}$ denotes the quadrature weight at $(t_k,x_l) \in S$ (which equals $\Delta t\ b_i$ for the $i^{th}$ node of $S$ along fixed $x$ and equals $\Delta x\ \tilde{b}_\alpha$ for the $\alpha^{th}$ node of $S$ along fixed $t$). Such a discrete Cartan form involving summing over nodes corresponding to boundary variations was introduced by \citet{MaPaSh1998} in the Lagrangian framework; in the discrete Hamiltonian setting which we constructed, (\ref{Discrete Cartan Form}) is the appropriate definition since $\Theta^{\partial \Box}_d$ precisely encodes such discrete boundary variations.

Then, the discrete Noether theorem (\ref{Discrete Noether's Theorem}) can be expressed as
$$ F_{H_d^+}^* (\Theta^{B(\Box)}_d ) \cdot \tilde{\xi} = \Theta^{A(\Box)}_d \cdot \tilde{\xi}. $$

Note also that the discrete multisymplectic form formula (\ref{Discrete MFF}) can be expressed as
$$ d\Theta_d^{\partial \Box}(\cdot,\cdot) = 0, $$
when evaluated on discrete first variations.
\end{remark}

\subsection{Galerkin Hamiltonian Variational Integrators}\label{GHVI}
The missing ingredient in our construction of a variational integrator is the discrete approximation of the action over $\Box \in \mathcal{T}(X)$, $S_d^{\Box}[\phi,p]$. We will extend the construction of Galerkin Hamiltonian variational integrators, introduced in \citet{LeZh2009} for Hamiltonian ODEs, to the case of Hamiltonian PDEs.

\begin{remark}
To be definitive, we will assume that the space(time) $X$ has the Euclidean metric. The discussion below is equally valid for the Minkowski metric, except one has to include the appropriate minus signs throughout.
\end{remark}

Consider for simplicity $[0,\Delta t]\times [0,\Delta x] = \Box \in \mathcal{T}(X)$. Fix quadrature rules in the temporal direction (weights $b_i$ and nodes $c_i$, $i=1,\dots,s$) and spatial direction (weights $\tilde{b}_\alpha$ and nodes $\tilde{c}_\alpha$, $\alpha = 1, \dots, \sigma$) as before. Note the action $S[\phi,p] = \int (p^\mu \partial_\mu\phi - H(\phi,p^0,p^1) )d^2x$ involves the fields $\phi,$ their derivatives $\partial_\mu\phi$, and the multimomenta $p^\mu$ ($\mu = 0,1$). For the field and their derivatives, we could either approximate the field using a finite-dimensional subspace and subsequently take derivatives; or conversely, approximate the derivatives and subsequently integrate to obtain the values of the field. We will take the latter approach (we will extremize over the internal stages at the end, so the two approaches are equivalent). Introduce basis functions $\{\chi_i(\tau)\}_{i=1}^s$, $\tau \in [0,1]$, for an $s$-dimensional function space and similarly $\{\tilde{\chi}_\alpha(\tau)\}_{\alpha=1}^\sigma$ for a $\sigma$-dimensional function space. We will use the tensor product basis $\{\chi_i(\tau \Delta t) \tilde{\chi}_\alpha(\rho \Delta x) \}_{i,\alpha}$ to discretize the derivatives of the field. Approximate the derivatives as
\begin{subequations}
\begin{align} \label{derivativeapprox1}
\partial_t\phi_d (\tau \Delta t, \rho \Delta x) = \sum_{i,\alpha} V^{i\alpha} \chi_i(\tau) \tilde{\chi}_\alpha(\rho), \\
\partial_x\phi_d (\tau \Delta t, \rho \Delta x) = \sum_{i,\alpha} W^{i\alpha} \chi_i(\tau) \tilde{\chi}_\alpha(\rho). \label{derivativeapprox2}
\end{align}
\end{subequations}
We can integrate in time or space to determine the field values. In particular, the internal stages are given by the field values at the nodes $(c_i\Delta t, \tilde{c}_\alpha \Delta x$):
\begin{align*}
\Phi_{i\alpha} &\equiv \phi(c_i\Delta t, \tilde{c}_\alpha \Delta x) = \phi(0, \tilde{c}_\alpha \Delta x) + \Delta t \sum_{j,\beta} V^{j\beta} \int_0^{c_i} \chi_j(s)ds\ \tilde{\chi}_\beta (\tilde{c}_\alpha) = \varphi_{0[\alpha]} + \Delta t \sum_{j,\beta} A_{i\alpha,j\beta}V^{j\beta}, \\
\Phi_{i\alpha} &\equiv \phi(c_i\Delta t, \tilde{c}_\alpha \Delta x) = \phi(c_i\Delta t, 0) + \Delta x \sum_{j,\beta} W^{j\beta} \chi_j(c_i) \int_0^{c_\alpha} \tilde{\chi}_\beta(s)ds = \varphi_{[i]0} + \Delta x \sum_{j,\beta} \tilde{A}_{i\alpha,j\beta}W^{j\beta},
\end{align*}
where $A_{i\alpha,j\beta} = \int_0^{c_i}\chi_j(s)ds\ \tilde{\chi}_\beta(\tilde{c}_\alpha)$ and $\tilde{A}_{i\alpha,j\beta} = \chi_j(c_i) \int_0^{c_\alpha} \tilde{\chi}_\beta(s)ds$. Note that $\Phi_{i\alpha}$ must of course be single-valued, so we have a relation between the two above equations:
\begin{equation}\label{Dependent Internal Stages}
\phi_{0[\alpha]} + \Delta t A_{i\alpha,j\beta}V^{j\beta} = \Phi_{i\alpha} = \phi_{[i]0} + \Delta x \tilde{A}_{i\alpha,j\beta}W^{j\beta}.
\end{equation}
We expect such a relation since extremizing over $\Phi_{i\alpha}$ is equivalent to extremizing over $V_{i\alpha}$ or $W_{i\alpha}$ (but not both; however, we will relax this assumption in the subsequent discussion).

Integrating to $1$ gives the unknown field boundary values,
\begin{align*}
\varphi_{1[\alpha]} &\equiv \phi(\Delta t, \tilde{c}_\alpha \Delta x) = \phi(0, \tilde{c}_\alpha \Delta x) + \Delta t \sum_{j,\beta} V^{j\beta} \int_0^{1} \chi_j(s)ds\ \tilde{\chi}_\beta (\tilde{c}_\alpha) = \varphi_{0[\alpha]} + \Delta t \sum_{j,\beta} B_{\alpha,j\beta}V^{j\beta}, \\
\varphi_{[i]1} &\equiv \phi(c_i\Delta t, \Delta x) = \phi(c_i\Delta t, 0) + \Delta x \sum_{j,\beta} W^{j\beta} \chi_j(c_i) \int_0^{1} \tilde{\chi}_\beta(s)ds = \varphi_{[i]0} + \Delta x \sum_{j,\beta} \tilde{B}_{i,j\beta}W^{j\beta},
\end{align*}
where $B_{\alpha,j\beta} = \int_0^{1}\chi_j(s)ds\ \tilde{\chi}_\beta(\tilde{c}_\alpha)$ and $\tilde{B}_{i,j\beta} = \chi_j(c_i) \int_0^{1} \tilde{\chi}_\beta(s)ds$.

We define the internal stages for the momenta $P^0_{i\alpha} = p^0(c_i\Delta t, \tilde{c}_\alpha \Delta x), P^1_{i\alpha} = p^1(c_i\Delta t, \tilde{c}_\alpha \Delta x)$. Unlike the field internal stage expansions, one does not need to introduce an approximating function space for the momenta internal stages, since the action only involves derivatives of the field and not the momenta. At this point, we could work directly with these internal stages; however, we will expand the momenta similarly to the fields,
\begin{align*}
P^0_{i\alpha} &= \pi^0_{1[\alpha]} - \Delta t \sum_{j,\beta}A'_{i\alpha,j\beta}X^{j\beta}, \\
P^1_{i\alpha} &= \pi^1_{[i]1} - \Delta x \sum_{j,\beta}\tilde{A}'_{i\alpha,j\beta} Y^{j\beta},
\end{align*}
where $A'_{i\alpha,j\beta}$ and $\tilde{A}'_{i\alpha,j\beta}$ are arbitrary expansion coefficients and $X^{j\beta}, Y^{j\beta}$ are internal variables representing $\partial_0p^0$ and $\partial_1p^1$ respectively. The unknown momenta boundary values are similarly defined as
\begin{align*}
\pi^0_{0[\alpha]} &= \pi^0_{1[\alpha]} - \Delta t \sum_{j,\beta}B'_{\alpha,j\beta}X^{j\beta}, \\
\pi^1_{[i]0} &= \pi^1_{[i]1} - \Delta x \sum_{j,\beta}\tilde{B}'_{i,j\beta} Y^{j\beta},
\end{align*}
where $B'_{\alpha,j\beta}$ and $\tilde{B}'_{i,j\beta}$ are again arbitrary expansion coefficients. We will see later that the expansion coefficients will have to satisfy symplecticity conditions in order for the method to be well-defined.

We then approximate the action integral $S[\phi,p] = \int (p^\mu \partial_\mu\phi - H(\phi,p^0,p^1))d^2x$ using quadrature and the above internal stages
$$ S_d^{\Box}[\Phi_{i\alpha},P_{i\alpha}] = \Delta t \Delta x \sum_{i,\alpha} b_i \tilde{b}_\alpha \Big(P^0_{i\alpha}\partial_t\phi_d (c_i\Delta t, \tilde{c}_\alpha \Delta x) + P^1_{i\alpha}\partial_x\phi_d (c_i\Delta t, \tilde{c}_\alpha \Delta x) - H(\Phi_{i\alpha},P^0_{i\alpha},P^1_{i\alpha}) \Big). $$
The discrete boundary Hamiltonian is obtained by extremizing over the internal stages $\Phi,P^0,P^1$, which are defined in terms of $V,X,Y$. Since we have already enforced the boundary conditions in the above field and momenta expansions, we can construct the discrete boundary Hamiltonian by extremizing over $V^{i\alpha}$, $X^{i\alpha}$, $Y^{i\alpha}$ (for every $i=1,\dots,s$ and $\alpha = 1,\dots,\sigma$),
$$ H_d^+(\{\varphi_{[i]0},\varphi_{0[\alpha]},\pi^1_{[i]1},\pi^0_{1[\alpha]}\}_{i,\alpha}) = \underset{V^{i\alpha},X^{i\alpha},Y^{i\alpha}}{\ext } \Big( \underbrace{ \sum_{i=1}^s \Delta t\ b_i \pi^1_{[i]1}\varphi_{[i]1} + \sum_{\alpha = 1}^\sigma \Delta x\ \tilde{b}_\alpha \pi^0_{1[\alpha]} \varphi_{1[\alpha]} - S_d^{\Box}[\Phi_{i\alpha},P_{i\alpha}]}_{\equiv K(\{\varphi_A,\pi_B, V^{i\alpha},X^{i\alpha},Y^{i\alpha}\})}\Big). $$
$H_d^+$ is then given by extremizing $K(\{\varphi_A,\pi_B, V^{i\alpha}, X^{i\alpha},Y^{i\alpha}\})$ with respect to $V^{i\alpha}, X^{i\alpha}$, and $Y^{i\alpha}$ (where again we denote $\varphi_A = \{\varphi_{[i]0},\varphi_{0[\alpha]}\}$ and $\pi_B = \{\pi^1_{[i]1},\pi^0_{1[\alpha]}\}$). Expanding $K$, we have
\begin{align*}
&K(\{\varphi_A,\pi_B, V^{i\alpha}, X^{i\alpha},Y^{i\alpha}\})\\
&\quad= \Delta t \sum_i b_i \pi^1_{[i]1}(\varphi_{[i]0} + \Delta x \sum_{j,\beta} \tilde{B}_{i,j\beta}W^{j\beta}) + \Delta x \sum_\alpha \tilde{b}_\alpha \pi^0_{1[\alpha]}(\varphi_{0[\alpha]} + \Delta t \sum_{j,\beta} B_{\alpha,j\beta}V^{j\beta}) \\
&\qquad- \Delta t \Delta x \sum_{i,\alpha} b_i \tilde{b}_\alpha \Big(\pi^0_{1[\alpha]} - \Delta t \sum_{k,\gamma}A'_{i\alpha,k\gamma}X^{k\gamma} \Big) \sum_{j,\beta}V^{j\beta}\chi_j(c_i)\tilde{\chi}_\alpha(\tilde{c}_\alpha) \\
&\qquad-  \Delta t \Delta x \sum_{i,\alpha} b_i\tilde{b}_\alpha \Big( \pi^1_{[i]1} - \Delta x \sum_{k,\gamma}\tilde{A}'_{i\alpha,k\gamma} Y^{k\gamma} \Big) \sum_{j,\beta}W^{j\beta}\chi_j(c_i)\tilde{\chi}_\alpha(\tilde{c}_\alpha) \\
&\qquad+ \Delta t \Delta x \sum_{i,\alpha}b_i\tilde{b}_\alpha H(\Phi_{i\alpha},P^0_{i\alpha},P^1_{i\alpha}).
\end{align*}
The stationarity conditions $\partial K/\partial V^{i\alpha} = 0$, $\partial K/\partial X^{i\alpha} = 0$, $\partial K/\partial Y^{i\alpha} = 0$, combined with the discrete forward Hamilton's equations (\ref{Discrete Hamilton's Equations, full quad, 1})-(\ref{Discrete Hamilton's Equations, full quad, 2}) define our multisymplectic variational integrator. 

Supposing that one solves the stationarity conditions for $V^{i\alpha}$, $X^{i\alpha}$, $Y^{i\alpha}$ in terms of $\varphi_A$ and $\pi_B$, this gives $H_d^+(\{\varphi_A,\pi_B\}) = K(\{\varphi_A, \pi_B, V^{i\alpha}(\varphi_A,\pi_B), X^{i\alpha}(\varphi_A,\pi_B),Y^{i\alpha}(\varphi_A,\pi_B) \}).$ The right hand side of the discrete forward Hamilton's equations, (\ref{Discrete Hamilton's Equations, full quad, 1})-(\ref{Discrete Hamilton's Equations, full quad, 2}), can then be computed in terms of $K$ via
\begin{align*}
\frac{\partial}{\partial \varphi_{[i]0}} H_d^+(\{\varphi_A,\pi_B\}) &= \frac{\partial}{\partial \varphi_{[i]0}} K(\{\varphi_A, \pi_B, V^{i\alpha}(\varphi_A,\pi_B), X^{i\alpha}(\varphi_A,\pi_B),Y^{i\alpha}(\varphi_A,\pi_B) \}) \\
&= \frac{\partial}{\partial \varphi_{[i]0}} K + \sum_{j,\alpha} \Big( \cancel{\frac{\partial K}{\partial V^{j\alpha}}}\frac{\partial V^{j\alpha}}{\partial \varphi_{[i]0}} + \cancel{\frac{\partial K}{\partial X^{j\alpha}}}\frac{\partial X^{j\alpha}}{\partial \varphi_{[i]0}} + \cancel{\frac{\partial K}{\partial Y^{j\alpha}}}\frac{\partial Y^{j\alpha}}{\partial \varphi_{[i]0}} \Big) \\
&= \frac{\partial}{\partial \varphi_{[i]0}} K,
\end{align*}
and similarly for the other specified boundary values. Hence, the derivatives of $H_d^+$ with respect to $\varphi_A$, $\pi_B$ can be computed using only the explicit dependence of $K$ on $\varphi_A, \pi_B$.

\subsection{Multisymplectic Partitioned Runge--Kutta Method}\label{MPRK subsection}
Let us suppose that instead of the basis $\{\chi_i\},\{\tilde{\chi}_\alpha\}$, we choose basis functions $\{\psi_i\}$, $\{\tilde{\psi}_\alpha\}$ that have the interpolating property $\psi_i(c_j) = \delta_{ij}$, $\tilde{\psi}_\alpha(\tilde{c}_\beta) = \delta_{\alpha\beta}$. Note that one can always transform the previous set of basis functions to a set of basis functions with this property, assuming that the original choice of basis functions $\chi_i, \tilde{\chi}_\alpha$ have the property that the matrices with entries $M_{ij} = \chi_i(c_j)$, $\tilde{M}_{\alpha\beta} = \tilde{\chi}_\alpha(\tilde{c}_\beta)$ are invertible. If they are not, then the expansion of the derivatives, equations (\ref{derivativeapprox1})-(\ref{derivativeapprox2}), does not depend independently on all of the $V^{i\alpha}$, $W^{i\alpha}$ and hence one needs to reduce the number of independent variables; to avoid this, ensure that the matrices with entries $\chi_i(c_j)$ and $\tilde{\chi}_\alpha(\tilde{c}_\beta$) are invertible. Letting $\chi(\cdot) = (\chi_1(\cdot),\dots,\chi_s(\cdot))^T$ and $\tilde{\chi}(\cdot) = (\tilde{\chi}_1(\cdot),\dots,\tilde{\chi}_\sigma(\cdot))^T$ (and similarly define $\psi$, $\tilde{\psi}$), a set of basis functions with the interpolating property can be constructed by $\psi = M^{-1}\chi$, $\tilde{\psi} = \tilde{M}^{-1} \tilde{\chi}$. In particular, the $\{\psi_i\}$, $\{\tilde{\psi}_\alpha\}$ span the same function spaces as the $\{\chi_i\}$, $\{\tilde{\chi}_\alpha\}$ respectively, so there is no loss of generality.

With this assumption, we approximate the derivatives of the fields  as
\begin{align*}
\partial_t\phi_d(c_i\Delta t, \tilde{c}_\alpha \Delta x) &= \sum_{j,\beta} V^{j\beta}\psi_j(c_i)\tilde{\psi}_\beta(\tilde{c}_\alpha) = V^{i\alpha}, \\
\partial_x\phi_d(c_i\Delta t, \tilde{c}_\alpha \Delta x) &= \sum_{j,\beta} W^{j\beta}\psi_j(c_i)\tilde{\psi}_\beta(\tilde{c}_\alpha) = W^{i\alpha}.
\end{align*}
Integrating gives the internal stages and the unknown boundary values,
\begin{align*}
\Phi_{i\alpha} &= \varphi_{0[\alpha]} + \Delta t \sum_{j}a_{ij}V^{j\alpha},\\
\Phi_{i\alpha} &= \varphi_{[i]0} + \Delta x \sum_\beta \tilde{a}_{\alpha\beta} W^{i\beta} \\
\varphi_{1[\alpha]} &= \varphi_{0[\alpha]} + \Delta t \sum_j b_j V^{j\alpha},\\
\varphi_{[i]1} &= \varphi_{[i]0} + \Delta x \sum_\beta \tilde{b}_\beta W^{i\beta},
\end{align*}
where $a_{ij} = \int_0^{c_i}\psi_j(s)ds$, $\tilde{a}_{\alpha\beta} = \int_0^{\tilde{c}_\alpha} \tilde{\psi}_\beta(s)ds$ and the quadrature weights $b_i = \int_0^1 \psi_i(s)ds$, $\tilde{b}_\alpha = \int_0^1 \tilde{\psi}_\alpha(s)ds$ are chosen so that quadrature is exact on the span of the basis functions. As before, expand the momenta using a different set of coefficients. 
\begin{align*}
X^{i\alpha} &= \partial_tp^0_d(c_i\Delta t, \tilde{c}_\alpha \Delta x),\\
Y^{i\alpha} &= \partial_xp^1_d(c_i\Delta t, \tilde{c}_\alpha \Delta x), \\
P^0_{i\alpha}&= \pi^0_{1[\alpha]} - \Delta t \sum_j a'_{ij}X^{j\alpha}, \\
P^1_{i\alpha} &= \pi^1_{[i]1} - \Delta x \sum_\beta \tilde{a}'_{\alpha\beta}Y^{i\beta}, \\
\pi^0_{0[\alpha]} &= \pi^0_{1[\alpha]} - \Delta t \sum_j b_j' X^{j\alpha},\\
\pi^1_{[i]0} &= \pi^1_{[i]1} - \Delta x \sum_\beta \tilde{b}'_\beta Y^{i\beta}.
\end{align*}
We impose that $b_j'>0$, $\tilde{b}_\beta'>0$ and that $\sum_j b_j' = 1$, $\sum_\beta \tilde{b}_\beta' = 1$ for the approximation to be consistent. We will later derive a condition on the coefficients $a'_{ij}$, $\tilde{a}'_{\alpha\beta}$, $b'_i$, $\tilde{b}'_\alpha $ in order for the method to be well-defined. For now, we proceed formally.

With these, $K$ can be expressed as
\begin{align*}
&K(\{\varphi_A,\pi_B, V^{i\alpha}, X^{i\alpha},Y^{i\alpha}\})\\
&\quad= \Delta x \sum_\alpha \tilde{b}_\alpha \pi^0_{1[\alpha]}(\varphi_{0[\alpha]} + \Delta t \sum_{j} b_j V^{j\alpha}) + \Delta t \sum_i b_i \pi^1_{[i]1}(\varphi_{[i]0} + \Delta x \sum_\beta \tilde{b}_\beta W^{i\beta}) \\
&\qquad- \Delta t \Delta x \sum_{i,\alpha} b_i \tilde{b}_\alpha \Big(\pi^0_{1[\alpha]} - \Delta t \sum_j a'_{ij}X^{j\alpha} \Big) V^{i\alpha} \\
&\qquad-  \Delta t \Delta x \sum_{i,\alpha} b_i\tilde{b}_\alpha \Big( \pi^1_{[i]1} - \Delta x \sum_{\beta}\tilde{a}'_{\alpha\beta} Y^{i\beta} \Big) W^{i\alpha}\\
&\qquad+ \Delta t \Delta x \sum_{i,\alpha}b_i\tilde{b}_\alpha H(\Phi_{i\alpha},P^0_{i\alpha},P^1_{i\alpha}).
\end{align*}
Now, we compute the stationarity conditions. First, note that $V$ and $W$ are not independent, since they are related by
$$ \varphi_{0[\alpha]} + \Delta t \sum_{j}a_{ij}V^{j\alpha} = \Phi_{i\alpha} = \varphi_{[i]0} + \Delta x \sum_\beta \tilde{a}_{\alpha\beta} W^{i\beta}; $$
then, taking the derivative with respect to $V^{j\beta}$,
$$ \Delta t a_{ij} \delta_{\alpha\beta} = \Delta x \sum_\gamma \tilde{a}_{\alpha \gamma} \frac{\partial W^{i\gamma}}{\partial V^{j\beta}}.  $$
Let us assume that the Runge--Kutta matrices $(a_{ij})$ and $(\tilde{a}_{\alpha\gamma})$ are invertible (however, in the subsequent section, we will show how to derive the stationarity conditions without this assumption using independent internal stages). Then, the above relation can be inverted to give
$$ \frac{\partial W^{i\sigma}}{\partial V^{j\beta}} = \frac{\Delta t}{\Delta x} a_{ij} (\tilde{a}^{-1})_{\sigma\beta}. $$

Extremizing $K$ with respect to $X^{j\alpha},$
\begin{align*}
0 &= \frac{\partial K}{\partial X^{j\alpha}} = \Delta t^2 \Delta x \sum_i b_i \tilde{b}_\alpha a'_{ij} V^{i\alpha} - \Delta t^2 \Delta x \sum_i b_i \tilde{b}_\alpha a'_{ij} \frac{\partial H}{\partial p^0}(\Phi_{i\alpha}, P^0_{i\alpha}, P^1_{i\alpha}).
\end{align*}
Dividing by $\Delta t^2 \Delta x \tilde{b}_\alpha$ gives
$$ \sum_i b_i a'_{ij} \Big( V^{i\alpha} -  \frac{\partial H}{\partial p^0}(\Phi_{i\alpha},P^0_{i\alpha}, P^1_{i\alpha}) \Big) = 0. $$
Similarly, extremizing $K$ with respect to $Y^{j\alpha}$ gives
$$ \sum_\beta \tilde{b}_\beta \tilde{a}'_{\beta\alpha} \Big(W^{j\beta} - \frac{\partial H}{\partial p^1}(\Phi_{j\beta}, P^0_{j\beta}, P^1_{j\beta}) \Big) = 0. $$
%Note that the matrix with $ji$ entry $b_ia_{ij}$ is invertible since $(a_{ij})$ is (the transpose of this matrix is obtained by multiplying the $i^{th}$ row of $(a_{ij})$ by a non-zero scalar $b_i$, so the rows are of course still linearly independent) and similarly for the matrix with $\alpha\beta$ entry $\tilde{b}_\beta \tilde{a}_{\beta\alpha}$. Hence, we have
% $$ V^{i\alpha} = \frac{\partial H}{\partial p^0}(\Phi_{i\alpha},P^0_{i\alpha}, P^1_{i\alpha}), \ \ W^{i\alpha} = \frac{\partial H}{\partial p^1}(\Phi_{i\alpha},P^0_{i\alpha}, P^1_{i\alpha}). $$ 
These are respectively the internal stage approximations to the De Donder--Weyl equations $\partial_t\phi = \partial H/\partial p^0$ and $\partial_x\phi = \partial H/\partial p^1$.

Extremizing $K$ with respect to $V^{j\beta}$,
\begin{align*}
0 = \frac{\partial K}{\partial V^{j\beta}} &= \Delta t \Delta x b_j \tilde{b}_\beta \pi^0_{1[\beta]} + \Delta t \Delta x \sum_{i,\sigma} b_i\tilde{b}_\sigma \pi^1_{[i]1} \frac{\Delta t}{\Delta x} a_{ij} (\tilde{a}^{-1})_{\sigma\beta} - \Delta t \Delta x b_j \tilde{b}_\beta P^0_{j\beta} \\
&\qquad- \Delta t \Delta x \sum_{i,\sigma} b_i \tilde{b}_\sigma P^1_{i\sigma} \frac{\Delta t}{\Delta x} a_{ij}(\tilde{a}^{-1})_{\sigma\beta} + \Delta t \Delta x \sum_i b_i \tilde{b}_\beta \Delta t a_{ij} \frac{\partial H}{\partial \phi}(\Phi_{i\beta}, P^0_{i\beta}, P^1_{i\beta}).
\end{align*}
Dividing by $\Delta t^2 \Delta x$ and grouping gives
\begin{align*}
b_j \tilde{b}_\beta \frac{\pi^0_{1[\beta]} - P^0_{j\beta}}{\Delta t} + \sum_{i,\sigma} b_i \tilde{b}_\sigma a_{ij}(\tilde{a}^{-1})_{\sigma\beta} \frac{\pi^1_{[i]1}-P^1_{i\sigma}}{\Delta x} = - \sum_i b_i \tilde{b}_\beta a_{ij} \frac{\partial H}{\partial \phi}(\Phi_{i\beta}, P^0_{i\beta},P^1_{i\beta}).
\end{align*}
Substitute $\frac{\pi^0_{1[\beta]} - P^0_{j\beta}}{\Delta t} = \sum_k a'_{jk}X^{k\beta} $ and $\frac{\pi^1_{[i]1}-P^1_{i\sigma}}{\Delta x} = \sum_\gamma \tilde{a}'_{\sigma\gamma}Y^{i\gamma}, $
\begin{align*}
\sum_k b_j \tilde{b}_\beta a'_{jk}X^{k\beta} + \sum_{i,\sigma,\gamma} b_i \tilde{b}_\sigma a_{ij}(\tilde{a}^{-1})_{\sigma\beta} \tilde{a}'_{\sigma\gamma}Y^{i\gamma} = - \sum_i b_i \tilde{b}_\beta a_{ij} \frac{\partial H}{\partial \phi}(\Phi_{i\beta}, P^0_{i\beta},P^1_{i\beta}).
\end{align*}
To symmetrize the above equations, multiply by $\tilde{a}_{\beta\delta}$ and sum over $\beta$, which yields
$$ \sum_{k,\beta} b_j \tilde{b}_\beta a'_{jk} \tilde{a}_{\beta\delta} X^{k\beta} + \sum_{i,\gamma} b_i \tilde{b}_\delta a_{ij} \tilde{a}'_{\delta \gamma} Y^{i\gamma} = - \sum_{i,\beta} b_i\tilde{b}_\beta a_{ij}\tilde{a}_{\beta\delta} \frac{\partial H}{\partial \phi}(\Phi_{i\beta}, P^0_{i\beta}, P^1_{i\beta}). $$
This is the internal stage approximation to the remaining De Donder--Weyl equation $\partial_tp^0 + \partial_xp^1 = - \partial H/\partial \phi$. Note that the above form of the stationarity condition does not involve $a^{-1}$ or $\tilde{a}^{-1}$, so it is plausible that one can derive these equations without assuming the invertibility of the Runge--Kutta matrices; later, we will show that this is the case using independent internal stages.

Now, we compute the discrete forward Hamilton's equations. We have
\begin{align*}
\varphi_{1[\alpha]} &= \frac{1}{\tilde{b}_\alpha \Delta x} \frac{\partial H_d^+}{\partial \pi^0_{1[\alpha]}}\\&= \frac{1}{\tilde{b}_\alpha \Delta x} \frac{\partial K}{\partial \pi^0_{1[\alpha]}} \\
&= \varphi_{0[\alpha]} + \Delta t \sum_j b_j V^{j\alpha} - \Delta t \sum_j b_j V^{j\alpha} + \Delta t \sum_j b_j \frac{\partial H}{\partial p^0} (\Phi_{j\alpha}, P^0_{j\alpha}, P^1_{j\alpha}) \\
&= \varphi_{0[\alpha]} + \Delta t \sum_j b_j \frac{\partial H}{\partial p^0} (\Phi_{j\alpha}, P^0_{j\alpha}, P^1_{j\alpha}).
\end{align*}
%This is of course consistent with the equations $\varphi_{1[\alpha]} = \varphi_{0[\alpha]} + \Delta t \sum_j b_j V^{j\alpha}$ and $V^{j\alpha} = \frac{\partial H}{\partial p^0}(\Phi_{j\alpha}, P^0_{j\alpha},P^1_{j\alpha})$.
Similarly, 
$$ \varphi_{[i]1} = \varphi_{[i]0} + \Delta x \sum_{\beta} \tilde{b}_\beta \frac{\partial H}{\partial p^1}(\Phi_{i\beta},P^0_{i\beta},P^1_{i\beta}). $$
Computing the discrete forward Hamilton's equations for the momenta gives
\begin{align*}
\pi^0_{0[\alpha]} &= \pi^0_{1[\alpha]} + \frac{\Delta t}{\tilde{b}_\alpha} \sum_{i,\beta} b_i \tilde{b}_\beta \frac{\partial H}{\partial \phi}(\Phi_{i\beta},P^0_{i\beta},P^1_{i\beta}) \frac{\partial \Phi_{i\beta}}{\partial \varphi_{0[\alpha]}},\\
\pi^1_{[i]0} &= \pi^1_{[i]1} + \frac{\Delta x}{b_i} \sum_{j,\alpha} b_j \tilde{b}_\alpha \frac{\partial H}{\partial \phi}(\Phi_{j\alpha}, P^0_{j\alpha},P^1_{j\alpha}) \frac{\partial \Phi_{j\alpha}}{\partial \varphi_{[i]0}}.
\end{align*}
We will postpone the discussion of the discrete forward Hamilton's equations until after discussing independent internal stages, which will give a more explicit characterization of these equations.

To summarize, our method is given by
\begin{subequations}
\begin{align} \label{MPRKfirst}
\Phi_{i\alpha} &= \varphi_{0[\alpha]} + \Delta t \sum_{j}a_{ij}V^{j\alpha},\\
P^0_{i\alpha}&= \pi^0_{1[\alpha]} - \Delta t \sum_j a'_{ij}X^{j\alpha}, \label{P0} \\
\varphi_{1[\alpha]} &= \varphi_{0[\alpha]} + \Delta t \sum_j b_j V^{j\alpha},\\
\pi^0_{0[\alpha]} &= \pi^0_{1[\alpha]} - \Delta t \sum_j b'_j X^{j\alpha}, \label{pi0}\\
\hspace{0.1cm}& \nonumber \\
\Phi_{i\alpha} &= \varphi_{[i]0} + \Delta x \sum_\beta \tilde{a}_{\alpha\beta} W^{i\beta} \\
P^1_{i\alpha} &= \pi^1_{[i]1} - \Delta x \sum_\beta \tilde{a}'_{\alpha\beta}Y^{i\beta}, \label{P1}\\
\varphi_{[i]1} &= \varphi_{[i]0} + \Delta x \sum_\beta \tilde{b}_\beta W^{i\beta}, \\
\pi^1_{[i]0} &= \pi^1_{[i]1} - \Delta x \sum_\beta \tilde{b}'_\beta Y^{i\beta}. \label{pi1} \\ 
\hspace{0.1cm}& \nonumber \\
\sum_i b_i a'_{ij} \Big(&V^{i\alpha} -  \frac{\partial H}{\partial p^0}(\Phi_{i\alpha},P^0_{i\alpha}, P^1_{i\alpha}) \Big) = 0,\\
\sum_\beta \tilde{b}_\beta \tilde{a}'_{\beta\alpha} \Big(&W^{j\beta} - \frac{\partial H}{\partial p^1}(\Phi_{j\beta}, P^0_{j\beta}, P^1_{j\beta}) \Big) = 0,\\
 \sum_{k,\beta} b_j \tilde{b}_\beta a'_{jk} \tilde{a}_{\beta\delta} X^{k\beta} + \sum_{i,\gamma} & b_i \tilde{b}_\delta a_{ij} \tilde{a}'_{\delta \gamma} Y^{i\gamma} = - \sum_{i,\beta} b_i\tilde{b}_\beta a_{ij}\tilde{a}_{\beta\delta} \frac{\partial H}{\partial \phi}(\Phi_{i\beta}, P^0_{i\beta}, P^1_{i\beta}). \label{MPRKfinal}
\end{align}	
\end{subequations}

\textbf{Independent Internal Stages.} We now reformulate the above construction using independent internal stages and derive explicit conditions on the coefficients for the momenta expansion for the method to be well-defined. Recall that in the above construction, we enforced the condition that the internal stages $\Phi_{i\alpha}$ produced by both $V^{i\alpha}$ and $W^{i\alpha}$ had to be the same; we now relax this assumption and let the internal stages be independent, but subsequently enforce that they are the same by using Lagrange multipliers. Compared to the previous formulation, the use of independent internal stages has the advantage that the discrete forward Hamilton's equations can be written explicitly. Furthermore, the generalization to higher spacetime dimensions is straight-forward as opposed to the previous formulation, which would involve inverting the condition that the internal stages obtained from the various spacetime derivative approximations, $\partial_\mu \phi_d$, are consistent.

%equal internal stage condition (relating the various spacetime derivative approximations, $\partial_\mu \phi_d$) corresponding to each  direction in spacetime.

Hence, we define independent internal stages corresponding to integration in each spacetime direction,
\begin{align*}
\Phi_{i\alpha} &\equiv \phi(c_i\Delta t, \tilde{c}_\alpha \Delta x) = \phi(0, \tilde{c}_\alpha \Delta x) + \Delta t \sum_{j,\beta} V^{j\beta} \int_0^{c_i} \psi_j(s)ds\ \tilde{\psi}_\beta (\tilde{c}_\alpha) = \varphi_{0[\alpha]} + \Delta t \sum_j a_{ij} V^{j\alpha}, \\
\tilde{\Phi}_{i\alpha} &\equiv \phi(c_i\Delta t, \tilde{c}_\alpha \Delta x) = \phi(c_i\Delta t, 0) + \Delta x \sum_{j,\beta} W^{j\beta} \psi_j(c_i) \int_0^{c_\alpha} \tilde{\psi}_\beta(s)ds = \varphi_{[i]0} + \Delta x \sum_\beta \tilde{a}_{\alpha\beta} W^{i\beta}.
\end{align*}
The expansion of the other quantities are the same as the previous discussion.

We will evaluate the Hamiltonian at the weighted combination $\Phi^\theta_{i\alpha} \equiv \theta \Phi_{i\alpha} + (1-\theta)\tilde{\Phi}_{i\alpha}$ for some arbitrary parameter $\theta \in \mathbb{R}$ and subsequently enforce that the two sets of internal stages are the same through a Lagrange multiplier term $\sum_{i,\alpha}\lambda_{i\alpha}(\Phi_{i\alpha} - \tilde{\Phi}_{i\alpha})$. Thus, after enforcing the stationarity conditions, $\Phi^\theta_{i\alpha} = \Phi_{i\alpha} = \tilde{\Phi}_{i\alpha}$. In this formulation, $K$ is
\begin{align*}
\hspace{-0.5cm} K(\{\varphi_A,\pi_B, V^{i\alpha}, W^{i\alpha}, X^{i\alpha},Y^{i\alpha},\lambda_{i\alpha}\}) &= \Delta x \sum_\alpha \tilde{b}_\alpha \pi^0_{1[\alpha]}(\varphi_{0[\alpha]} + \Delta t \sum_{j} b_j V^{j\alpha})\\
&\qquad + \Delta t \sum_i b_i \pi^1_{[i]1}(\varphi_{[i]0} + \Delta x \sum_\beta \tilde{b}_\beta W^{i\beta}) \\
&\qquad- \Delta t \Delta x \sum_{i,\alpha} b_i \tilde{b}_\alpha \Big(\pi^0_{1[\alpha]} - \Delta t \sum_j a'_{ij}X^{j\alpha} \Big) V^{i\alpha} \\
&\qquad-  \Delta t \Delta x \sum_{i,\alpha} b_i\tilde{b}_\alpha \Big( \pi^1_{[i]1} - \Delta x \sum_{\beta}\tilde{a}'_{\alpha\beta} Y^{i\beta} \Big) W^{i\alpha} \\
&\qquad+ \Delta t \Delta x \sum_{i,\alpha}b_i\tilde{b}_\alpha H(\Phi^\theta_{i\alpha},P^0_{i\alpha},P^1_{i\alpha}) + \sum_{i,\alpha}\lambda_{i\alpha}(\Phi_{i\alpha} - \tilde{\Phi}_{i\alpha});
\end{align*}
where now both $\{V^{i\alpha}\}$ and $\{W^{i\alpha}\}$ are independent. The discrete boundary Hamiltonian $H_d^+$ is given by extremizing $K$ with respect to all of the internal variables, $\{V^{i\alpha}, W^{i\alpha}, X^{i\alpha}, Y^{i\alpha}, \lambda_{i\alpha}\}$.

Extremizing $K$ with respect to $X^{j\alpha}$ and $Y^{j\alpha}$ gives the same stationarity conditions as the previous case of equal internal stages, since the momenta expansions were unchanged, except with $H$ evaluated at $\Phi^\theta_{i\alpha}$. Namely,
\begin{subequations}
\begin{align} \label{stationarityX}
\sum_i b_i a'_{ij} \Big( V^{i\alpha} -  \frac{\partial H}{\partial p^0}(\Phi^\theta_{i\alpha},P^0_{i\alpha}, P^1_{i\alpha}) \Big) &= 0, \\
\sum_\beta \tilde{b}_\beta \tilde{a}'_{\beta\alpha} \Big(W^{j\beta} - \frac{\partial H}{\partial p^1}(\Phi^\theta_{j\beta}, P^0_{j\beta}, P^1_{j\beta}) \Big) &= 0. \label{stationarityY}
\end{align}
\end{subequations}
Extremizing $K$ with respect to $V^{j\beta}$,
\begin{align*}
0 = \frac{\partial K}{\partial V^{j\beta}} &= \Delta t \Delta x b_j \tilde{b}_\beta \pi^0_{1[\beta]} - \Delta t \Delta x b_j \tilde{b}_\beta P^0_{j\beta} + \Delta t^2 \Delta x \sum_i b_i\tilde{b}_\beta a_{ij} \theta \frac{\partial H}{\partial \phi}(\Phi^\theta_{i\beta},P^0_{i\beta},P^1_{i\beta}) + \Delta t \sum_i \lambda_{i\beta} a_{ij} \\
&= \Delta t^2 \Delta x b_j \tilde{b}_\beta \sum_k a'_{jk} X^{k\beta}  + \Delta t^2 \Delta x \sum_i b_i\tilde{b}_\beta a_{ij} \theta \frac{\partial H}{\partial \phi}(\Phi^\theta_{i\beta},P^0_{i\beta},P^1_{i\beta}) + \Delta t \sum_i \lambda_{i\beta} a_{ij}.
\end{align*}
Dividing by $\Delta t^2 \Delta x$,
\begin{equation} \label{stationarityV}
\sum_k b_j \tilde{b}_\beta a'_{jk} X^{k\beta}  + \sum_i b_i\tilde{b}_\beta a_{ij} \theta \frac{\partial H}{\partial \phi}(\Phi^\theta_{i\beta},P^0_{i\beta},P^1_{i\beta}) + \frac{1}{\Delta t \Delta x} \sum_i \lambda_{i\beta} a_{ij} = 0.
\end{equation}
Similarly, extremizing $K$ with respect to $W^{j\beta}$ (and dividing by $\Delta t \Delta x^2$) gives
\begin{equation} \label{stationarityW}
\sum_\alpha b_j\tilde{b}_\beta \tilde{a}'_{\beta\alpha}Y^{j\alpha} + \sum_\alpha b_j \tilde{b}_\alpha \tilde{a}_{\alpha\beta} (1-\theta) \frac{\partial H}{\partial \phi}(\Phi^\theta_{j\alpha},P^0_{j\alpha},P^1_{j\alpha}) - \frac{1}{\Delta t\Delta x}\sum_\alpha \lambda_{j\alpha}\tilde{a}_{\alpha\beta} = 0.
\end{equation}
Let us combine these two stationarity conditions to eliminate $\theta$ and the Lagrange multiplier terms. Multiply equation (\ref{stationarityV}) by $\tilde{a}_{\beta\delta}$ and sum over $\beta$; multiply equation (\ref{stationarityW}) by $a_{ji}$ and sum over $j$. Subsequently, add the two resulting equations. This gives
\begin{equation}\label{stationarityVW}
\sum_{k,\beta} b_j \tilde{b}_\beta a'_{jk} \tilde{a}_{\beta\delta} X^{k\beta} + \sum_{i,\gamma} b_i \tilde{b}_\delta a_{ij} \tilde{a}'_{\delta \gamma} Y^{i\gamma} = - \sum_{i,\beta} b_i\tilde{b}_\beta a_{ij}\tilde{a}_{\beta\delta} \frac{\partial H}{\partial \phi}(\Phi^\theta_{i\beta}, P^0_{i\beta}, P^1_{i\beta}). 
\end{equation}
Finally, extremizing $K$ with respect to $\lambda_{i\alpha}$ enforces that the independent internal stages are the same, $0 = \partial K/\partial \lambda_{i\alpha} = \Phi_{i\alpha} - \tilde{\Phi}_{i\alpha}$, and hence $\Phi^\theta_{i\alpha} = \Phi_{i\alpha} = \tilde{\Phi}_{i\alpha}$. We have rederived the stationarity conditions that we saw in the case of equal internal stages, without the assumption of invertibility of the Runge--Kutta matrices, $(a_{ij})$, $(\tilde{a}_{\alpha\beta})$. 

Now, we aim to provide a more explicit characterization of the discrete forward Hamilton's equations. We will assume again that the Runge--Kutta matrices $(a_{ij})$, $(\tilde{a}_{\alpha\beta})$ are invertible. Computing the discrete forward Hamilton's equations for the field boundary values,
\begin{align*}
\varphi_{1[\alpha]} &= \frac{1}{\tilde{b}_\alpha \Delta x} \frac{\partial H_d^+}{\partial \pi^0_{1[\alpha]}} = \frac{1}{\tilde{b}_\alpha \Delta x} \frac{\partial K}{\partial \pi^0_{1[\alpha]}} = \varphi_{0[\alpha]} + \Delta t \sum_j b_j \frac{\partial H}{\partial p^0} (\Phi^\theta_{j\alpha}, P^0_{j\alpha}, P^0_{j\alpha}), \\
\varphi_{[i]1} &= \frac{1}{b_i \Delta t} \frac{\partial H_d^+}{\partial \pi^1_{[i]1}} = \frac{1}{b_i \Delta t} \frac{\partial K}{\partial \pi^1_{[i]1}} = \varphi_{[i]0} + \Delta x \sum_{\beta} \tilde{b}_\beta \frac{\partial H}{\partial p^1}(\Phi^\theta_{i\beta},P^0_{i\beta},P^1_{i\beta}).
\end{align*}
Recall that we also have the expansion for the field boundary values 
\begin{align*}
\varphi_{1[\alpha]} &= \varphi_{0[\alpha]} + \Delta t \sum_j b_j V^{j\alpha},\\
\varphi_{[i]1} &= \varphi_{[i]0} + \Delta x \sum_\beta \tilde{b}_\beta W^{i\beta}.
\end{align*}
We will see shortly that, with a particular condition on the coefficients of the momenta expansion, the discrete forward Hamilton's equations for the field values are consistent with the field expansions, i.e., that $V^{j\alpha} = \frac{\partial H}{\partial p^0} (\Phi^\theta_{j\alpha}, P^0_{j\alpha}, P^0_{j\alpha})$ and similarly $W^{i\beta} = \frac{\partial H}{\partial p^1}(\Phi^\theta_{i\beta},P^0_{i\beta},P^1_{i\beta})$.

First, we compute the discrete forward Hamilton's equations for the momenta boundary values,
\begin{align*}
\pi^0_{0[\alpha]} &= \frac{1}{\tilde{b}_\alpha \Delta x} \frac{\partial H_d^+}{\partial \varphi_{0[\alpha]}} = \frac{1}{\tilde{b}_\alpha \Delta x} \frac{\partial K}{\partial \varphi_{0[\alpha]}} = \pi^0_{1[\alpha]} + \Delta t \sum_i b_i \theta \frac{\partial H}{\partial \phi}(\Phi^\theta_{i\alpha},P^0_{i\alpha},P^1_{i\alpha}) + \frac{1}{\tilde{b}_\alpha \Delta x} \sum_i \lambda_{i\alpha}, \\
\pi^1_{[i]0} &= \frac{1}{b_i \Delta t} \frac{\partial H_d^+}{\partial \varphi_{[i]0}} = \frac{1}{b_i \Delta t} \frac{\partial K}{\partial \varphi_{[i]0}} = \pi^1_{[i]1} + \Delta x \sum_\alpha \tilde{b}_\alpha (1-\theta) \frac{\partial H}{\partial \theta}(\Phi^\theta_{i\alpha},P^0_{i\alpha},P^1_{i\alpha}) - \frac{1}{b_i\Delta t} \sum_\alpha \lambda_{i\alpha}.
\end{align*}
For our method to be well-defined, these are required to be consistent with the momenta expansions,
\begin{align*}
\pi^0_{0[\alpha]} &= \pi^0_{1[\alpha]} - \Delta t \sum_j b_j' X^{j\alpha},\\
\pi^1_{[i]0} &= \pi^1_{[i]1} - \Delta x \sum_\beta \tilde{b}'_\beta Y^{i\beta}.
\end{align*}
To do this, we solve the stationarity conditions (\ref{stationarityV}) and (\ref{stationarityW}) for the Lagrange multipliers. Multiply equation (\ref{stationarityV}) by $(a^{-1})_{jl}$ and sum over $j$; multiply equation (\ref{stationarityW}) by $(\tilde{a}^{-1})_{\beta\gamma}$ and sum over $\beta$. This gives
\begin{align*}
\lambda_{l\beta} &= -\Delta t \Delta x b_l \tilde{b}_\beta \theta \frac{\partial H}{\partial \phi}(\Phi^\theta_{l\beta},P^0_{l\beta},P^1_{l\beta}) - \Delta t \Delta x \sum_{j,k} b_j \tilde{b}_\beta a'_{jk} (a^{-1})_{jl} X^{k\beta}, \\
\lambda_{j\gamma} &= \Delta t \Delta x b_j \tilde{b}_\gamma (1-\theta) \frac{\partial H}{\partial \phi}(\Phi^\theta_{j\gamma},P^0_{j\gamma},P^1_{j\gamma}) + \Delta t \Delta x \sum_{\alpha,\beta} b_j \tilde{b}_\beta \tilde{a}'_{\beta\alpha} (\tilde{a}^{-1})_{\beta\gamma}Y^{j\alpha}.
\end{align*}
Plugging these into the respective discrete forward Hamilton's equations for the momenta boundary values, we have
\begin{align*}
\pi^0_{0[\alpha]} &= \pi^0_{1[\alpha]} - \Delta t \sum_{j,k,l} b_j a'_{jk} (a^{-1})_{jl} X^{k\beta} \overset{!}{=} \pi^0_{1[\alpha]} - \Delta t \sum_k b_k' X^{k\alpha},\\
\pi^1_{[i]0} &= \pi^1_{[i]1} - \Delta x \sum_{\alpha,\beta,\gamma} \tilde{b}_\beta \tilde{a}'_{\beta\alpha} (\tilde{a}^{-1})_{\beta\gamma}Y^{i\alpha} \overset{!}{=} \pi^1_{[i]1} - \Delta x \sum_\alpha \tilde{b}'_\alpha Y^{i\alpha}.
\end{align*}

\begin{prop}\label{Symplecticity-Consistency}
The method arising from approximating the internal stages with the partitioned Runge--Kutta expansion is well-defined if and only if the partitioned Runge--Kutta method is symplectic in both space and time, i.e.
\begin{align*}
\sum_{j,l}b_ja_{jk}' (a^{-1})_{jl} &= b_k',\\
\sum_{\beta,\gamma} \tilde{b}_\beta \tilde{a}'_{\beta\alpha} (\tilde{a}^{-1})_{\beta\gamma} &= \tilde{b}'_\alpha.
\end{align*}
A sufficient condition is the usual choice of symplectic partitioned Runge--Kutta coefficients,
\begin{align*}
a'_{jk} &= \frac{b'_k a_{kj}}{b_j}, \\
\tilde{a}'_{\beta\alpha} &= \frac{\tilde{b}'_\alpha \tilde{a}_{\alpha\beta}}{\tilde{b}_\beta}.
\end{align*}
(We will see after expressing the momenta internal stages in terms of $\pi_A$ instead of $\pi_B$ that these are the usual choice of symplectic partitioned Runge--Kutta coefficients).
\begin{proof}
By comparing the momenta expansions to the discrete forward Hamilton's equations for the momenta, we must have 
\begin{subequations}
\begin{align}\label{SymplecticCoefficients1}
\sum_{j,k,l} b_j a'_{jk} (a^{-1})_{jl} X^{k\beta} &= \sum_k b_k' X^{k\alpha},\\
\sum_{\alpha,\beta,\gamma} \tilde{b}_\beta \tilde{a}'_{\beta\alpha} (\tilde{a}^{-1})_{\beta\gamma}Y^{i\alpha} &= \sum_\alpha \tilde{b}'_\alpha Y^{i\alpha}. \label{SymplecticCoefficients2}
\end{align}
\end{subequations}
Since the internal variables $\{X^{i\alpha},Y^{i\alpha}\}$ are generally arbitrary (depending on the choice of Hamiltonian and the supplied boundary data), the above must hold for arbitrary choices of $\{X^{i\alpha}\}$ and $\{Y^{i\alpha}\}$; hence, we have the necessary and sufficient conditions
\begin{align*}
\sum_{j,l}b_ja_{jk}' (a^{-1})_{jl} &= b_k',\\
\sum_{\beta,\gamma} \tilde{b}_\beta \tilde{a}'_{\beta\alpha} (\tilde{a}^{-1})_{\beta\gamma} &= \tilde{b}'_\alpha. 
\end{align*}
Plugging in the choice (\ref{SymplecticCoefficients1}) and (\ref{SymplecticCoefficients2}) to the left hand sides of the above conditions,
\begin{align*}
\sum_{j,l}b_ja_{jk}' (a^{-1})_{jl} &= \sum_{j,l}b_k'a_{kj} (a^{-1})_{jl} = \sum_{l}b_k'\delta_{kl} = b'_k,\\
\sum_{\beta,\gamma} \tilde{b}_\beta \tilde{a}'_{\beta\alpha} (\tilde{a}^{-1})_{\beta\gamma} &= \sum_{\beta,\gamma} \tilde{b}'_\alpha \tilde{a}_{\alpha\beta} (\tilde{a}^{-1})_{\beta\gamma} = \sum_{\gamma}\tilde{b}'_\alpha \delta_{\alpha\gamma} = \tilde{b}'_\alpha;
\end{align*}
so this choice is sufficient for the method to be well-defined.
\end{proof}
\end{prop}

Now, consider the stationarity conditions (\ref{stationarityX}) and (\ref{stationarityY}). Plugging in the choice of coefficients (\ref{SymplecticCoefficients1}) and (\ref{SymplecticCoefficients2}), we have
\begin{align*}
\sum_i b_j' a_{ji} \Big( V^{i\alpha} -  \frac{\partial H}{\partial p^0}(\Phi^\theta_{i\alpha},P^0_{i\alpha}, P^1_{i\alpha}) \Big) &= 0, \\
\sum_\beta \tilde{b}'_\alpha \tilde{a}_{\alpha\beta} \Big(W^{j\beta} - \frac{\partial H}{\partial p^1}(\Phi^\theta_{j\beta}, P^0_{j\beta}, P^1_{j\beta}) \Big) &= 0.
\end{align*}
Since $(a_{ji})$ and $(\tilde{a}_{\alpha\beta})$ are invertible, we have $V^{j\alpha} = \frac{\partial H}{\partial p^0} (\Phi^\theta_{j\alpha}, P^0_{j\alpha}, P^0_{j\alpha})$ and $W^{i\beta} = \frac{\partial H}{\partial p^1}(\Phi^\theta_{i\beta},P^0_{i\beta},P^1_{i\beta})$ so that the discrete forward Hamilton's equations for the field boundary values are also consistent with the their expansions. Similarly, plugging this choice of coefficients into the stationarity condition (\ref{stationarityVW}) gives
$$ \sum_{k,\beta} b'_k \tilde{b}_\beta a_{kj} \tilde{a}_{\beta\delta} X^{k\beta} + \sum_{i,\gamma} b_i \tilde{b}'_\gamma a_{ij} \tilde{a}_{\gamma\delta} Y^{i\gamma} = - \sum_{i,\beta} b_i\tilde{b}_\beta a_{ij}\tilde{a}_{\beta\delta} \frac{\partial H}{\partial \phi}(\Phi^\theta_{i\beta}, P^0_{i\beta}, P^1_{i\beta}). $$
To invert this relation, we impose $b_k' = b_k, \tilde{b}'_\gamma = \tilde{b}_\gamma$. Note that the matrix with $jk$ entry $b_ka_{kj}$ is invertible since $(a_{jk})$ is (its transpose is obtained by multiplying the $i^{th}$ row of $(a_{ij})$ by $b_i \neq 0$, so the rows are still linearly independent) and similarly for the matrix with $\delta \gamma$ entry $\tilde{b}_\gamma \tilde{a}_{\gamma\delta}$. Hence, this stationarity condition can be inverted to give
$$ X^{i\alpha} + Y^{i\alpha} = - \frac{\partial H}{\partial \phi}(\Phi^\theta_{i\beta}, P^0_{i\beta},P^1_{i\beta}). $$

Finally, to write our method in the traditional form of a partitioned Runge--Kutta method, we express the internal stages $P^0_{i\alpha}$ and $P^1_{i\alpha}$ in terms of $\pi_A$ instead of $\pi_B$, by plugging equations (\ref{pi0}) and (\ref{pi1}) into equations (\ref{P0}) and (\ref{P1}) respectively, 
\begin{align*}
P^0_{i\alpha} &= \pi^0_{0[\alpha]} + \Delta t \sum_j (b_j - a_{ij}')X^{j\alpha} = \pi^0_{0[\alpha]} + \Delta t \sum_j \underbrace{\frac{b_jb_i - b_ja_{ji}}{b_i}}_{\equiv a^{(2)}_{ij}}X^{j\alpha}, \\
P^1_{i\alpha} &= \pi^1_{[i]0} + \Delta x \sum_\beta (\tilde{b}_\beta - \tilde{a}'_{\alpha\beta} )Y^{i\beta} = \pi^1_{[i]0} + \Delta x \sum_\beta \underbrace{\frac{\tilde{b}_\beta \tilde{b}_\alpha - \tilde{b}_\beta \tilde{a}_{\beta\alpha}}{\tilde{b}_\alpha}}_{\equiv \tilde{a}^{(2)}_{\alpha\beta}} Y^{i\beta}.
\end{align*}
To summarize, our method is
\begin{subequations}
\begin{align} \label{MPRK2first}
\Phi_{i\alpha} &= \varphi_{0[\alpha]} + \Delta t \sum_{j}a_{ij}V^{j\alpha},\\
P^0_{i\alpha}&= \pi^0_{0[\alpha]} + \Delta t \sum_j a^{(2)}_{ij}X^{j\alpha}, \\
\varphi_{1[\alpha]} &= \varphi_{0[\alpha]} + \Delta t \sum_j b_j V^{j\alpha},\\
\pi^0_{1[\alpha]} &= \pi^0_{0[\alpha]} + \Delta t \sum_j b_j X^{j\alpha}, \label{pi0b} \\
\hspace{0.1cm}& \nonumber \\
\Phi_{i\alpha} &= \tilde{\Phi}_{i\alpha} = \varphi_{[i]0} + \Delta x \sum_\beta \tilde{a}_{\alpha\beta} W^{i\beta} \\
P^1_{i\alpha} &= \pi^1_{[i]0} + \Delta x \sum_\beta \tilde{a}^{(2)}_{\alpha\beta}Y^{i\beta}, \\
\varphi_{[i]1} &= \varphi_{[i]0} + \Delta x \sum_\beta \tilde{b}_\beta W^{i\beta}, \\
\pi^1_{[i]1} &= \pi^1_{[i]0} + \Delta x \sum_\beta \tilde{b}_\beta Y^{i\beta}. \label{pi1b} \\ 
\hspace{0.1cm}& \nonumber \\
V^{i\alpha} &=  \frac{\partial H}{\partial p^0}(\Phi_{i\alpha},P^0_{i\alpha}, P^1_{i\alpha}) ,\\
W^{i\alpha} &= \frac{\partial H}{\partial p^1}(\Phi_{i\alpha}, P^0_{i\alpha}, P^1_{i\alpha}),\\
X^{i\alpha} + Y^{i\alpha} &= -\frac{\partial H}{\partial \phi}(\Phi_{i\alpha},P^0_{i\alpha},P^1_{i\alpha}), \label{MPRK2final}
\end{align}
\end{subequations}
where $a^{(2)}_{ij} = \frac{b_jb_i - b_ja_{ji}}{b_i}$ and $\tilde{a}^{(2)}_{\alpha\beta} = \frac{\tilde{b}_\beta \tilde{b}_\alpha - \tilde{b}_\beta \tilde{a}_{\beta\alpha}}{\tilde{b}_\alpha}$. This is the usual form of a multisymplectic partitioned Runge--Kutta method. Note that our choice of $a^{(2)}_{ij}$ and $\tilde{a}^{(2)}_{\alpha\beta}$ (or equivalently our choice of $a'_{ij}, \tilde{a}'_{\alpha\beta}$) is the usual choice for the coefficients in the momenta expansion for a partitioned Runge--Kutta method to be multisymplectic (see, for example, \citet{HoLiSu2006}, \citet{Re2000}, \citet{RyMcFr2007}). Interestingly, however, from our perspective, our method based on the discrete boundary Hamiltonian is guaranteed to be multisymplectic so we had to impose no such conditions on the coefficients to ensure multisymplecticity; rather, the conditions for the coefficients arose from the necessity of the method to be well-defined, i.e., that the expansions of the field and momenta boundary values agreed with the discrete forward Hamilton's equations.

\begin{remark}\label{RKmatrixinvertibility}
In the above construction, we saw that the Runge--Kutta matrices $(a_{ij})$ and $(\tilde{a}_{\alpha\beta})$ were required to be invertible. We can see this directly from the internal stage expansions
\begin{align*}
\Phi_{i\alpha} &= \varphi_{0[\alpha]} + \Delta t \sum_j a_{ij}V^{j\alpha},\\
\tilde{\Phi}_{i\alpha} &= \varphi_{[i]0} + \Delta x \sum_\beta \tilde{a}_{\alpha\beta}W^{i\beta},
\end{align*}
since only when $(a_{ij})$ and $(\tilde{a}_{\alpha\beta})$ are invertible is extremizing $K$ over $V^{i\alpha}$ and $W^{i\alpha}$ equivalent to extremizing $K$ over $\Phi_{i\alpha}$ and $\tilde{\Phi}_{i\alpha}$, respectively. In the case of non-invertible Runge--Kutta matrices, the internal stages $\Phi_{i\alpha}$ and $\tilde{\Phi}_{i\alpha}$ do not depend independently on all of the $V^{i\alpha}$, $W^{i\alpha}$. For collocation Runge--Kutta methods, non-invertibility arises from the choice of the first quadrature point $c_1 = 0$. In our construction, if we choose $c_1 = 0$, then we are specifying an internal stage at a quadrature point where the field boundary value $\varphi_A$ is already specified; thus, the internal stage at this point is not free to extremize over. Hence, in the non-invertible case, one has to use the specified boundary values to eliminate the degeneracy in the internal variables $V^{i\alpha}$ and $W^{i\alpha}$, reducing the number of internal variables to an independent subcollection of internal variables. Subsequently, one extremizes only over this independent subcollection of internal variables.
\end{remark}

\begin{remark}
It should also be remarked that while certain types of Galerkin multisymplectic Hamiltonian variational integrators recover multisymplectic partitioned Runge--Kutta methods, it remains to see whether there is a more general correspondence between Galerkin multisymplectic Hamiltonian variational integrators with a class of modified multisymplectic partitioned Runge--Kutta methods for the case of spacetime tensor product (hyper)rectangular meshes. This would generalize the connection between Galerkin variational integrators and modified symplectic Runge--Kutta methods in the ODE setting that was observed in~\cite{Ob2016}.
\end{remark}
%\begin{remark}
\textbf{Momenta Internal Stages.}
In the above construction, we saw that we had to enforce consistency conditions on the momenta expansion coefficients in order for the method (\ref{MPRK2first})-(\ref{MPRK2final}) to be well-defined. The issue is that we over-constrained the form of the momenta internal stages via our particular choice of expansion, since ultimately our goal was to derive the class of multisymplectic partitioned Runge--Kutta methods within our variational framework. One can avoid this problem altogether by working directly with the momenta internal stages $P^0_{i\alpha}$ and $P^1_{i\alpha}$ instead of the internal variables $X^{i\alpha}$ and $Y^{i\alpha}$, although the method will not ultimately be in the form of a multisymplectic partitioned Runge--Kutta method. This is possible for the momenta internal stages since the action does not depend on the derivatives of the momenta, unlike the field variable. We outline this procedure. 

Assume the same expansions of $\Phi_{i\alpha}, \tilde{\Phi}_{i\alpha}, \varphi_{1[\alpha]}, \varphi_{[i]1}$ in terms of $\{V^{i\alpha}\}$ and $\{W^{i\alpha}\}$. For the momenta, we work directly with the internal stages $P^0_{i\alpha}$, $P^1_{i\alpha}$ instead of using an expansion. In this case, $K$ is
\begin{align*}
K(\{\varphi_A,\pi_B,V^{i\alpha}, W^{i\alpha}, P^0_{i\alpha},P^1_{i\alpha}, \lambda_{i\alpha} \}) &= \Delta x \sum_\alpha \tilde{b}_\alpha \pi^0_{1[\alpha]}(\varphi_{0[\alpha]} + \Delta t \sum_{j} b_j V^{j\alpha})\\
&\qquad + \Delta t \sum_i b_i \pi^1_{[i]1}(\varphi_{[i]0} + \Delta x \sum_\beta \tilde{b}_\beta W^{i\beta}) \\
&\qquad- \Delta t \Delta x \sum_{i,\alpha} b_i \tilde{b}_\alpha P^0_{i\alpha} V^{i\alpha} -  \Delta t \Delta x \sum_{i,\alpha} b_i\tilde{b}_\alpha P^1_{i\alpha} W^{i\alpha} \\
&\qquad+ \Delta t \Delta x \sum_{i,\alpha}b_i\tilde{b}_\alpha H(\Phi^\theta_{i\alpha},P^0_{i\alpha},P^1_{i\alpha}) + \sum_{i,\alpha}\lambda_{i\alpha}(\Phi_{i\alpha} - \tilde{\Phi}_{i\alpha}).
\end{align*}
$H_d^+$ is obtained by extremizing $K$ over the internal variables, $\{V^{i\alpha},W^{i\alpha},P^0_{i\alpha},P^1_{i\alpha}, \lambda_{i\alpha}\}$. 

The stationarity condition $\partial K/\partial P^0_{i\alpha} = 0$ (divided by $\Delta t \Delta x b_i\tilde{b}_\alpha)$ gives
$$ V^{i\alpha} = \frac{\partial H}{\partial p^0}(\Phi^\theta_{i\alpha},P^0_{i\alpha},P^1_{i\alpha}). $$
Similarly, the stationarity condition $\partial K/\partial P^1_{i\alpha} = 0$ (divided by $\Delta t \Delta x b_i \tilde{b}_\alpha$) gives
$$ W^{i\alpha} = \frac{\partial H}{\partial p^1}(\Phi^\theta_{i\alpha},P^0_{i\alpha},P^1_{i\alpha}). $$
The stationarity condition $\partial K/\partial \lambda_{i\alpha} = 0$ gives $\Phi_{i\alpha} = \tilde{\Phi}_{i\alpha}$. The stationarity conditions $\partial K/\partial V^{j\beta} = 0$ and $\partial K/\partial W^{j\beta} = 0$ give respectively
\begin{align*}
\Delta t \Delta x b_j \tilde{b}_\beta (\pi^0_{1[\beta]} - P^0_{j\beta}) + \Delta t^2 \Delta x \sum_i b_i\tilde{b}_\beta a_{ij}\theta \frac{\partial H}{\partial \phi}(\Phi^\theta_{i\beta}, P^0_{i\beta}, P^1_{i\beta}) + \Delta t \sum_i \lambda_{i\beta}a_{ij} &= 0, \\
\Delta t \Delta x b_j \tilde{b}_\beta (\pi^1_{[j]1} - P^1_{j\beta}) + \Delta x^2 \Delta t \sum_\alpha b_j\tilde{b}_\alpha \tilde{a}_{\alpha\beta} (1-\theta) \frac{\partial H}{\partial \phi}(\tilde{\Phi}_{j\alpha},P^0_{j\alpha},P^1_{j\alpha}) - \Delta x \sum_\alpha \lambda_{j\alpha} \tilde{a}_{\alpha\beta} &= 0.
\end{align*}
Performing the same procedure we used to combine equations (\ref{stationarityV}) and (\ref{stationarityW}) to eliminate $\theta$ and the Lagrange multipliers, these two stationarity conditions can be combined to give
$$ \sum_\beta b_j \tilde{b}_\beta \tilde{a}_{\beta\delta} \frac{(\pi^0_{1[\beta]} - P^0_{j\beta})}{\Delta x} + \sum_i b_i \tilde{b}_\delta a_{ij} \frac{(\pi^1_{[i]1} - P^1_{i\delta})}{\Delta t} = - \sum_{i,\beta}b_i\tilde{b}_\beta a_{ij}\tilde{a}_{\beta\delta} \frac{\partial H}{\partial \phi}(\Phi^\theta_{i\beta},P^0_{i\beta},P^1_{i\beta}). $$
This combined condition, together with the other stationarity conditions $V^{i\alpha} = \partial H/\partial p^0(\Phi^\theta_{i\alpha},P^0_{i\alpha},P^1_{i\alpha})$, $W^{i\alpha} = \partial H/\partial p^1(\Phi^\theta_{i\alpha},P^0_{i\alpha},P^1_{i\alpha}),$ and $\Phi_{i\alpha} = \tilde{\Phi}_{i\alpha}$ (ranging over all free indices) can be used to solve for the collection of internal variables $\{V^{i\alpha},W^{i\alpha},P^0_{i\alpha},P^1_{i\alpha}\}_{i,\alpha}$ in terms of the supplied boundary data. 

To conclude, we compute the discrete forward Hamilton's equations. For the field boundary values,
\begin{align*}
\varphi_{1[\alpha]} &= \frac{1}{\tilde{b}_\alpha \Delta x} \frac{\partial K}{\partial \pi^0_{1[\alpha]}} = \varphi_{0[\alpha]} + \Delta t \sum_j b_j V^{j\alpha},\\
\varphi_{[i]1} &= \frac{1}{b_i \Delta t} \frac{\partial K}{\partial \pi^1_{[i]1}} = \varphi_{[i]0} + \Delta x \sum_\beta b_\beta W^{i\beta}.
\end{align*}
Note that these equations already agree with the field expansion. For the momenta boundary values,
\begin{align*}
\pi^0_{0[\alpha]} &= \frac{1}{\tilde{b}_\alpha \Delta x} \frac{\partial K}{\partial \varphi_{0[\alpha]}} = \pi^0_{1[\alpha]} + \Delta t \sum_i b_i \theta \frac{\partial H}{\partial \phi}(\Phi^\theta_{i\alpha},P^0_{i\alpha},P^1_{i\alpha}) + \frac{1}{\tilde{b}_\alpha \Delta x}\sum_i \lambda_{i\alpha},\\
\pi^1_{[i]0} &= \frac{1}{b_i\Delta t} \frac{\partial K}{\partial \varphi_{[i]0}} = \pi^1_{[i]1} + \Delta x \sum_\alpha \tilde{b}_\alpha (1-\theta) \frac{\partial H}{\partial \phi}(\Phi^\theta_{i\alpha},P^0_{i\alpha},P^1_{i\alpha}) - \frac{1}{b_i \Delta t}\sum_\alpha \lambda_{i\alpha}.
\end{align*}
As we did before for the partitioned Runge--Kutta method, we can act on the stationarity conditions $\partial K/\partial V^{j\beta} = 0 = \partial K/\partial W^{j\beta}$ by the inverses of the Runge--Kutta matrices to solve for the Lagrange multipliers and substitute them into the discrete forward Hamilton's equations for the momenta, ultimately eliminating $\theta$ and the Lagrange multipliers. The discrete forward Hamilton's equations for the momenta are then
\begin{align*}
\pi^0_{0[\alpha]} &= \pi^0_{1[\alpha]} - \Delta t \sum_{j,l} b_j (a^{-1})_{jl} \frac{\pi^0_{1[\alpha]} - P^0_{j\alpha}}{\Delta x}, \\
\pi^1_{[i]0} &= \pi^1_{[i]1} - \Delta x \sum_{\alpha,\beta} \tilde{b}_\alpha (\tilde{a}^{-1})_{\alpha\beta} \frac{\pi^1_{[i]1} - P^1_{i\alpha}}{\Delta t}.
\end{align*}
Hence, by working with the internal stages for the momenta directly, as opposed to utilizing an expansion, we see that the method we derived is already well-defined (and also automatically multisymplectic), although it is not directly in the form of a multisymplectic partitioned Runge--Kutta method. 

These various approaches demonstrate the versatility of our variational framework; once one chooses an approximation for the fields, its derivatives, and the momenta (as well as some approximation for the various integrals involved), one can construct the discrete boundary Hamiltonian and subsequently the variational framework produces a multisymplectic integrator. If one over-constrains the form of the momenta expansion, as opposed to using the internal stages directly, one must also check whether the method is well-defined. Another approach that is possible within this framework is to discretize at the level of the field using some (possibly non-tensor product) function space and subsequently take derivatives of the basis functions to obtain an approximation of the derivatives of the fields. For example, we expect that utilizing spectral element bases to discretize at the level of the field within our framework will produce multisymplectic spectral discretizations like those obtained in \citet{BrRe2001b}, \citet{IsSc2004,IsSc2006}. Another interesting application of our construction would be to construct multisymplectic discretizations of the total exterior algebra bundle (see \citet{BrRe2006}) using Galerkin discretizations arising from the Finite Element Exterior Calculus framework (\citet{ArFaWi2006, ArFaWi2010}, \citet{Hi2002}), allowing one to discretize Hamiltonian PDEs with more general configuration bundles.
%\end{remark}

\subsection{Multisymplecticity Revisited}\label{Multisymplecticity Revisited}
Now, we discuss in what sense the discrete multisymplectic form formula (\ref{Discrete MFF}) corresponds to our discretization of the field equations. Consider the integral form of the De Donder--Weyl equations over $\Box = [0,\Delta t]\times [0,\Delta x]$,
\begin{subequations}
\begin{align} \label{integralDDW1}
\int_\Box \Big(\partial_\mu p^\mu + \frac{\partial H}{\partial\phi}(\phi,p^0,p^1)\Big) d^2x &= 0,\\ \label{integralDDW2}
\int_\Box \Big(\partial_0\phi - \frac{\partial H}{\partial p^0}(\phi,p^0,p^1)\Big) d^2x &= 0., \\ \label{integralDDW3}
\int_\Box \Big(\partial_1\phi - \frac{\partial H}{\partial p^1}(\phi,p^0,p^1)\Big) d^2x &= 0.
\end{align}
\end{subequations}
Applying our quadrature approximation to equation (\ref{integralDDW1}),
\begin{align*}
\hspace{-1.5cm} 0 &= \int_0^{\Delta t} \int_0^{\Delta x}\Big(\partial_0p^0 + \partial_1p^1 + \frac{\partial H}{\partial\phi}(\phi,p^0,p^1)\Big) dxdt \\
&= \int_0^{\Delta x}(p^0|_{t=\Delta t} - p^0|_{t=0}) dx + \int_0^{\Delta t}(p^1|_{x=\Delta x} - p^0|_{x=0})dt + \int_0^{\Delta t} \int_0^{\Delta x} \frac{\partial H}{\partial\phi}(\phi,p^0,p^1) dxdt \\
&\approx \Delta x \sum_{\alpha} \tilde{b}_\alpha (p^0|_{(\Delta t, \tilde{c}_\alpha \Delta x)} - p^0|_{(0, \tilde{c}_\alpha \Delta x)} ) + \Delta t \sum_i b_i (p^1|_{(c_i\Delta t, \Delta x)} - p^1|_{(c_i\Delta t,0)})\\
&\qquad+ \Delta t \Delta x \sum_{i,\alpha} \frac{\partial H}{\partial \phi} (\phi,p^0,p^1)|_{(c_i\Delta t, \tilde{c}_\alpha \Delta x)}.
\end{align*}
Consider the multisymplectic partitioned Runge--Kutta method (\ref{MPRK2first})-(\ref{MPRK2final}); if we multiply equation (\ref{pi0b}) by $\tilde{b}_\alpha$ and sum over $\alpha$, multiply equation (\ref{pi1b}) by $b_i$ and sum over $i$, and add the resulting equations together, we have
\begin{equation}\label{integralapproxDDW}
0 = \Delta x \sum_\alpha \tilde{b}_\alpha (\pi^0_{1[\alpha]} - \pi^1_{0[\alpha]}) + \Delta t \sum_i b_i (\pi^1_{[i]1} - \pi^1_{[i]0}) + \Delta t \Delta x \sum_{i,\alpha} b_i \tilde{b}_\alpha \frac{\partial H}{\partial\phi}(\Phi_{i\alpha},P^0_{i\alpha},P^1_{i\alpha}),
\end{equation}
where we used $X^{i\alpha} + Y^{i\alpha} = \partial H/\partial \phi (\Phi_{i\alpha},P^0_{i\alpha},P^1_{i\alpha}).$ Comparing these two, we see that the discrete method satisfies an approximation of the integral form of the De Donder--Weyl equation (\ref{integralDDW1}) and that the error in the approximation of the field equations is directly related to the quadrature error and the field and momenta expansions. Similar statements can be made about the other De Donder--Weyl equations, (\ref{integralDDW2}) and (\ref{integralDDW3}).

Now, let's write our approximation (\ref{integralapproxDDW}) of the integral De Donder--Weyl equations as a difference equation. For a quantity $f$ defined on the nodes of the edges $\{0\} \times [0, \Delta x]$ and $\{\Delta t\} \times [0,\Delta x]$ (and similarly a quantity $g$ defined on the nodes of the edges $[0,\Delta t] \times \{0\}$ and $[0,\Delta t] \times \{\Delta x\}$), define
\begin{align*}
\delta^0_{[\alpha]} f &= f_{1[\alpha]} - f_{0[\alpha]}, \\
\delta^1_{[i]} g &= g_{[i]1} - g_{[i]0}.
\end{align*}
Define the discrete difference operators
\begin{align*}
\partial^\Box_0 &= \frac{1}{\Delta t} \sum_\alpha \tilde{b}_\alpha \delta^0_{[\alpha]}, \\
\partial^\Box_1 &= \frac{1}{\Delta x} \sum_i b_i \delta^1_{[i]}.
\end{align*}
Dividing equation (\ref{integralapproxDDW}) by $\Delta t \Delta x$, we see that it satisfies
$$ \partial_0^\Box \pi^0 + \partial_1^\Box \pi^1 = -\sum_{i\alpha} b_i \tilde{b}_\alpha \frac{\partial H}{\partial\phi}(\Phi_{i\alpha},P^0_{i\alpha},P^1_{i\alpha}) \equiv -\Big\langle \frac{\partial H}{\partial \phi} \Big\rangle_{\Box}, $$
where $\langle \frac{\partial H}{\partial \phi} \rangle_{\Box}$ denotes our quadrature approximation of the average value of $\partial H/\partial\phi$ on $\Box$. Similarly, the other discrete equations satisfy
\begin{align*}
\partial_0^\Box \varphi &= \Big\langle \frac{\partial H}{\partial p^0} \Big\rangle_{\Box}, \\
\partial_1^\Box \varphi &= \Big\langle \frac{\partial H}{\partial p^1} \Big\rangle_{\Box}.
\end{align*}
These difference equations correspond to our discretization of (the integral form) of the DDW equations $\partial_0p^0 + \partial_1p^1 = -\partial H/\partial\phi$, $\partial_\mu\phi = \partial H/p^\mu$. As mentioned in Section \ref{Multisymplectic Integrators Subsection}, a method is called multisymplectic if the difference operators used in the discretization of the field equations are the same difference operators which appear in the discrete multisymplectic form formula that the method admits. In our case, if we divide the discrete multisymplectic form formula (\ref{Discrete MFF}) by $\Delta t \Delta x$, we see that it satisfies
$$ \partial_0^\Box \omega^0 + \partial_1^\Box \omega^1 = 0 $$
(when evaluated on discrete first variations), where $\omega^0 = d\varphi \wedge d\pi^0, \omega^1 = d\varphi \wedge d\pi^1$. Hence, our method is multisymplectic in the sense that the difference operators which appear in the difference equation that the discrete solution satisfies over $\Box \in \mathcal{T}(X)$ are the same difference operators which appear in the discrete multisymplectic form formula.

\section{Numerical Example}\label{Numerical Section}
For our numerical example, we will study the $(1+1)-$dimensional sine--Gordon equation,
\begin{equation}\label{sine-Gordon}
\partial_0^2 \phi(t,x) - \partial_1^2 \phi(t,x) = - \sin\phi(t,x).
\end{equation}
The Hamiltonian for this equation is given by
\begin{equation}\label{SG Hamiltonian}
H(\phi,p^0,p^1) = \frac{1}{2}(p^0)^2 - \frac{1}{2}(p^1)^2 - \cos\phi.
\end{equation}
The De Donder--Weyl equations corresponding to this Hamiltonian are
\begin{subequations}
\begin{align}
\partial_0\phi &= \partial H/\partial p^0 = p^0, \label{SG DDW a} \\
\partial_1\phi &= \partial H/\partial p^1 = -p^1, \label{SG DDW b} \\
\partial_0p^0 + \partial_1p^1 &= -\partial H/\partial \phi = -\sin\phi. \label{SG DDW c}
\end{align}
\end{subequations}
Note that substituting (\ref{SG DDW a}) and (\ref{SG DDW b}) into (\ref{SG DDW c}) recovers (\ref{sine-Gordon}).

With this example, we aim to qualitatively show the preservation of multisymplecticity by considering the family of soliton solutions,
\begin{equation}\label{soliton}
\phi_v(t,x) = 4 \arctan\left(\exp\left( \frac{x-vt}{\sqrt{1-v^2}} \right) \right),
\end{equation}
where the family of solutions is indexed by a parameter $v \in (0,1)$. Consider the following curve on the space of sections of the restricted dual jet bundle
$$ v \mapsto (\phi_v,p^0_v,p^1_v) \equiv (\phi_v, \partial_0\phi_v, -\partial_1\phi_v) $$
and note that it is differentiable for $v \in (0,1)$. Thus, the associated vector field, given by differentiating the above map with respect to $v$, defines a vector field on the space of sections of the restricted dual jet bundle. The associated vector field is a first variation on the space of soliton solutions, since its flow maps soliton solutions to other soliton solutions.

To visualize multisymplecticity for this example, we observe the following. Each soliton solution (\ref{soliton}) propagates to the right at speed $v$ in time, without changing form. Thus, the shape of a soliton solution in the $(\phi,p^0)$ plane does not change with respect to time. Hence, for a family of soliton solutions, the associated area in the $(\phi,p^0)$ plane will not expand or contract as the system evolves in time. In other words, restricting to the above first variations, this means that
$$ \partial_0\omega^0 = 0. $$
By the multisymplectic form formula $\partial_0\omega^0 + \partial_1\omega^1 = 0$, we also have that 
$$ \partial_1\omega^1 = 0, $$
and hence the family of soliton solutions, occupying an area in the $(\phi,p^1)$ plane, will not expand or contract as the system evolves in space. This example then provides an intuitive way to visualize multisymplecticity as symplecticity in each spacetime direction, since the multisymplectic conservation law splits into two symplectic conservation laws, for this given family of solutions. This is a multisymplectic analogue of the visualization of symplecticity in the literature for symplectic integrators, where one evolves a family of initial conditions occupying an area in phase space; for symplectic integrators, this area is preserved under the flow of the integrator, unlike a generic method (see, for example, \citet{HaLuWa2006}).

\textbf{Explicit Methods for Separable Hamiltonians.} Recall that in the above derivation of the multisymplectic partitioned Runge--Kutta method, we used that the Runge--Kutta matrices $a,\tilde{a}$ were invertible and hence the variational construction in Section \ref{MPRK subsection} does not directly apply to explicit Runge--Kutta matrices $a$ and $\tilde{a}$ (since explicit Runge--Kutta methods have strictly lower triangular Runge--Kutta matrices and hence are non-invertible). Since equation (\ref{sine-Gordon}) is nonlinear, using an implicit  method would be computationally expensive and hence an explicit method would be preferable. However, for separable Hamiltonians of the form $H(\phi,p^\mu) = K(p^\mu) + V(\phi)$ (as is the case for the sine--Gordon Hamiltonian (\ref{SG Hamiltonian})), we can derive an explicit method as follows. Let $a$ be an explicit Runge--Kutta matrix such that its symplectic pair $a^{(2)}$ is invertible, where again the symplectic pair is given by
$$ a^{(2)}_{ij} = \frac{b_jb_i - b_ja_{ji}}{b_i}. $$
Then, it follows that the symplectic pair of the symplectic pair of $a$ equals $a$, i.e., $(a^{(2)})^{(2)} = a$, since
\begin{align*}
(a^{(2)})^{(2)}_{ij} = \frac{b_jb_i - b_ja^{(2)}_{ji}}{b_i} = \frac{b_jb_i - b_j \frac{b_ib_j - b_ia_{ij}}{b_j}}{b_i} = a_{ij}.
\end{align*}
Thus, we can choose $a$ to instead be $a^{(2)}$, so that the symplectic pair of $a$ becomes an explicit Runge--Kutta matrix. The variational construction now applies with this choice of $a$, since it is invertible. For a separable Hamiltonian, the integration scheme splits and the system can be evolved explicitly in time. 

\textbf{Numerical Scheme.} We take a one-stage Runge--Kutta matrix in the temporal direction $a=1$ (with $b=1$, $c=1$) so that $a^{(2)} = 0$, and similarly in the spatial direction $\tilde{a}=1$ (with $\tilde{b}=1$, $\tilde{c}=1$) so that $\tilde{a}^{(2)} = 0$. Let $\Phi_{a,b},P^0_{a,b},P^1_{a,b},V_{a,b},W_{a,b},X_{a,b},Y_{a,b}$ denote the internal stages associated to $\Box_{ab} = \{t_a, t_a + \Delta t\}\times \{x_b,x_b+\Delta x\}$. Letting $\varphi_{a,b}$ denote the value of $\varphi$ at $\{t_a,x_b\}$ (and similarly for the momenta), the multisymplectic partitioned Runge--Kutta method (\ref{MPRK2first})-(\ref{MPRK2final}) gives, with the choice of the sine--Gordon Hamiltonian (\ref{SG Hamiltonian}),
\begin{subequations}
\begin{align}
P^0_{a,b} &= \pi^0_{a,b+1}, \\
\varphi_{a+1,b+1} &= \varphi_{a,b+1} + \Delta t V_{a,b} = \Phi_{a,b},\\
\pi^0_{a+1,b+1} &= \pi^0_{a,b+1} + \Delta t X_{a,b}, \\
P^1_{a,b} &= \pi^1_{a+1,b}, \label{SG scheme P1} \\
\varphi_{a+1,b+1} &= \varphi_{a+1,b} + \Delta x W_{a,b}, \label{SG scheme W0} \\
\pi^1_{a+1,b+1} &= \pi^1_{a+1,b} + \Delta x Y_{a,b}, \\
V_{a,b} &= \frac{\partial H}{\partial p^0}(\Phi_{a,b}, P^0_{a,b},P^1_{a,b}) = P^0_{a,b}, \label{SG scheme V}\\
W_{a,b} &= \frac{\partial H}{\partial p^1}(\Phi_{a,b}, P^0_{a,b},P^1_{a,b}) = -P^1_{a,b}, \label{SG scheme W}\\
X_{a,b} + Y_{a,b} &= \frac{\partial H}{\partial \phi}(\Phi_{a,b}, P^0_{a,b},P^1_{a,b}) = -\sin(\Phi_{a,b}). \label{SG Scheme XY}
\end{align}
\end{subequations}
Eliminating the internal stage variables, equations (\ref{SG scheme V}) and (\ref{SG Scheme XY}) can be expressed as an integration scheme in time
\begin{align*}
\frac{\varphi_{a+1,b} - \varphi_{a,b}}{\Delta t} &= \pi^0_{a,b}, \\
\frac{\pi^0_{a+1,b} - \pi^0_{a,b}}{\Delta t} + \frac{\pi^1_{a+1,b} - \pi^1_{a+1,b-1}}{\Delta x} &= -\sin(\varphi_{a+1,b})
\end{align*}
(where we shifted $b \mapsto b-1$). Further eliminating the $\pi^1$ variables using (\ref{SG scheme P1}), (\ref{SG scheme W0}), (\ref{SG scheme W}), the second equation above can be expressed as
$$ \frac{\pi^0_{a+1,b} - \pi^0_{a,b}}{\Delta t} - \frac{\varphi_{a+1,b+1} - 2\varphi_{a+1,b} + \varphi_{a+1,b-1}}{\Delta x^2} = -\sin(\varphi_{a+1,b}). $$
Thus, the corresponding numerical scheme is
\begin{subequations}
\begin{align}
\varphi_{a+1,b} &= \varphi_{a,b} + \Delta t\, \pi^0_{a,b}, \label{SG explicit scheme a} \\
\pi^0_{a+1,b} &= \pi^0_{a,b} + \Delta t \frac{\varphi_{a+1,b+1} - 2\varphi_{a+1,b} + \varphi_{a+1,b-1}}{\Delta x^2} - \Delta t\sin(\varphi_{a+1,b}). \label{SG explicit scheme b}
\end{align}
\end{subequations}
The scheme corresponds to discretizing the first-order formulation of the sine--Gordon equation,
\begin{align*}
\partial_0p^0 &= \partial_1^2\phi - \sin\phi,\\
\partial_0\phi &= p^0,
\end{align*}
in space using the standard discrete Laplacian and in time using the (adjoint) symplectic Euler method. We refer to the method (\ref{SG explicit scheme a})-(\ref{SG explicit scheme b}) as MSE (multisymplectic Euler).

As discussed in Section \ref{Discrete Hamiltonian Field Theory}, this scheme can be computed in a time marching fashion, given supplied initial conditions for $\varphi$ and $\pi^0$, as well as spatial boundary conditions. This scheme is explicit, since the values of the field can first be updated using (\ref{SG explicit scheme a}), followed by updating the temporal momenta using (\ref{SG explicit scheme b}). For the numerical experiment, we will compare this scheme to the scheme which uses the standard discrete Laplacian in space and the forward Euler method in time,
\begin{subequations}
\begin{align}
\varphi_{a+1,b} &= \varphi_{a,b} + \Delta t\, \pi^0_{a,b}, \label{SG fwd Euler scheme a} \\
\pi^0_{a+1,b} &= \pi^0_{a,b} + \Delta t \frac{\varphi_{a,b+1} - 2\varphi_{a,b} + \varphi_{a,b-1}}{\Delta x^2} - \Delta t\sin(\varphi_{a,b}). \label{SG fwd Euler scheme b}
\end{align}
\end{subequations}
We refer to the method (\ref{SG fwd Euler scheme a})-(\ref{SG fwd Euler scheme b}) as FE (forward Euler).

For our numerical experiment, we consider a family of initial conditions given by interpolating the soliton solutions $(\varphi_v(x),\pi^0_v(x)) =(\phi_v(0,x),\partial_0v(0,x))$ onto the spatial grid, for several values of $v$ ($v=0.50,0.47,0.45$), on a spatial domain $[-L,L]$. We choose Neumann boundary conditions $\pi^1(-L) = 0 = \pi^1(L)$ and choose $L$ sufficiently large so that the Neumann conditions are satisfied initially, up to a desired level of error (since $\pi^1(x) = -\partial_1\phi_v(0,x)$ decays monotonically to $0$ as $|x|$ increases), say $L = 20$ (so that $\pi^1(L) = \pi^1(-L) \sim 10^{-10}$). To demonstrate the robustness of MSE, we take a large spatial step $\Delta x = 0.1$ and a time step $\Delta t = \Delta x/2$; the experiment is run until a final time $T=20$.

\begin{figure}[h!]
\begin{center}
\includegraphics[width=100mm]{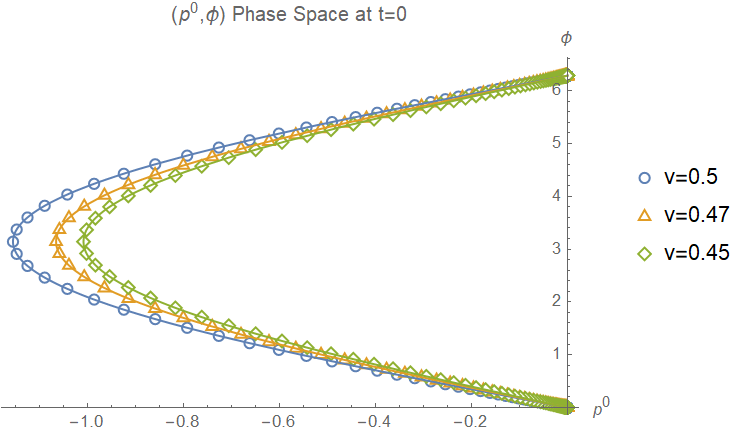}
\caption{The $(p^0,\phi)$ phase space distribution of the initial conditions (running over all spatial nodes in $[-L,L]$). The solid curves indicate the exact distribution.}
\label{initial}
\end{center}
\end{figure}

\begin{figure}[h!]
\begin{center}
\includegraphics[width=100mm]{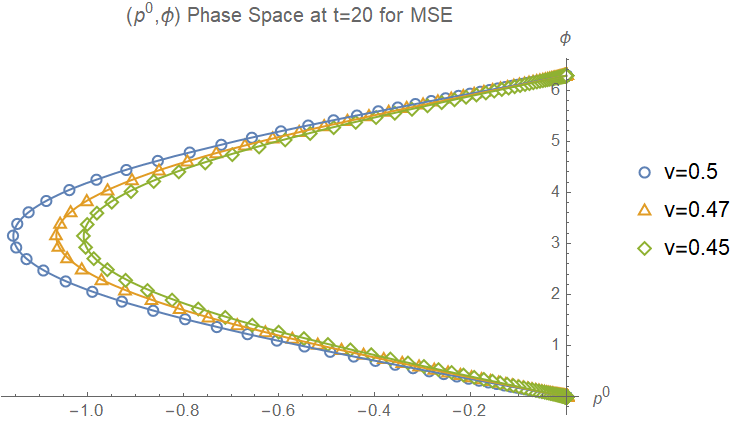}
\caption{The $(p^0,\phi)$ phase space distribution at $t=20$ using MSE (running over all spatial nodes in $[-L,L]$). The solid curves indicate the exact distribution.}
\label{mse20}
\end{center}
\end{figure}

\begin{figure}[h!]
\begin{center}
\includegraphics[width=100mm]{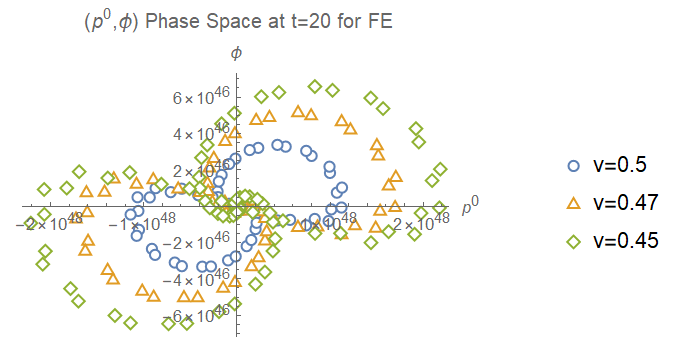}
\caption{The $(p^0,\phi)$ phase space distribution at $t=20$ using FE (running over all spatial nodes in $[-L,L]$).}
\label{fe20}
\end{center}
\end{figure}

The initial $(p^0,\phi)$ phase space distribution is shown in Figure \ref{initial}. The $(p^0,\phi)$ phase space distribution at $t=20$ is shown in Figures \ref{mse20} and \ref{fe20} for MSE and FE, respectively. Comparison of Figures \ref{initial} and \ref{mse20} shows the preservation of symplecticity in the temporal direction for the method MSE, whereas it is clearly not preserved for FE. 

\begin{figure}[h!]
\begin{center}
\includegraphics[width=100mm]{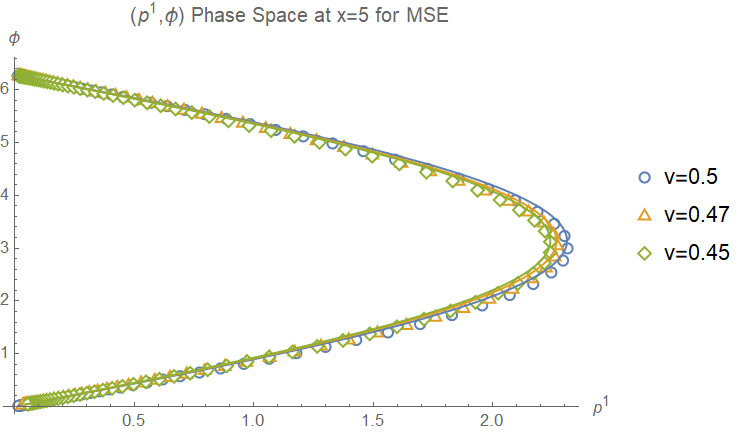}
\caption{The $(p^1,\phi)$ phase space distribution at $x=5$ using MSE (running over all timesteps in $[0,T]$). The solid curves indicate the exact distribution.}
\label{mse5}
\end{center}
\end{figure}

\begin{figure}[h!]
\begin{center}
\includegraphics[width=100mm]{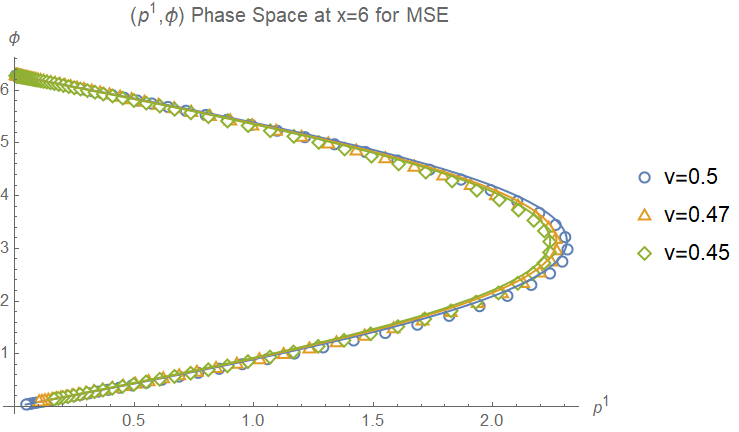}
\caption{The $(p^1,\phi)$ phase space distribution at $x=6$ using MSE (running over all timesteps in $[0,T]$). The solid curves indicate the exact distribution.}
\label{mse6}
\end{center}
\end{figure}

Similarly, Figures \ref{mse5} and \ref{mse6} shown the preservation of symplecticity in the $x$ direction for MSE.

\newpage
\section{Conclusion and Future Directions}
In this paper, we extended the construction of Hamiltonian variational integrators to the setting of multisymplectic Hamiltonian PDEs. Our construction is based on a discrete approximation of the boundary Hamiltonian, introduced in \citet{VaLiLe2011}. Through the Type II variational principle, this discrete boundary Hamiltonian is a generating function for the discrete Hamilton's equations that define our multisymplectic integrator. The discrete variational principle automatically yields integrators which are multisymplectic and satisfy a discrete Noether's theorem for group-invariant discretizations. As an application of this variational framework, we derived the class of multisymplectic partitioned Runge--Kutta methods; however, our construction is more general and is not limited to this class of multisymplectic integrators. Finally, we showed that the discrete multisymplecticity which arose from the discrete variational principle agrees with the notion of discrete multisymplecticity introduced in \citet{BrRe2001}.

Perhaps the most natural research direction is to establish a variational error analysis result which demonstrates that a computable discrete Hamiltonian that approximates the boundary Hamiltonian to a given order of accuracy will result in a numerical method for the Hamiltonian partial differential equation with the same order of accuracy. It should be observed that this poses two main challenges as compared to the case for ordinary differential equations. The first is that the boundary of the spacetime domain is in general curved, and the space of boundary data (and boundary momentum) is infinite-dimensional. As such, one would first have to approximate the spacetime domain with a spacetime mesh, and choose a finite-dimensional subspace for sections of the dual jet bundle that is subordinate to this spacetime mesh. Then, the error between the computable discrete Hamiltonian and the boundary Hamiltonian can be decomposed into three terms, the first of which can be bounded by assuming that the boundary-value problem is well-posed and therefore has continuous dependence on the boundary data, the second is associated with the variational crime of replacing the spacetime domain with a spacetime mesh, and the third is a term that is analogous to what arises in the usual variational error analysis for ordinary differential equations.

The second natural direction would be to establish a quasi-optimality result which demonstrates that the variational error in the construction of a Galerkin boundary Hamiltonian is bounded from above by a multiple of the best approximation error of the finite-dimensional function space used to approximate sections of the configuration bundle.

Finally, it was established in \citet{McAr2020} that many hybridizable discontinuous Galerkin methods are multisymplectic when applied to semilinear elliptic PDEs in mixed form, and it would be interesting to see the kind of multisymplectic Hamiltonian variational integrators that would arise for Hamiltonian time-evolution PDEs when using spacetime discontinuous Galerkin finite element spaces to discretize the dual jet bundle.

\section*{Acknowledgements}
The authors would like to thank the referees for their careful review of this paper and their helpful suggestions. BT was supported by the NSF Graduate Research Fellowship DGE-2038238, and by NSF under grants DMS-1411792, DMS-1813635. ML was supported by NSF under grants DMS-1411792, DMS-1345013, DMS-1813635, by AFOSR under grant FA9550-18-1-0288, and by the DoD under grant 13106725 (Newton Award for Transformative Ideas during the COVID-19 Pandemic). 

\appendix
\section{Higher Spacetime Dimensions}\label{HigherDimensions}
In this appendix, we treat the case of a spacetime tensor product (hyper)rectangular mesh in $(n+1)$-spacetime dimensions, where the coordinates on spacetime are given by $\{x^\mu\}_{\mu=0}^n$. Let $\mathcal{T}(X)$ be a regular (hyper)rectangular mesh, with $\Delta x^\mu$ the spacing in the $x^\mu$ direction. We index the nodes of this mesh by $x^\mu_a = a \Delta x^\mu$ (where $a$ is an integer) and consider $\Box_{a^0\dots a^n} \in \mathcal{T}(X)$ given by $\Box_{a^0\dots a^n} = \prod_{\mu=0}^n [x^\mu_{a^\mu}, x^\mu_{a^\mu} + \Delta x^\mu]$, where $\prod$ denotes the (ordered) Cartesian product. Fix a spacetime direction $x^\mu$. For this direction, there are two $(n-1)$-dimensional faces of $\Box_{a^0\dots a^n}$, located along the hyperplanes $x^\mu = x^\mu_{a^\mu}$ and $x^\mu = x^\mu_{a^\mu+1}$, to which the unit vector in the $x^\mu$ direction is normal. To each pair of such faces, we associate a quadrature rule (for simplicity, we consider the case of one quadrature point). The field values at this pair of quadrature points are denoted $\varphi_{[a^0]\dots [a^{\mu-1}]a^\mu[a^{\mu+1}]\dots [a^{n}]}$ and $\varphi_{[a^0]\dots [a^{\mu-1}]a^\mu +1[a^{\mu+1}]\dots [a^{n}]}$, where the unbracketed indices $a^\mu$ and $a^\mu +1$ indices denote the faces with smaller and larger $x^\mu$, respectively. Similarly, the corresponding normal momenta to these faces are denoted $\pi^\mu_{[a^0]\dots [a^{\mu-1}]a^\mu[a^{\mu+1}]\dots [a^{n}]}$ and $\pi^\mu_{[a^0]\dots [a^{\mu-1}]a^\mu +1[a^{\mu+1}]\dots [a^{n}]}$. Note that, in the $(1+1)$-dimensional case, this notation agrees with the notation that we used in Section \ref{Discrete Hamiltonian Field Theory} (where $a^0 = a, a^1 = b$). 

We take $B(\Box_{a^0\dots a^n})$ to consist of the ``forward" faces; that is, $B(\Box_{a^0\dots a^n})$ is the union, over all $\mu$, of the face in the $x^\mu$ direction with larger $x^\mu$ coordinate, $x^\mu = x^\mu_{a^\mu+1}$ (and similarly $A(\Box_{a^0\dots a^n})$ is the union, over all $\mu$ of the face in the $x^\mu$ direction with smaller $x^\mu$ coordinate, $x^\mu = x^\mu_{a^\mu}$). For brevity in the following equations, let
\begin{align*}
\pi^\mu_{]a^\mu[} &\equiv \pi^\mu_{[a^0]\dots [a^{\mu-1}]a^\mu[a^{\mu+1}]\dots [a^{n}]},\\
\varphi_{]a^\mu[} &\equiv \varphi_{[a^0]\dots [a^{\mu-1}]a^\mu[a^{\mu+1}]\dots [a^{n}]}, \\
\pi^\mu_{]a^\mu+1[} &\equiv \pi^\mu_{[a^0]\dots [a^{\mu-1}]a^\mu+1[a^{\mu+1}]\dots [a^{n}]},\\
\varphi_{]a^\mu+1[} &\equiv \varphi_{[a^0]\dots [a^{\mu-1}]a^\mu+1[a^{\mu+1}]\dots [a^{n}]},
\end{align*}
where we implicitly understand that $(a^0,\dots,a^n)$ are fixed. Then the quadrature approximation of the integral over $B$ is given by 
$$ \sumint_{B(\Box_{a^0\dots a^n})} \pi_B \varphi_B = \sum_{\mu=0}^{n} \left[ \pi^\mu_{]a^\mu +1[} \varphi_{]a^\mu +1[} \Delta^{n}x_\mu \right], $$
where $\Delta^{n}x_\mu \equiv \prod_{\nu \neq \mu}\Delta x^\nu$. Letting $\varphi_A$ denote the collection of values of $\varphi$ on the quadrature points on $A(\Box_{a^0\dots a^n})$ (and similarly for $\pi_B$), the associated discrete boundary Hamiltonian is 
\begin{align*}
H_d^+&(\varphi_A,\pi_B) = \ext\left( \sum_{\mu=0}^{n} \left[ \pi^\mu_{]a^\mu +1[} \varphi_{]a^\mu +1[} \Delta^nx_\mu \right] - S_d^{\Box_{a^0\dots a^n}}[\phi,p]\right),
\end{align*}
where $S_d^{\Box_{a^0\dots a^n}}$ is some discrete approximation of the action and the extremization is over all $(\phi,p)$ in the discrete approximating space satisfying the prescribed $(\varphi_A,\pi_B)$ boundary conditions. The Type II variational principle yields the discrete forward Hamilton's equations: for each $\mu$,
\begin{align*}
\pi^\mu_{]a^\mu[} &= \frac{1}{\Delta^nx_\mu} D_{\varphi,A,\mu} H_d^+(\varphi_A,\pi_B), \\
 \varphi_{]a^\mu +1[} &= \frac{1}{\Delta^nx_\mu} D_{\pi,B,\mu} H_d^+(\varphi_A,\pi_B),
\end{align*}
where $D_{\phi,A,\mu}$ denotes differentiation with respect to the value of $\varphi$ on the node on $A$ in the $\mu$ direction, i.e., $\partial/\partial \varphi_{]a^\mu[}$ (and similarly $D_{\pi,B,\mu} = \partial/\partial \pi^\mu_{]a^\mu+1[}$). 

Analogous results to the main body of the paper can be derived for the case of higher spacetime dimensions. For example, the multisymplectic conservation law $d^2H_d^+ = 0$ (when evaluated on first variations) gives
\begin{align*}
\sum_\mu \left( d\varphi_{]a^\mu +1[} \wedge d\pi^\mu_{]a^\mu+1[} - d\varphi_{]a^\mu[} \wedge d\pi^\mu_{]a^\mu[} \right) \Delta^nx_\mu = 0
\end{align*}
(which formally is the quadrature approximation to $\int_{\Box_{a^0\dots a^n}} \omega^\mu(\cdot,\cdot)d^nx_\mu = 0$).

Similarly, the generalization to multiple quadrature points is straight-forward; for each pair of forward and backward $(n-1)$-dimensional faces in the $\mu$ direction, we can choose multiple quadrature points and weights on the faces (the quadrature rules on the forward and backward faces in the same direction must be the same, but the quadrature rules can differ among the spacetime directions, as was the case in $(1+1)$-spacetime dimensions). Associated to these quadrature points are the field and normal momenta values. Then, the discrete forward Hamilton's equations just states that the value of $\varphi$ on a quadrature node in $B$ is given by differentiating $H_d^+$ with respect to the normal momenta $\pi$ on that node, divided by the product of $\Delta^nx_\mu$ and the quadrature weight at that node (and similarly, the value of the normal momenta $\pi$ on a quadrature node in $A$ is given by differentiating $H_d^+$ with respect to the field value on that node, divided by the product of $\Delta^nx_\mu$ and the quadrature weight at that node). As one can verify, in the $(1+1)$-dimensional case, this precisely reproduces (\ref{Discrete Hamilton's Equations, full quad, 1})-(\ref{Discrete Hamilton's Equations, full quad, 2}).

One can then proceed as we did in the main body of the paper, in using the Galerkin construction as a discrete approximation for the action. Utilizing analogous expansions to those in the main body of the paper (with an expansion in each spacetime direction), the resulting variational integrator would then give a multisymplectic partitioned Runge--Kutta method, where the integrator would formally be a symplectic partitioned Runge--Kutta method in each spacetime direction with the internal stages satisfying the De Donder--Weyl equations. 

Finally, it is worth noting that, at the start of Section \ref{Discrete Hamiltonian Field Theory}, we laid a general formulation for unstructured meshes, arbitrary finite element spaces, and arbitrary spacetime dimensions. However, in general, the form of the discrete Hamilton's equations arising from the Type II variational principle can not be written explicitly, which is why we specialized to the case of spacetime tensor product (hyper)rectangular meshes. It would be interesting to determine the form of the discrete forward Hamilton's equations in other settings for particular choices of meshes, finite element spaces, and spacetime dimensions. For example, although a fully unstructured spacetime mesh would be challenging, one could consider a spacetime tensor product mesh which is the tensor product of an unstructured spatial mesh and a regular temporal mesh. Even in the case of a (hyper)rectangular mesh, it would be interesting to consider finite element spaces of differential forms (such as the $Q_r^{-}\Lambda^k$ spaces arising in finite element exterior calculus~\cite{ArFaWi2006}), which could be interesting in physical applications such as lattice field theory. 

\section{Relation to Galerkin Lagrangian Variational Integrators}\label{LagrangianGalerkin}
In this appendix, we discuss the relation between Galerkin Hamiltonian and Lagrangian variational integrators. From the Lagrangian perspective, the appropriate generating functional is the boundary Lagrangian (see \citet{VaLiLe2011}),
$$ L_{\partial U}(\varphi) = \int_U L(\phi,\partial_\mu\phi)d^{n+1}x, $$
where the expression on the right hand side is extremized over all $\phi$ such that $\phi|_{\partial U} = \varphi$. In general, methods derived from discretizing the boundary Hamiltonian and boundary Lagrangian are not expected to be equivalent, even in the hyperregular case, as was shown in \citet{ScLe2017} for the case of mechanics (where the boundary Hamiltonian and boundary Lagrangian are referred to as the exact discrete Hamiltonian and the exact discrete Lagrangian, respectively).

However, for the case of Galerkin Lagrangian variational integrators on a $(1+1)-$dimensional rectangular mesh, they are equivalent (in the hyperregular case), for a suitable choice of discrete boundary Lagrangian and using the same (Galerkin based) expansions that we utilized in our construction of Galerkin Hamiltonian variational integrators. We assume the same field expansions that we used in Section \ref{MPRK subsection} (when discussing independent internal stages). Unlike the boundary Hamiltonian where $\varphi$ is specified on $A$ and $\pi$ is specified on $B$, the Lagrangian perspective specifies $\varphi$ on both $A$ and $B$. One can define a discrete boundary Lagrangian as
\begin{align*}
L_d^{\partial\Box}(\varphi_A,\varphi_B) &= \ext_{\stackrel{\phi}{\phi|_A = \varphi_A, \phi|_B = \varphi_B}} \Delta t \Delta x \sum_{i\alpha} L\big(\phi(c_i\Delta t, \tilde{c}_\alpha \Delta x), \partial_0\phi(c_i\Delta t, \tilde{c}_\alpha \Delta x), \partial_1\phi(c_i\Delta t, \tilde{c}_\alpha \Delta x) \big) \\
&= \ext_{V^{i\alpha},W^{i\alpha},\lambda_\alpha,\lambda_i, \lambda_{i\alpha}} \Delta t \Delta x \Biggl[ \sum_{i,\alpha}  b_i\tilde{b}_\alpha L(\Phi^\theta_{i\alpha}, V^{i\alpha}, W^{i\alpha}) \\ 
&\qquad\qquad\qquad\qquad\qquad\qquad\quad + \sum_\alpha \lambda_\alpha\Biggl(\varphi_{1[\alpha]} - \varphi_{0[\alpha]} - \Delta t \sum_j b_jV^{j\alpha} \Biggr) \\
&\qquad\qquad\qquad\qquad\qquad\qquad\quad + \sum_i \lambda_i\Biggl(\varphi_{[i]1} - \varphi_{[i]0} - \Delta x \sum_\beta \tilde{b}_\beta W^{i\beta}\Biggr) \\
&\qquad\qquad\qquad\qquad\qquad\qquad\quad + \sum_{i,\alpha} \lambda_{i\alpha} (\Phi_{i\alpha} + \tilde{\Phi}_{i\alpha})\Biggr],
\end{align*}
where in the first line, the right hand side is extremized over the finite-dimensional function space chosen in the Galerkin construction (to obtain a discrete boundary Lagrangian instead of the exact discrete boundary Lagrangian which extremizes over an infinite-dimensional space). The second equality follows from substituting the chosen expansion and explicitly enforcing that the boundary condition $\phi|_B = \varphi_B$ are satisfied by the Lagrange multipliers $\lambda_\alpha$ and $\lambda_i$. The normal momenta are then obtained by enforcing the variational principle, which gives the normal momenta $\pi_A,\pi_B$ in terms of the derivatives of $L_d^{\partial \Box}$ with respect to $\varphi_A,\varphi_B$. This defines a Galerkin Lagrangian variational integrator.

\begin{prop}
If the continuous Hamiltonian $H$ is hyperregular and the associated Lagrangian $L$ is constructed by the Legendre transform, then the Galerkin Hamiltonian variational integrator and the Galerkin Lagrangian variational integrator are equivalent, for the same choice of expansion (i.e., specified by the basis functions $\psi_i$, $\tilde{\psi}_\beta$ and quadrature rules). 
\begin{proof}
The proof follows from using the Legendre transform to express 
\begin{align*}
\partial_\mu\phi = \frac{\partial H(\phi,p^0,p^1)}{\partial p^\mu},
\end{align*}
which is invertible by assumption of hyperregularity (i.e., one can express the momenta in terms of the field and their derivatives). The computation then follows analogously to the $1$-dimensional (mechanics) case, as shown in \citet{LeZh2009}, noting that the Legendre transform holds at the internal stages. 
\end{proof}
\end{prop}

It is expected that this equivalence holds in the case of higher-dimensional spacetime tensor product (hyper)rectangular meshes, although it is still unclear to what degree this holds for general unstructured spacetime meshes and general finite element spaces. We aim to explore this in future work. 

%\nocite{*}

\bibliographystyle{plainnat}
\bibliography{hamiltonianvarint.bib}

\begin{thebibliography}{44}
\providecommand{\natexlab}[1]{#1}
\providecommand{\url}[1]{\texttt{#1}}
\expandafter\ifx\csname urlstyle\endcsname\relax
  \providecommand{\doi}[1]{doi: #1}\else
  \providecommand{\doi}{doi: \begingroup \urlstyle{rm}\Url}\fi

\bibitem[Abraham and Marsden(1978)]{AbMa1978}
R.~Abraham and {J.~E.} Marsden.
\newblock \emph{Foundations of mechanics}.
\newblock Benjamin/Cummings Publishing Co. Inc. Advanced Book Program, Reading,
  Mass., 1978.
\newblock Second edition, revised and enlarged, With the assistance of Tudor
  Ra{\c{t}}iu and Richard Cushman.

\bibitem[Arnold et~al.(2006)Arnold, Falk, and Winther]{ArFaWi2006}
{D. N.} Arnold, {R. S.} Falk, and R.~Winther.
\newblock Finite element exterior calculus, homological techniques, and
  applications.
\newblock \emph{Acta Numer.}, 15:\penalty0 1--155, 2006.

\bibitem[Arnold et~al.(2010)Arnold, Falk, and Winther]{ArFaWi2010}
{D. N.} Arnold, {R. S.} Falk, and R.~Winther.
\newblock Finite element exterior calculus: from {H}odge theory to numerical
  stability.
\newblock \emph{Bull. Amer. Math. Soc}, 47\penalty0 (2):\penalty0 281--354,
  2010.

\bibitem[Benettin and Giorgilli(1994)]{BeGi1994}
G.~Benettin and A.~Giorgilli.
\newblock On the {H}amiltonian interpolation of near to the identity symplectic
  mappings with application to symplectic integration algorithms.
\newblock \emph{J. Stat. Phys.}, 74:\penalty0 1117--1143, 1994.

\bibitem[Bridges(1997)]{Br1997}
{T. J.} Bridges.
\newblock Multi-symplectic structures and wave propagation.
\newblock \emph{Math. Proc. Cambridge Philos. Soc.}, 121\penalty0 (1):\penalty0
  147--190, 1997.

\bibitem[Bridges and Reich(2001{\natexlab{a}})]{BrRe2001}
{T. J.} Bridges and S.~Reich.
\newblock Multi-symplectic integrators: numerical schemes for {H}amiltonian
  {PDE}s that conserve symplecticity.
\newblock \emph{Physics Letters A}, 284\penalty0 (4):\penalty0 184 -- 193,
  2001{\natexlab{a}}.

\bibitem[Bridges and Reich(2001{\natexlab{b}})]{BrRe2001b}
T.~J. Bridges and S.~Reich.
\newblock Multi-symplectic spectral discretizations for the
  {Z}akharov--{K}uznetsov and shallow water equations.
\newblock \emph{Physica D: Nonlinear Phenomena}, 152-153:\penalty0 491 -- 504,
  2001{\natexlab{b}}.
\newblock Advances in Nonlinear Mathematics and Science: A Special Issue to
  Honor Vladimir Zakharov.

\bibitem[Bridges and Reich(2006)]{BrRe2006}
{T. J.} Bridges and S.~Reich.
\newblock Numerical methods for {H}amiltonian {PDEs}.
\newblock \emph{Journal of Physics A: Mathematical and General}, 39\penalty0
  (19):\penalty0 5287--5320, 2006.

\bibitem[Chen(2005)]{Ch2005}
J.~Chen.
\newblock Variational formulation for the multisymplectic {H}amiltonian
  systems.
\newblock \emph{Lett. Math. Phys.}, 71\penalty0 (3):\penalty0 243--253, 2005.

\bibitem[Chen(2008)]{Ch2008}
{J.-B.} Chen.
\newblock Variational integrators and the finite element method.
\newblock \emph{Applied Mathematics and Computation}, 196\penalty0
  (2):\penalty0 941--958, 2008.

\bibitem[Demoures et~al.(2016)Demoures, {Gay-Balmaz}, and Ratiu]{DeGBRa2016}
F.~Demoures, F.~{Gay-Balmaz}, and {T. S.} Ratiu.
\newblock Multisymplectic variational integrators for nonsmooth {L}agrangian
  cntinuum mechanics.
\newblock \emph{Forum of Mathematics, Sigma}, 4:\penalty0 e19, 2016.

\bibitem[Duruisseaux et~al.(2020)Duruisseaux, Schmitt, and Leok]{DuScLe2020}
V.~Duruisseaux, J.~Schmitt, and M.~Leok.
\newblock Adaptive {H}amiltonian variational integrators and symplectic
  accelerated optimization.
\newblock (preprint, \href{https://arxiv.org/abs/1709.01975}{\tt
  arxiv/1709.01975 [math.NA]}), 2020.

\bibitem[Gotay et~al.(1998)Gotay, Isenberg, Marsden, and
  Montgomery]{GoIsMaMo1998}
{M. J.} Gotay, J.~Isenberg, {J. E.} Marsden, and R.~Montgomery.
\newblock Momentum maps and classical relativistic fields. {P}art {I}:
  {C}ovariant field theory.
\newblock (preprint, \href{http://arxiv.org/abs/physics/9801019}{\tt
  arXiv:physics/9801019 [math-ph]}), 1998.

\bibitem[Gotay et~al.(2004)Gotay, Isenberg, Marsden, and
  Montgomery]{GoIsMaMo2004}
{M. J.} Gotay, J.~Isenberg, {J. E.} Marsden, and R.~Montgomery.
\newblock Momentum maps and classical relativistic fields. {P}art {II}:
  {C}anonical analysis of field theories.
\newblock (preprint, \href{http://arxiv.org/abs/physics/0411032}{\tt
  arXiv:physics/0411032[math-ph]}), 2004.

\bibitem[Hairer(1997)]{Ha1997}
E.~Hairer.
\newblock Variable time step integration with symplectic methods.
\newblock \emph{Applied Numerical Mathematics}, 25\penalty0 (2-3):\penalty0
  219--227, 1997.

\bibitem[Hairer et~al.(2006)Hairer, Lubich, and Wanner]{HaLuWa2006}
E.~Hairer, C.~Lubich, and G.~Wanner.
\newblock \emph{Geometric Numerical Integration: Structure-preserving
  algorithms for ordinary differential equations}, volume~31 of \emph{Springer
  Series in Computational Mathematics}.
\newblock Springer-Verlag, Berlin, second edition, 2006.

\bibitem[Hall and Leok(2015)]{HaLe2015}
J.~Hall and M.~Leok.
\newblock Spectral variational integrators.
\newblock \emph{Numer. Math.}, 130\penalty0 (4):\penalty0 681--740, 2015.

\bibitem[Hiptmair(2002)]{Hi2002}
R.~Hiptmair.
\newblock Finite elements in computational electromagnetism.
\newblock \emph{Acta Numerica}, 11:\penalty0 237--339, 2002.

\bibitem[Hong et~al.(2006)Hong, Liu, and Sun]{HoLiSu2006}
J.~Hong, H.~Liu, and G.~Sun.
\newblock The multi-symplecticity of partitioned {R}unge-{K}utta methods for
  {H}amiltonian {PDE}s.
\newblock \emph{Math. Comp.}, 75:\penalty0 167--181, 2006.

\bibitem[Islas and Schober(2004)]{IsSc2004}
A.L. Islas and C.M. Schober.
\newblock On the preservation of phase space structure under multisymplectic
  discretization.
\newblock \emph{J. Comp. Phys.}, 197\penalty0 (2):\penalty0 585 -- 609, 2004.

\bibitem[Islas and Schober(2006)]{IsSc2006}
A.L. Islas and C.M. Schober.
\newblock Conservation properties of multisymplectic integrators.
\newblock \emph{Future Generation Computer Systems}, 22\penalty0 (4):\penalty0
  412 -- 422, 2006.

\bibitem[Lall and West(2006)]{LaWe2006}
S.~Lall and M.~West.
\newblock Discrete variational {H}amiltonian mechanics.
\newblock \emph{Journal of Physics A: Mathematical and General}, 39\penalty0
  (19):\penalty0 5509--5519, 2006.

\bibitem[Leitz et~al.(2021)Leitz, {Sato Martin de Aimagro}, and
  Leyendecker]{LeSaLe2021}
T.~Leitz, {R. T.} {Sato Martin de Aimagro}, and S.~Leyendecker.
\newblock Multisymplectic {G}alerkin {L}ie group variational integrators for
  geometrically exact beam dynamics based on unit dual quaternion
  interpolation.
\newblock \emph{Computer Methods in Applied Mechanics and Engineering},
  374:\penalty0 113475, 2021.

\bibitem[Leok(2019)]{Le2019}
M.~Leok.
\newblock Variational discretizations of gauge field theories using
  group-equivariant interpolation.
\newblock \emph{Foundations of Computational Mathematics}, 19\penalty0
  (5):\penalty0 965--989, 2019.

\bibitem[Leok and Ohsawa(2011)]{LeOh2008}
M.~Leok and T.~Ohsawa.
\newblock Variational and geometric structures of discrete {D}irac mechanics.
\newblock \emph{Found. Comput. Math.}, 11\penalty0 (5):\penalty0 529--562,
  2011.

\bibitem[Leok and Zhang(2011)]{LeZh2009}
M.~Leok and J.~Zhang.
\newblock Discrete {H}amiltonian variational integrators.
\newblock \emph{IMA J. Numer. Anal.}, 31\penalty0 (4):\penalty0 1497--1532,
  2011.

\bibitem[Le{\'o}n et~al.(2017)Le{\'o}n, Prieto-Mart\'{i}nez, Rom\'{a}n-Roy, and
  Vilari{\~n}o]{LePrRoVi2017}
M.~Le{\'o}n, P.~D. Prieto-Mart\'{i}nez, N.~Rom\'{a}n-Roy, and S.~Vilari{\~n}o.
\newblock {H}amilton-{J}acobi theory in multisymplectic classical field
  theories.
\newblock \emph{J. Math. Phys.}, 58:\penalty0 092901, 36 pp., 2017.

\bibitem[Lew et~al.(2003)Lew, Marsden, Ortiz, and West]{LeMaOrWe2003}
A.~Lew, {J. E.} Marsden, M.~Ortiz, and M.~West.
\newblock Asynchronous variational integrators.
\newblock \emph{Archive for Rational Mechanics and Analysis}, 167\penalty0
  (2):\penalty0 85--146, 2003.

\bibitem[Marsden and Shkoller(1999)]{MaSh1999}
{J. E.} Marsden and S.~Shkoller.
\newblock Multisymplectic geometry, covariant {H}amiltonians, and water waves.
\newblock \emph{Mathematical Proceedings of the Cambridge Philosophical
  Society}, 125\penalty0 (3):\penalty0 553--575, 1999.

\bibitem[Marsden and West(2001)]{MaWe2001}
{J. E.} Marsden and M.~West.
\newblock Discrete mechanics and variational integrators.
\newblock \emph{Acta Numer.}, 10:\penalty0 317--514, 2001.

\bibitem[Marsden et~al.(1998)Marsden, Patrick, and Shkoller]{MaPaSh1998}
{J. E.} Marsden, {G. W.} Patrick, and S.~Shkoller.
\newblock Multisymplectic geometry, variational integrators, and nonlinear
  {PDE}s.
\newblock \emph{Commun. Math. Phys.}, 199\penalty0 (2):\penalty0 351--395,
  1998.

\bibitem[Marsden et~al.(2001)Marsden, Pekarsky, Shkoller, and
  West]{MaPeShWe2001}
{J. E.} Marsden, S.~Pekarsky, S.~Shkoller, and M.~West.
\newblock Variational methods, multisymplectic geometry and continuum
  mechanics.
\newblock \emph{J. Geom. Phys.}, 38\penalty0 (3-4):\penalty0 253--284, 2001.

\bibitem[McLachlan and Stern(2020)]{McAr2020}
{R. I.} McLachlan and A.~Stern.
\newblock Multisymplecticity of hybridizable discontinuous {G}alerkin methods.
\newblock \emph{Foundations of Computational Mathematics}, 20\penalty0
  (1):\penalty0 35---69, 2020.

\bibitem[Ober-Bl\"{o}baum(2016)]{Ob2016}
S.~Ober-Bl\"{o}baum.
\newblock {Galerkin variational integrators and modified symplectic
  Runge--Kutta methods}.
\newblock \emph{IMA Journal of Numerical Analysis}, 37\penalty0 (1):\penalty0
  375--406, 02 2016.

\bibitem[Ohsawa et~al.(2010)Ohsawa, Bloch, and Leok]{OhBlLe2010}
T.~Ohsawa, {A. M.} Bloch, and M.~Leok.
\newblock Discrete {H}amilton--{J}acobi theory and discrete optimal control.
\newblock \emph{Proc. IEEE Conf. on Decision and Control}, pages 5438--5443,
  2010.

\bibitem[Reich(1999)]{Re1999}
S.~Reich.
\newblock Backward error analysis for numerical integrators.
\newblock \emph{SIAM J. Numer. Anal.}, 36:\penalty0 1549--1570, 1999.

\bibitem[Reich(2000{\natexlab{a}})]{Re2000}
S.~Reich.
\newblock Multi-symplectic {R}unge--{K}utta collocation methods for
  {H}amiltonian wave equations.
\newblock \emph{J. Comp. Phys.}, 157\penalty0 (2):\penalty0 473 -- 499,
  2000{\natexlab{a}}.

\bibitem[Reich(2000{\natexlab{b}})]{Re2000b}
S.~Reich.
\newblock Finite volume methods for multi-symplectic {PDE}s.
\newblock \emph{BIT Numerical Mathematics}, 40\penalty0 (3):\penalty0 559--582,
  2000{\natexlab{b}}.

\bibitem[Ryland et~al.(2007)Ryland, Mclachlan, and Frank]{RyMcFr2007}
B.~N. Ryland, R.~I. Mclachlan, and J.~Frank.
\newblock On the multisymplecticity of partitioned {R}unge--{K}utta and
  splitting methods.
\newblock \emph{International Journal of Computer Mathematics}, 84\penalty0
  (6):\penalty0 847--869, 2007.

\bibitem[Schmitt and Leok(2017)]{ScLe2017}
J.~M. Schmitt and M.~Leok.
\newblock Properties of {H}amiltonian variational integrators.
\newblock \emph{IMA Journal of Numerical Analysis}, 36\penalty0 (2), 2017.

\bibitem[Schmitt et~al.(2018)Schmitt, Shingel, and Leok]{ScShLe2017}
J.~M. Schmitt, T.~Shingel, and M.~Leok.
\newblock {L}agrangian and {H}amiltonian {T}aylor variational integrators.
\newblock \emph{{BIT} Numerical Mathematics}, 58\penalty0 (2):\penalty0
  457--488, 2018.

\bibitem[Tang(1994)]{Ta1994}
Y.-F. Tang.
\newblock Formal energy of a symplectic scheme for {H}amiltonian systems and
  its applications (i).
\newblock \emph{Computers \& Mathematics with Applications}, 27\penalty0
  (7):\penalty0 31 -- 39, 1994.

\bibitem[Vankerschaver et~al.(2012)Vankerschaver, Yoshimura, and
  Leok]{VaYoLeMa2012}
J.~Vankerschaver, H.~Yoshimura, and M.~Leok.
\newblock The {H}amilton--{P}ontryagin principle and multi-{D}irac structures
  for classical field theories.
\newblock \emph{J. Math. Phys.}, 53\penalty0 (7):\penalty0 072903 (25 pages),
  2012.

\bibitem[Vankerschaver et~al.(2013)Vankerschaver, Liao, and Leok]{VaLiLe2011}
J.~Vankerschaver, C.~Liao, and M.~Leok.
\newblock Generating functionals and {L}agrangian partial differential
  equations.
\newblock \emph{J. Math. Phys.}, 54\penalty0 (8):\penalty0 082901 (22 pages),
  2013.

\end{thebibliography}

\end{document}